\documentclass[11pt]{article}
\usepackage{amsfonts,amsmath,amssymb,accents,amsthm,bm}
\usepackage{verbatim}
\usepackage{hyperref}
\usepackage{graphicx}
\usepackage{subcaption}
\usepackage{xcolor}
\usepackage[normalem]{ulem}

\usepackage{tikz}
\usetikzlibrary{calc,arrows}

\newcommand{\file}{}

\newcommand{\detail}[1]{\par\noi{\bf [Proof detail\ }{#1}
\hfill{\bf ]}\par\noi\hspace{-4pt}}
\renewcommand{\detail}[1]{}

\newcommand{\dis}{\displaystyle}

\newcommand{\noi}{\noindent}
\newcommand{\halmos}{\rule{1ex}{1.4ex}}
\newcommand{\QED}{\nopagebreak{\hspace*{\fill}$\halmos$\medskip}}
\newcommand{\ruit}{\nopagebreak{\hspace*{\fill}\tikz[scale=0.18,baseline=-5pt]{\draw (0,-1) -- (1,0) -- (0,1) -- (-1,0) -- cycle;}\medskip}}
\newcommand{\med}{\medskip}
\newcommand{\quand}{\quad\mbox{and}\quad}

\newtheoremstyle{mythm}
  {}
  {}
  {\itshape}
  {}
  {\bfseries}
  {}
  {.5em}
  {#1 #2 \thmnote{(#3)}}

\theoremstyle{mythm}
\newtheorem{theorem}{Theorem}
\newtheorem{proposition}[theorem]{Proposition}
\newtheorem{lemma}[theorem]{Lemma}
\newtheorem{exercise}[theorem]{Exercise}
\newtheorem{corollary}[theorem]{Corollary}
\newtheorem{conjecture}[theorem]{Conjecture}
\newtheorem{example}[theorem]{Example}
\newtheorem{counterex}[theorem]{Counterexample}

\newcommand{\bt}{\begin{theorem}}
\newcommand{\et}{\end{theorem}}
\newcommand{\bl}{\begin{lemma}}
\newcommand{\el}{\end{lemma}}
\newcommand{\bp}{\begin{proposition}}
\newcommand{\ep}{\end{proposition}}
\newcommand{\bcor}{\begin{corollary}}
\newcommand{\ecor}{\end{corollary}}
\newcommand{\br}{\begin{remark}\rm}
\newcommand{\er}{\end{remark}}
\newcommand{\bcon}{\begin{conjecture}}
\newcommand{\econ}{\end{conjecture}}
\newcommand{\bex}{\begin{exercise}}
\newcommand{\eex}{\end{exercise}}
\newcommand{\bcou}{\begin{counterex}}
\newcommand{\ecou}{\end{counterex}}
\newcommand{\bexa}{\begin{example}}
\newcommand{\eexa}{\end{example}}
\newenvironment{remark}{\noi\textbf{Remark}}{\ruit}
\newcommand{\brm}{\begin{remark}}
\newcommand{\erm}{\end{remark}}

\newenvironment{Proof}[1][]{\noi\textbf{Proof #1}}{\QED}
\newcommand{\bpro}{\begin{Proof}}
\newcommand{\epro}{\end{Proof}}

\newcommand{\be}{\begin{equation}}
\newcommand{\ee}{\end{equation}}
\newcommand{\ba}{\begin{array}}
\newcommand{\ea}{\end{array}}
\newcommand{\bc}{\be\begin{array}{r@{\,}c@{\,}l}}
\newcommand{\ec}{\end{array}\ee}
\newcommand{\bac}{\begin{array}{r@{\,}c@{\,}l}}

\newcommand{\al}{\alpha}
\newcommand{\bet}{\beta}
\newcommand{\ga}{\gamma}

\newcommand{\de}{\delta}

\newcommand{\eps}{\varepsilon}
\newcommand{\la}{\lambda}

\newcommand{\sig}{\sigma}
\newcommand{\Sig}{\Sigma}

\newcommand{\om}{\omega}
\newcommand{\Om}{\Omega}

\newcommand{\si}{\ensuremath{\sigma}}

\newcommand{\Ai}{{\cal A}}

\newcommand{\Ci}{{\cal C}}
\newcommand{\Di}{{\cal D}}

\newcommand{\Fi}{{\cal F}}
\newcommand{\Gi}{{\cal G}}
\newcommand{\Hi}{{\cal H}}

\newcommand{\Li}{{\cal L}}
\newcommand{\Mi}{{\cal M}}

\newcommand{\Pc}{{\cal P}}

\newcommand{\Ti}{{\cal T}}

\newcommand{\E}{{\mathbb E}}

\newcommand{\N}{{\mathbb N}}
\newcommand{\Op}{{\mathbb O}}
\renewcommand{\P}{{\mathbb P}}

\newcommand{\R}{{\mathbb R}}
\let\Parag\S
\let\S\relax
\newcommand{\S}{{\mathbb S}}
\newcommand{\T}{{\mathbb T}}
\newcommand{\U}{{\mathbb U}}
\newcommand{\V}{{\mathbb V}}

\newcommand{\li}{\langle}
\newcommand{\re}{\rangle}

\newcommand{\desd}{\ensuremath{\Leftrightarrow}}
\newcommand{\volgt}{\ensuremath{\Rightarrow}}
\newcommand{\up}{\uparrow}

\newcommand{\sub}{\subset}
\newcommand{\beh}{\backslash}

\newcommand{\isd}{\stackrel{\scriptstyle\rm d}{=}}
\newcommand{\asto}[1]{\underset{{#1}\to\infty}{\longrightarrow}}

\newcommand{\Asto}[1]{\underset{{#1}\to\infty}{\Longrightarrow}}

\newcommand{\ti}{\tilde}
\newcommand{\dgg}{\dagger}
\newcommand{\ov}{\overline}
\newcommand{\un}{\underline}

\newcommand{\lvec}[1]{\accentset{\leftarrow}{#1}}

\newcommand{\pa}{\partial}
\newcommand{\ffrac}[2]{{\textstyle\frac{{#1}}{{#2}}}}
\newcommand{\dif}[1]{\ffrac{\partial}{\partial{#1}}}
\newcommand{\diff}[1]{\ffrac{\partial^2}{{\partial{#1}}^2}}

\newcommand{\di}{\mathrm{d}}
\newcommand{\half}{{[0,\infty)}}
\newcommand{\expo}{\mbox{\large\it e}}
\newcommand{\ex}[1]{\expo^{\,\textstyle{#1}}}

\newcommand{\var}{{\rm Var}}

\newcommand{\nab}{\nabla}

\newcommand{\ha}{\ffrac{1}{2}}

\newcommand{\Xb}{{\mathbf X}}

\newcommand{\qb}{{\mathbf q}}
\newcommand{\rb}{{\mathbf r}}

\newcommand{\cob}{{\mbox{\tt cob}}}
\newcommand{\dth}{{\mbox{\tt dth}}}
\newcommand{\bth}{{\mbox{\tt bth}}}
\newcommand{\bra}{{\mbox{\tt bra}}}

\newcommand{\ibf}{\mathbf{i}}
\newcommand{\jbf}{\mathbf{j}}

\newcommand{\lib}{{\lvec{\mathbf i}}}

\newcommand{\ch}{\check}

\newcommand{\poa}{\pa}

\newcommand{\omb}{\bm{\omega}}

\newcommand{\wurz}{\varnothing}

\begin{document}

\makeatletter\@addtoreset{equation}{section}
\makeatother\def\theequation{\thesection.\arabic{equation}} 

\renewcommand{\labelenumi}{{\rm (\roman{enumi})}}

\title{Recursive tree processes and the\\ mean-field limit of stochastic flows}
\author{
Tibor~Mach
\footnote{The Czech Academy of Sciences,
Institute of Information Theory and Automation,
Pod vod\'arenskou v\v e\v z\' i 4,
18208 Praha 8,
Czech Republic;
swart@utia.cas.cz}
\and
Anja~Sturm
\footnote{Institute for Mathematical Stochastics,
Georg-August-Universit\"at G\"ottingen,
Goldschmidtstr.~7,
37077 G\"ottingen,
Germany;
asturm@math.uni-goettingen.de}
\and
Jan~M.~Swart$\;^\ast$
}

\date{{\small\file}\today}

\maketitle

\begin{abstract}\noi
Interacting particle systems can often be constructed from a graphical
representation, by applying local maps at the times of associated Poisson
processes. This leads to a natural coupling of systems started in different
initial states. We consider interacting particle systems on the complete graph
in the mean-field limit, i.e., as the number of vertices tends to infinity. We
are not only interested in the mean-field limit of a single process, but
mainly in how several coupled processes behave in the limit. This turns out to
be closely related to recursive tree processes as studied by Aldous and
Bandyopadyay in discrete time. We here develop an analogue theory for
recursive tree processes in continuous time.  We illustrate the abstract
theory on an example of a particle system with cooperative branching. This
yields an interesting new example of a recursive tree process that is not
endogenous.
\end{abstract}
\vspace{.5cm}

\noi
{\it MSC 2010.} Primary: 82C22, Secondary: 60J25, 60J80, 60K35.\\
{\it Keywords.} Mean-field limit, recursive tree process, recursive
distributional equation, endogeny, interacting particle systems, cooperative
branching.\\
{\it Acknowledgement.} Work sponsored by grant 16-15238S of the Czech Science
Foundation (GA CR) 
and by grant STU\ 527/1-2 of the German Research Foundation (DFG) within
the Priority Programme 1590 ``Probabilistic Structures in Evolution''.
 
%

\newpage

{\setlength{\parskip}{-2pt}\tableofcontents}

\newpage

\section{Introduction and main results}\label{S:intro}

\subsection{Introduction}\label{S:firstintro}

Let $\Om$ and $S$ be Polish spaces, let $\rb$ be a finite measure on
$\Om$ with total mass $|\rb|:=\rb(\Om)>0$, and let $\ga:\Om\times S^{\N_+}\to S$
be measurable, where $\N_+:=\{1,2,\dots\}$. Let $T$ be the operator acting on probability measures on $S$
defined as
\be\label{Tdef}
T(\mu):=\mbox{ the law of }\ga[\omb](X_1,X_2,\ldots),
\ee
where $\omb$ is an $\Om$-valued random variable with law $|\rb|^{-1}\rb$ and
$(X_i)_{i\geq 1}$ are i.i.d.\ with law $\mu$. In this paper, we will be
interested in the differential equation
\be\label{mean}
\dif{t}\mu_t=|\rb|\big\{T(\mu_t)-\mu_t\big\}\qquad(t\geq 0).
\ee
In Theorem~\ref{T:mean} below, we will prove existence and uniqueness of
solutions to (\ref{mean}) under the assumption that there exists a measurable
function $\kappa:\Om\to\N$ such that
\be\label{kappa}
\ga[\om](x_1,x_2,\ldots)=\ga[\om](x_1,\ldots,x_{\kappa(\om)})
\qquad(\om\in\Om,\ x\in S^{\N_+})
\ee
depends only on the first $\kappa(\om)\in\N$ coordinates, and
\be\label{sumr}
\int_\Om\rb(\di\om)\,\kappa(\om)<\infty.
\ee
Our interest in equation (\ref{mean}) stems from the fact that, as we will
prove in Theorem~\ref{T:meanlim} below, the \emph{mean-field limits} of a large
class of interacting particle systems are described by equations of the form
(\ref{mean}). In view of this, we call (\ref{mean}) a \emph{mean-field
  equation}. The analysis of this sort of equations is commonly the first step
towards understanding a given interacting particle system. Some illustrative
examples of mean-field equations in the literature are
\cite[(1.1)]{DN97}, \cite[(1.2)]{NP99}, and \cite[(4)]{FL17}.

In the special case that $\kappa(\om)=1$ for all $\om\in\Om$, we observe that
$T(\mu)=\int_S\mu(\di x)K(x,\,\cdot\,)$, where $K$ is the probability kernel
on $S$ defined as
\be\label{randmap}
K(x,A):=\P\big[\ga[\omb](x)\in A\big]\qquad(x\in S,\ A\sub S\mbox{ measurable}).
\ee
In view of this, if $\omb_1,\omb_2,\ldots$ are i.i.d.\ with law
$|\rb|^{-1}\rb$ and $X_0$ has law $\mu$, then setting
$X_k:=\ga[\omb_k](X_{k-1})$ $(k\geq 1)$ inductively defines a Markov chain
with transition kernel $K$, such that $X_k$ has law $T^k(\mu)$, where $T^k$
denotes the $k$-th iterate of the map $T$. Also, (\ref{mean}) describes
the forward evolution of a continuous-time Markov chain where random maps
$\ga[\om]$ are applied with Poisson rate $\rb(\di\om)$. A representation of a
probability kernel $K$ in terms of a random map $\ga[\omb]$ as in
(\ref{randmap}) is called a \emph{random mapping representation}.

More generally, when the function $\kappa$ is not identically one, Aldous and
Bandyopadhyay \cite{AB05} have shown that the iterates $T^k$ of the map $T$
from (\ref{Tdef}) can be represented in terms of a \emph{Finite Recursive Tree
  Process} (FRTP), which is a generalization of a discrete-time Markov chain
where time has a tree-like structure. More precisely, they construct a finite
tree of depth $k$ where the state of each internal vertex is a random function
of the states of its offspring. If the states of the leaves are i.i.d.\ with
law $\mu$, they show that the state at the root has law $T^k(\mu)$. They are
especially interested in fixed points of $T$, which generalize the concept of
an invariant law of a Markov chain. They show that each such fixed point $\nu$
gives rise to a \emph{Recursive Tree Process} (RTP), which is a process on an
infinite tree where the state of each vertex has law $\nu$. One can think of
such an RTP as a generalization of a stationary backward Markov chain
$(\ldots,X_{-2},X_{-1},X_0)$. A fixed point equation of the form $T(\nu)=\nu$
is called a \emph{Recursive Distributional Equation} (RDE).  Studying RDEs and
their solutions is of independent interest as they appear naturally in many
applications, see for example \cite{AB05, Als12}.

In the present paper, we develop an analogue theory in continuous time,
generalizing the concept of a continuous-time Markov chain to chains 
where time has a tree-like structure. Let $(T_t)_{t\geq 0}$ be the semigroup
defined by
\be\label{Tt}
T_t(\mu):=\mu_t\quad\mbox{where }(\mu_t)_{t\geq 0}\mbox{ solves (\ref{mean})
with }\mu_0=\mu.
\ee
In Theorem~\ref{T:meanrep}, we show that $T_t$ has a representation similar to
(\ref{Tdef}), namely
\be\label{Ttrep}
T_t(\mu)=\mbox{ the law of }G_t\big((X_\ibf)_{\ibf\in\T}\big),
\ee
where $\T$ is a countable set, $G_t:S^\T\to S$ is a random map, and the
$(X_\ibf)_{\ibf\in\T}$ are i.i.d.\ with law $\mu$ and independent of
$G_t$. Similar to what we have in (\ref{kappa}), the map $G_t$ does not depend
on all coordinates in $\T$ but only on a finite subcollection
$(X_\ibf)_{\ibf\in\nab\S_t}$. Here $(\nab\S_t)_{t\geq 0}$ turns out to be a
branching process and condition (\ref{sumr}) (which is not needed in the
discrete-time theory) guarantees that the offspring distribution of
this branching process has finite mean. 
  Similarly to (\ref{randmap}), we can view (\ref{Ttrep}) as a random mapping
representation of the operator in (\ref{Tt}).

As we have already mentioned, in Theorem~\ref{T:meanlim} below, we prove that
the mean-field limits of a large class of interacting particle systems are
described by equations of the form (\ref{mean}). These interacting particle 
systems are constructed by applying local maps at the times of
associated Poisson processes, which are introduced in detail in Section \ref{S:flow}.

We are not only interested in the mean-field limit of a
single process, but mainly in the mean-field limit of $n$ coupled processes
that are constructed from the same Poisson processes. For each $n\geq 1$, a
measurable map $g:S^k\to S$ gives rise to \emph{$n$-variate map}
$g^{(n)}:(S^n)^k\to S^n$ defined as
\be\label{nvar}
g^{(n)}\big(x_1,\ldots,x_k)=
g^{(n)}\big(x^1,\ldots,x^n):=\big(g(x^1),\ldots,g(x^n)\big)
\qquad(x^1,\ldots,x^n\in S^k),
\ee
where we denote an element of $(S^n)^k$ as
$(x^m_i)^{m=1,\ldots,n}_{i=1,\ldots,k}$ with $x_i=(x^1_i,\ldots,x^n_i)$
and $x^m=(x^m_1,\ldots,x^m_k)$.
Let $\Pc(S)$ denote the space of probability measures on a space $S$.
Letting $\ga^{(n)}[\om]$ denote the $n$-variate map associated with
  $\ga[\om]$ then, in analogy to (\ref{Tdef}),
\be\label{Tndef}
T^{(n)}(\mu):=\mbox{ the law of }\ga^{(n)}[\omb](X_1,\ldots,X_{\kappa(\omb)}),
\ee
defines an \emph{$n$-variate map} $T^{(n)}:\Pc(S^n)\to\Pc(S^n)$, which as in 
(\ref{mean}) gives rise to an \emph{$n$-variate mean-field equation},
which describes the mean-field limit of $n$ coupled processes.

If $X$ is an $S$-valued random variable whose law $\nu:=\P[X\in\,\cdot\,]$ is
a fixed point of $T$, then $\ov\nu^{(n)}:=\P[(X,\ldots,X)\in\,\cdot\,]$ is a
fixed point of $T^{(n)}$ that describes $n$ perfectly coupled processes. We
will be interested in the stability (or instability) of $\ov\nu^{(n)}$ under
the $n$-variate mean-field equation. In other words, for our mean-field
interacting particle systems, we fix the Poisson processes used in the
construction and want to know if small changes in the initial state lead to
small (or large) changes in the final state. Aldous and Bandyopadhyay
\cite{AB05} define an RTP to be \emph{endogenous} if the state at the root is
a measurable function of the random maps attached to all vertices of the
tree. They showed, in some precise sense (see Theorem~\ref{T:bivar} below),
that endogeny is equivalent to stability of $\ov\nu^{(n)}$. In
Theorem~\ref{T:bivar2}, we generalize their result to the continuous-time
setting.

The $n$-variate map $T^{(n)}$ is well-defined even for $n=\infty$, and
$T^{(\infty)}$ maps the space of all exchangeable probability laws on
$S^{\N_+}$ into itself. Let $\xi$ be a $\Pc(S)$-valued random variable with
law $\rho\in\Pc(\Pc(S))$, and conditional on $\xi$, let $(X^m)^{m=1,2,\ldots}$
be i.i.d.\ with common law $\xi$. Then the unconditional law of
$(X^m)^{m=1,2,\ldots}$ is exchangeable, and by De Finetti, each exchangeable
law on $S^{\N_+}$ is of this form. In view of this, $T^{(\infty)}$ naturally
gives rise to a map $\ch T:\Pc(\Pc(S))\to\Pc(\Pc(S))$ which is the
\emph{higher-level map} defined in \cite{MSS18}, and which analogously to
(\ref{mean}) gives rise to a \emph{higher-level mean-field equation}.
For any $\nu\in\Pc(S)$, let $\Pc(\Pc(S))_\nu$ denote the set of all
$\rho\in\Pc(\Pc(S))$ with mean $\int\rho(\di\mu)\mu=\nu$. In \cite{MSS18} it
is shown that if $\nu$ is a fixed point of $T$, then the corresponding
higher-level map $\ch T$ has two fixed points $\un\nu$ and $\ov\nu$ in
$\Pc(\Pc(S))_\nu$ that are minimal and maximal with respect to the
\emph{convex order}, defined in Theorem~\ref{T:Stras} below. Moreover,
$\un\nu=\ov\nu$ if and only if the RTP corresponding to $\nu$ is endogenous.

We will apply the theory developed here as well as in \cite{MSS18} to
the higher-level mean-field equation for a particular
interacting particle system with cooperative branching and deaths; see also
\cite{SS15,Mac17,BCH18} for several different variants of the model. To
formulate this properly, it is useful to introduce some more general notation.
Recall that for each $\om\in\Om$, $\ga[\om]$ is a map from $S^{\kappa(\om)}$
into $S$. We let
\be\label{Gi}
\Gi:=\big\{\ga[\om]:\om\in\Om\big\}
\ee
denote the set of all maps that can be obtained by varying $\om$. Here,
elements of $\Gi$ are measurable maps $g:S^k\to S$ where $k=k_g\geq 0$ may
depend on $g$. If $k=0$, then $S^0$ is defined to be a set with just one
element, which we denote by $\wurz$ (the empty sequence, which we distinguish
notationally from the empty set $\emptyset$). We equip $\Gi$
with the final \si-field for the map $\om\mapsto\ga[\om]$ and
let $\pi$ denote the image of the measure $\rb$ under this map.
Then the mean-field equation (\ref{mean}) can be rewritten as
\be\label{pimean}
\dif{t}\mu_t=\int_\Gi\pi(\di g)\big\{T_g(\mu_t)-\mu_t\big\}\qquad(t\geq 0),
\ee
where for any measurable map $g:S^k\to S$,
\be\label{Tgdef}
T_g(\mu):=\mbox{ the law of }g(X_1,\ldots,X_k),\quad
\mbox{where }(X_i)_{i=1,\ldots,k}\mbox{ are i.i.d.\ with law }\mu.
\ee

In the concrete example that we are interested in, $S:=\{0,1\}$ and
$\Gi:=\{\cob,\dth\}$ each have just two elements. Here $\cob:S^3\to S$ and
$\dth:S^0\to S$ are maps defined as
\be\label{cobdth}
\cob(x_1,x_2,x_3):=x_1\vee(x_2\wedge x_3)
\quand
\dth(\wurz):=0.
\ee
We choose
\be\label{picob}
\pi\big(\{\cob\}\big):=\al\geq 0\quand\pi\big(\{\dth\}\big):=1.
\ee
Then the mean-field equation (\ref{pimean}) takes the form
\be\label{coopmu11}
\dif{t}\mu_t=\al\big\{T_\cob(\mu_t)-\mu_t\big\}+\big\{T_\dth(\mu_t)-\mu_t\big\}.
\ee
which describes the mean-field limit of a particle system
with cooperative branching (with rate $\al$) and deaths (with rate 1). We will
see that for $\al>4$, (\ref{pimean}) has two stable fixed points $\nu_{\rm
  low},\nu_{\rm upp}$, and an unstable fixed point $\nu_{\rm mid}$ that
separates the domains of attraction of the stable fixed points.

In Theorem~\ref{T:coblev} below, we find all fixed points of the corresponding
higher-level mean-field equation, and determine their domains of
attraction. Note that solutions of the higher-level mean-field equation take
values in the probability measures on $\Pc(\{0,1\})\cong[0,1]$. As mentioned
before, each fixed point $\nu$ of the original mean-field equation gives rise
to two fixed points $\un\nu,\ov\nu$ of the higher-level mean-field equation,
which are minimal and maximal in $\Pc(\Pc(S))_\nu$ with respect to the convex
order. Moreover, $\un\nu=\ov\nu$ if and only if the RTP corresponding to $\nu$
is endogenous. In our example, we find that the stable fixed points $\nu_{\rm
  low},\nu_{\rm upp}$ give rise to endogenous RTPs, but the RTP associated
with $\nu_{\rm mid}$ is not endogenous. The higher-level equation has no other
fixed points in $\Pc(\Pc(S))_{\nu_{\rm mid}}$ except for $\un\nu_{\rm mid}$ and
$\ov\nu_{\rm mid}$, of which the former is stable and the latter
unstable. Numerical data for the nontrivial fixed point $\un\nu_{\rm mid}$
(viewed as a probability measure on $[0,1]$) are plotted in
Figure~\ref{fig:numid}.

\subsection{The mean-field equation}\label{S:meaneq}

In this subsection, we collect some basic results about the mean-field
equation (\ref{mean}) that form the basis for all that follows.
We interpret (\ref{mean}) in the following sense: letting
$\li\mu,\phi\re:=\int\phi\,\di\mu$, we say that a process $(\mu_t)_{t\geq 0}$
solves (\ref{mean}) if for each bounded measurable function
$\phi:S\to\R$, the function $t\mapsto\li\mu_t,\phi\re$ is continuously
differentiable and
\be\label{weakmean}
\dif{t}\li\mu_t,\phi\re=\int_\Om\rb(\di\om)
\big\{\li T_{\ga[\om]}(\mu_t),\phi\re-\li\mu_t,\phi\re\big\}\qquad(t\geq 0).
\ee
Our first result gives sufficient conditions for existence and uniqueness of
solutions to (\ref{mean}).

\bt[Mean-field equation]
Let\label{T:mean} $S$ and $\Om$ be Polish spaces, let $\rb$ be a nonzero
finite measure on $\Om$, and let $\ga:\Om\times S^{\N_+}\to S$ be measurable.
Assume that there exists a measurable function $\kappa:\Om\to\N$ such that
(\ref{kappa}) and (\ref{sumr}) hold. Then the mean-field equation (\ref{mean})
has a unique solution $(\mu_t)_{t\geq 0}$ for each initial state
$\mu_0\in\Pc(S)$.
\et

Theorem~\ref{T:mean} allows us to define a semigroup $(T_t)_{t\geq 0}$ of
operators $T_t:\Pc(S)\to\Pc(S)$ as in (\ref{Tt}). It is often useful to know
that solutions to (\ref{mean}) are continuous as a function of their initial
state. The following proposition gives continuity w.r.t.\ the total variation
norm $\|\,\cdot\,\|$ and moreover shows that if the constant $K$ from
(\ref{Kdef}) is negative, then the operators $(T_t)_{t\geq 0}$ form a
contraction semigroup.

\bp[Continuity in total variation norm]
Under\label{P:TVcont} the assumptions of Theorem~\ref{T:mean}, one has
\be\label{Ttcont}
\big\|T_t(\mu)-T_t(\nu)\big\|\leq e^{Kt}\|\mu-\nu\|
\qquad(\mu,\nu\in\Pc(S),\ t\geq 0),
\ee
where
\be\label{Kdef}
K:=\int_\Om\rb(\di\om)\,\big(\kappa(\om)-1\big).
\ee
\ep

Continuity w.r.t.\ weak convergence needs an additional assumption.

\bp[Continuity w.r.t.\ weak convergence] Assume\label{P:gacon}
that
\be\label{asco}
\rb\big(\{\om:\kappa(\om)=k,
\ \ga[\om]\mbox{ is discontinuous at }x\}\big)=0
\qquad(k\geq 0,\ x\in S^k).
\ee
Then the operator $T$ in (\ref{Tdef}) and the operators $T_t$ $(t\geq
  0)$  in (\ref{Tt}) are continuous w.r.t.\ the
topology of weak convergence.
\ep

The condition (\ref{asco}) is considerably weaker than the condition that
$\ga[\om]$ is continuous for all $\om\in\Om$. A simple example is
$\Om=S=[0,1]$, $\rb$ the Lebesgue measure, $\kappa\equiv 1$, and
$\ga[\om](x):=1_{\{x\geq\om\}}$.

\subsection{The mean-field limit}\label{S:flow}

\newcommand{\Omm}{\Om'}
\newcommand{\omi}{\om}

In this subsection, we show that equations of the form (\ref{mean}) arise as
the mean-field limits of a large class of interacting particle systems. In
order to be reasonably general, and in particular to allow for systems in
which more than one site can change its value at the same time, we will
introduce quite a bit of notation that will not be needed anywhere else in
Section~\ref{S:intro}, so impatient readers can just glance at
Theorem~\ref{T:meanlim} and the discussion surrounding (\ref{coopmean})
and skip the rest of this subsection.

Let $S$ be a Polish space as before, and let $N\in \N_+$. We will be
interested in continuous-time Markov processes taking values in $S^N$, where
$N$ is large. Denoting an element of $S^N$ by $x=(x_1,\ldots,x_N)$, we will
focus on processes with a high degree of symmetry, in the sense that their
dynamics are invariant under a permutation of the coordinates. It is
instructive, though not necessary for what follows, to view $\{1,\ldots,N\}$
as the vertex set of a complete graph, where all vertices are neighbors of
each other.
The basic ingredients we will use to describe our processes are:
\begin{enumerate}
\item a Polish space $\Omm$ equipped with a finite nonzero measure $\qb$,
\item a measurable function $\la:\Omm\to\N_+$,
\end{enumerate}
as well as, for each $\omi\in\Omm$ and $1\leq i\leq\la(\om)$,
\begin{enumerate}\addtocounter{enumi}{2}
\item a function $\ga_i[\omi]:S^{\la(\omi)}\to S$,
\item a finite set $K_i(\om)\sub\{1,\ldots,\la(\om)\}$ such that
$\ga_i[\omi](x_1,\ldots,x_{\la(\om)})
=\ga_i[\omi]\big((x_i)_{i\in K_i(\om)}\big)$ depends only on the coordinates
in $K_i(\om)$.
\end{enumerate}
Setting $\Omm_l:=\{\omi\in\Omm:\la(\omi)=l\}$, we assume that the functions
\be\label{vecmeas}
\Omm_l\times S^l\ni(\omi,x)\mapsto\ga_i[\omi](x)
\quand
\Omm_l\ni\omi\mapsto 1_{\{j\in K_i(\om)\}}
\quad\mbox{are measurable}
\ee
for each $1\leq i,j\leq l$. We let $\vec\ga[\omi]:S^{\la(\omi)}\to
S^{\la(\omi)}$ denote the function
\be
\vec\ga[\omi](x)
:=\big(\ga_1[\omi](x),\ldots,\ga_{\la(\om)}[\omi](x)\big)
\ee
and let $\kappa_i(\om):=|K_i(\om)|$ denote the cardinality of the set
$K_i(\om)$.

The space $\Omm$, measure $\qb$, and functions $\la$ and $\vec\ga$ play roles
similar, but not quite identical to $\Om,\rb,\kappa$, and $\ga$ from
Subsection~\ref{S:intro}. We can use $\Omm,\qb,\la$, and $\vec\ga$
to define the following mean-field equation:
\be\label{vecmean}
\dif{t}\mu_t
=\int_{\Omm}\!\qb(\di\om)\sum_{i=1}^{\la(\om)}
\big\{T_{\ga_i[\om]}(\mu_t)-\mu_t\big\}.
\ee
The following lemma says that (\ref{vecmean}) is really a mean-field equation
of the form we have already seen in (\ref{mean}). This is why in subsequent
sections we will only work with equations of this form.

\bl[Simplified equation]
Assume\label{L:simpmean} that
\be\label{rmom2}
{\rm(i)}\ \int_{\Omm\!}\qb(\di\om)\,\la(\om)<\infty
\quand
{\rm(ii)}\ \int_{\Omm\!}\qb(\di\om)\sum_{i=1}^{\la(\om)}\kappa_i(\om)<\infty.
\ee
Then equation (\ref{vecmean}) can be cast in the simpler form
(\ref{mean}) for a suitable choice of $\Om,\rb,\kappa$, and $\ga$, where
(\ref{rmom2})~(i) guarantees that $\rb$ is a finite measure and
(\ref{rmom2})~(ii) implies that (\ref{sumr}) holds. If
\be\label{qasco}
\qb\big(\{\omi:\la(\omi)=l,\ \ga_i[\om]\mbox{ is discontinuous at x}\}\big)=0
\qquad(1\leq i\leq l,\ x\in S^l),
\ee
then moreover (\ref{asco}) can be satisfied.
\el

We now use the ingredients $\Omm,\qb,\la$, and $\vec\ga$ to define the class
of Markov processes we are interested in. We construct these processes by
applying local maps, that affect only finitely many coordinates, at the times
of associated Poisson processes. In the context of interacting particle
systems, such constructions are called \emph{graphical representations}.

For any $N\in\N_+$ we set $[N]:=\{1,\ldots,N\}$. We let $[N]^{\li
  l\re}$ denote the set of all sequences $\ibf=(i_1,\ldots,i_l)$ for which
$i_1,\ldots,i_l\in[N]$ are all different. Note that $[N]^{\li l\re}$ has
$N^{\li l\re}:=N(N-1)\cdots(N-l+1)$ elements. We will consider
Markov processes $X=(X(t))_{t\geq 0}$ with values in $S^N$ that evolve in the
following way:
\begin{enumerate}
\item At the times of a Poisson process with intensity $|\qb|:=\qb(\Omm)$, an
  element $\omi\in\Omm$ is chosen according to the probability law
  $|\qb|^{-1}\qb$.
\item If $\la(\omi)>N$, nothing happens.
\item Otherwise, an element $\ibf\in[N]^{\li\la(\omi)\re}$ is selected
  according to the uniform distribution on $[N]^{\li\la(\omi)\re}$, and the
  previous values $\big(X_{i_1}(t-),\ldots,X_{i_{\la(\omi)}}(t-)\big)$ of $X$
  at the coordinates $i_1,\ldots,i_{\la(\omi)}$ are replaced by
  $\big(X_{i_1}(t),\ldots,X_{i_{\la(\omi)}}(t)\big)
  =\vec\ga[\omi]\big(X_{i_1}(t-),\ldots,X_{i_{\la(\omi)}}(t-)\big)$.
\end{enumerate}

More formally, we can construct our Markov process $X=(X(t))_{t\geq 0}$ as
follows. For each $\omi\in\Omm$ with $\la(\omi)\leq N$, and for each
$\ibf\in[N]^{\li\la(\omi)\re}$, define a map $m_{\omi,\ibf}:S^N\to S^N$ by
\be\label{momij}
m_{\omi,\ibf}(x)_j:=\left\{\ba{ll}
\ga_j[\omi](x_{i_1},\ldots,x_{i_{\la(\omi)}})
\quad\mbox{if }j\in\{i_1,\ldots,i_{\la(\omi)}\},\\
x_j\quad\mbox{otherwise,}\ea\right.
\qquad(x\in S^N).
\ee
Let $\Pi$ be a Poisson point set on
\be
\big\{(\omi,\ibf,t)
:\omi\in\Omm,\ \ibf\in[N]^{\li\la(\omi)\re},\ t\in\R\big\}
\ee
with intensity
\be
\qb(\di\omi)\,\frac{1_{\{\la(\omi)\leq N\}}}{N^{\li\la(\omi)\re}}\,\di t.
\ee
Since $\qb$ is a finite measure, the set
$\Pi_{s,u}:=\{(\omi,\ibf,t)\in\Pi:s<t\leq u\}$ is a.s.\ finite for each
$-\infty<s\leq u<\infty$, so we can order its elements as
\be
\Pi_{s,u}=\big\{(\omi_1,\ibf_1,t_1),\ldots,(\omi_n,\ibf_n,t_n)\big\}
\quad\mbox{with}\quad t_1<\cdots<t_n
\ee
and use this to define
\be\label{flowdef}
\Xb_{s,u}=m_{\omi_n,\ibf_n}\circ\cdots\circ m_{\omi_1,\ibf_1}
\ee
In words, $\Pi$ is a list of triples $(\omi,\ibf,t)$. Here $\omi$ represents
some external input that tells us that we need to apply the map
$\vec\ga[\omi]$. The coordinates where and the time when this map needs to be
applied are given by $\ibf$ and $t$, respectively. It is easy to see that
the random maps $(\Xb_{s,u})_{s\leq u}$ form a \emph{stochastic flow}, i.e.,
\be
\Xb_{s,s}=1\quand\Xb_{t,u}\circ\Xb_{s,t}=\Xb_{s,u}\qquad(s\leq t\leq u),
\ee
where $1$ denotes the identity map. Moreover $(\Xb_{s,u})_{s\leq u}$ has
\emph{independent increments} in the sense that
\be
\Xb_{t_1,t_2},\ldots,\Xb_{t_{k-1},t_k}\quad\mbox{are independent}
\ee
for each $t_1<\cdots<t_k$. It is well-known (see, e.g.,
\cite[Lemma~1]{SS18}) that if $X(0)$ is an $S^N$-valued random variable,
independent of the Poisson set $\Pi$, then setting
\be\label{XNdef}
X(t):=\Xb_{0,t}\big(X(0)\big)\qquad(t\geq 0)
\ee
defines a Markov process $X=(X(t))_{t\geq 0}$ with values in $S^N$. Note that
$(X(t))_{t\geq 0}$ has piecewise constant sample paths, which are
right-continuous because of the way we have defined $\Pi_{s,u}$.

We now formulate our result about the mean-field limit of
Markov processes as defined in (\ref{XNdef}). For any $x\in S^N$, we define an
\emph{empirical measure} $\mu\{x\}$ on $S$ by
\be\label{empir0}
\mu\{x\}:=\frac{1}{N}\sum_{i\in[N]}\de_{x_i}.
\ee
Below, $\mu^{\otimes n}:=\mu\otimes\cdots\otimes\mu$ denotes the product
measure of $n$ copies of $\mu$. The expectation $\E[\mu]$ of a random
measure $\mu$ on a Polish space $S$ is defined in the usual way, i.e.,
$\E[\mu]$ is the deterministic measure defined by
$\int\!\phi\,\di\E[\mu]:=\E[\int\!\phi\,\di\mu]$ for any bounded measurable
$\phi:S\to\R$. 
\bt[Mean-field limit]
Let\label{T:meanlim} $S$ be a Polish space, let $\Omm,\qb,\la$, and $\vec\ga$
be as above, and assume (\ref{rmom2}). For each $N\in\N_+$, let
$(X^{(N)}(t))_{t\geq 0}$ be Markov processes with state space $S^N$ as defined
in (\ref{XNdef}), and let 
$\mu^N_t:=\mu\{X^{(N)}(t)\}$
denote their associated
empirical measures. Let $d$ be any metric on $\Pc(S)$ that generates the
topology of weak convergence. Fix some (deterministic) $\mu_0\in\Pc(S)$ and
assume that (at least) one of the following two conditions is satisfied.
\begin{enumerate}
\item $\dis\P\big[d(\mu^N_0,\mu_0)\geq\eps]\asto{N}0$ for all $\eps>0$, and
 (\ref{qasco}) holds.
\item $\big\|\E[(\mu^N_0)^{\otimes n}]-\mu_0^{\otimes n}\big\|\asto{N}0$ for
  all $n\geq 1$, where $\|\,\cdot\,\|$ denotes the total variation norm.
\end{enumerate}
Then
\be\label{convmu2}
\P\big[\sup_{0\leq t\leq T}d(\mu^N_{Nt},\mu_t)\geq\eps\big]\asto{N}0
\qquad(\eps>0,\ T<\infty),
\ee
where $(\mu_t)_{t\geq 0}$ is the unique solution to the mean-field equation
(\ref{vecmean}) with initial state $\mu_0$.
\et

Condition~(ii) is in particular satisfied if
$X^N_1(0),\ldots,X^N_N(0)$ are i.i.d.\ with common law $\mu_0$.
Note that in (\ref{convmu2}), we rescale time by a factor $N$.

It is instructive to demonstrate the general set-up on our concrete example of
a particle system with cooperative branching and deaths. As before, we have
$S=\{0,1\}$. We choose for $\Omm$ a set with just two elements, say
$\Omm=\{1,2\}$, and we set $\qb(\{1\}):=\al\geq 0$ and $\qb(\{2\}):=1$. We let
$\la(1):=3$, $\la(2):=1$, and define $\vec\ga[1]:S^3\to S^3$ and
$\vec\ga[2]:S\to S$ by
\be
\vec\ga[1](x_1,x_2,x_3):=\big(x_1\vee(x_2\wedge x_3),x_2,x_3\big)
\quand
\vec\ga[2](x_1):=0.
\ee
Then the particle system in (\ref{XNdef}) has the following description. Let
us say that a site $i$ is occupied at time $t$ if $X_i(t)=1$. Then, with rate
$\al$, three sites $(i_1,i_2,i_3)\in[N]^{\li 3\re}$ are
selected at random. If the sites $i_2$ and $i_3$ are both occupied, then the
particles at these sites cooperate to produce a third particle at $i_1$,
provided this site is empty. In addition, with rate 1, a site $i$ is selected
at random, and any particle that is present there dies.

It is not hard to see that for our choice of $\Omm,\qb,\la$, and $\vec\ga$,
the mean-field equation (\ref{vecmean}) simplifies to (\ref{coopmu11}), 
Note that since $\ga_2[1]$ and $\ga_3[1]$ are the identity map, they drop out
of (\ref{vecmean}), so only $\ga_1[1]=\cob$ and $\ga_1[2]=\dth$
remain. Since $\ga_1[2](x_1)=0$ regardless of the value of $x_1$, we can
choose for $K_1(1)$ the empty set and view $\ga_1[2]=\dth$ as a function
$\dth:S^0\to S$.

Solutions of (\ref{coopmu11}) take values in the probability measures on
$S=\{0,1\}$, which are uniquely characterized by their value at 1. Rewriting
(\ref{coopmu11}) in terms of $p_t:=\mu_t(\{1\})$ yields the equation
\be\label{coopmean}
\dif{t}p_t=\al p_t^2(1-p_t)-p_t\qquad(t\geq 0).
\ee
This equation can also be found in \cite[(1.11)]{Nob92}, \cite[(1.2)]{Neu94},
\cite[(3.1)]{BW97}, \cite[(4)]{FL17}, and \cite[(2.1)]{BCH18}.
It is not hard to check that for $\al<4$, the only fixed
point of (\ref{coopmean}) is $z_{\rm low}:=0$, while for $\al\geq 4$, there are
additional fixed points
\be\label{fixpts}
z_{\rm mid}:=\ffrac{1}{2}-\sqrt{\ffrac{1}{4}-\ffrac{1}{\al}}
\quand
z_{\rm upp}:=\ffrac{1}{2}+\sqrt{\ffrac{1}{4}-\ffrac{1}{\al}}.
\ee
If $\al<4$, then solutions to (\ref{coopmean}) converge to $z_{\rm low}$
regardless of the initial state. On the other hand, for $\al\geq 4$, solutions
to (\ref{coopmean}) with $p_0>z_{\rm mid}$ converge to the upper fixed point
$z_{\rm upp}$ while solutions to (\ref{coopmean}) with $p_0<z_{\rm mid}$
converge to the lower fixed point $z_{\rm low}$. In particular, if $\al>4$,
then $z_{\rm low}$ and $z_{\rm upp}$ are stable fixed points while $z_{\rm
  mid}$ is an unstable fixed point separating the domains of attraction of
$z_{\rm low}$ and $z_{\rm upp}$.

\subsection{A recursive tree representation}\label{S:treerep}

In this subsection we formally introduce Finite Recursive Tree Processes
(FRTPs) and state the random mapping representation of solutions to the
mean-field equation (\ref{mean}) anticipated in (\ref{Ttrep}).

\newcommand{\ovT}{{\mathbb T}}
\newcommand{\St}{{\mathbb S}}

For $d\in\N_+$, let $\ovT^d$ denote the space of all finite words
$\ibf=i_1\cdots i_n$ $(n\in\N)$ made up from the alphabet $\{1,\ldots,d\}$,
and define $\ovT^\infty$ similarly, using the alphabet $\N_+$. If
$\ibf,\jbf\in\ovT^d$ with $\ibf=i_1\cdots i_m$ and $\jbf=j_1\cdots j_n$, then
we define the concatenation $\ibf\jbf\in\ovT^d$ by $\ibf\jbf:=i_1\cdots
i_mj_1\cdots j_n$. We denote the length of a word $\ibf=i_1\cdots i_n$ by
$|\ibf|:=n$ and let $\wurz$ denote the word of length zero. We view
$\ovT^d$ as a tree with root $\wurz$, where each vertex $\ibf\in\ovT^d$ has
$d$ children $\ibf 1,\ibf 2,\ldots$, and each vertex $\ibf=i_1\cdots i_n$
except the root has precisely one ancestor $\lib:=i_1\cdots i_{n-1}$.
For each \emph{rooted subtree} of $\ovT^d$, i.e., a subtree $\U\sub\ovT^d$ that
contains $\wurz$, we let
$\pa\U:=\{\ibf\in\ovT^d:\lib\in\U,\ \ibf\not\in\U\}$ denote the
\emph{boundary} of $\U$ relative to $\ovT^d$. We write
\be\label{Tdistn}
\ovT^d_{(n)}:=\{\ibf\in\ovT^d:|\ibf|<n\}
\quand
\pa\ovT^d_{(n)}=\{\ibf\in\ovT^d:|\ibf|=n\}\qquad(n\geq 1),
\ee
and use the convention $\pa\emptyset:=\{\wurz\}$, so that (\ref{Tdistn}) holds
also for $n=0$.

We return to the set-up of Subsection~\ref{S:firstintro}, i.e., $S$ and $\Om$
are Polish spaces, $\rb$ is a nonzero finite measure on $\Om$, and
$\ga:\Om\times S^{\N_+}\to S$ and $\kappa:\Om\to\N$ are measurable functions
such that (\ref{kappa}) holds. We fix some $d\in\N_+\cup\{\infty\}$ such that
$\kappa(\om)\leq d$ for all $\om\in\Om$ and set $\ovT:=\ovT^d$. Let 
$(\omb_\ibf)_{\ibf\in\T}$ be an i.i.d.\ collection of $\Om$-valued r.v.'s with
common law $|\rb|^{-1}\rb$. Fix $n\geq 1$ and assume that
\be\ba{rl}\label{FRTP}
{\rm(i)}&\mbox{the $(X_\ibf)_{\ibf\in{\pa\T_{(n)}}}$ are i.i.d.\ with common
 law $\mu$ and independent of $(\omb_\ibf)_{\ibf\in\T_{(n)}}$,}\\[5pt]
{\rm(ii)}&\dis X_\ibf:=\ga[\omb_\ibf](X_{\ibf 1},\ldots,X_{\ibf\kappa(\omb_\ibf)})
\qquad(\ibf\in\T_{(n)}).
\ec
Then it is easy to see that the law of $X_\wurz$ is given by $T^n(\mu)$, where
$T^n$ is the $n$-th iterate of the operator in (\ref{Tdef}). We call the
collection of random variables
\be
\big((\omb_\ibf)_{\ibf\in\T_{(n)}},(X_\ibf)_{\ibf\in\T_{(n)}\cup\pa\T_{(n)}}\big)
\ee
a \emph{Finite Recursive Tree Process} (FRTP). We can think of
$(X_\ibf)_{\ibf\in\T_{(n)}\cup\pa\T_{(n)}}$ as a generalization of a Markov
chain, where time has a tree-like structure.

We now aim to give a similar representation of the semigroup $(T_t)_{t\geq 0}$
from (\ref{Tt}). To do this, we let $(\sig_\ibf)_{\ibf\in\T}$ be
i.i.d.\ exponentially distributed random variables with mean $|\rb|^{-1}$. We
interpret $\sig_\ibf$ as the lifetime of the individual with index $\ibf$ and
let
\be\label{taudef}
\tau^\ast_{\ibf}:=\sum_{m=1}^{n-1}\sig_{i_1\cdots i_m}
\quand\tau^\dgg_\ibf:=\tau^\ast_\ibf+\sig_\ibf
\qquad(\ibf=i_1\cdots i_n).
\ee
denote the times when the individual $\ibf$ is born and dies, respectively.
Then
\be\label{Ttdef}
\ovT_t:=\big\{\ibf\in\ovT:\tau^\dgg_\ibf\leq t\big\}
\quand
\pa\ovT_t=\big\{\ibf\in\ovT:\tau^\ast_\ibf\leq t<\tau^\dgg_\ibf\big\}
\qquad(t\geq 0)
\ee
are the (random) subtrees of $\ovT$ consisting of all individuals that have
died before time $t$, resp.\ are alive at time $t$. If the function $\kappa$
from (\ref{kappa}) is bounded, then we can choose $\ovT:=\ovT^d$ with
$d<\infty$. Now it is easy to check that $(\pa\ovT_t)_{t\geq 0}$ is a
continuous-time branching process where each particle is with rate $|\rb|$
replaced by $d$ new particles. In particular, $\ovT_t$ is a.s.\ finite for
each $t>0$. On the other hand, when $\kappa$ is unbounded, we need to choose
$\ovT:=\ovT^\infty$, and this has the consequence that $\ovT_t$ is
a.s.\ infinite for each $t>0$. Nevertheless, under the assumption
(\ref{sumr}), it turns out that only a finite subtree of $\ovT_t$ is relevant
for the state at the root $X_\wurz$, as we explain now.

Let $\St$ be the random subtree of $\ovT$ defined as
\be\label{Sdef}
\St:=\big\{i_1\cdots i_n\in\ovT:
i_m\leq\kappa(\omb_{i_1\cdots i_{m-1}})\ \forall 1\leq m\leq n\big\},
\ee
and for each subtree $\U\sub\St$, let $\nab\U:=\{\ibf\in\St:\lib\in\U,
\ \ibf\not\in\U\}$ denote the outer boundary of $\U$ relative to $\St$,
where again we use the convention that $\nab\U:=\{\wurz\}$ if $\U$ is the
empty set. Then, under condition (\ref{sumr}),
\be\label{Stdef}
\St_t:=\ovT_t\cap\St
\quand
\nab\St_t=\big\{\ibf\in\St:\tau^\ast_\ibf\leq t<\tau^\dgg_\ibf\big\}
\qquad(t\geq 0)
\ee
are a.s.\ finite for all $t\geq 0$. Indeed, $(\nab\St_t)_{t\geq 0}$ is a
branching process where for each individual $\ibf$, with Poisson rate
$\rb(\di\om)$, an element $\om\in\Om$ is selected and $\ibf$ is replaced by
new individuals $\ibf 1,\ldots,\ibf\kappa(\om)$. The condition on the rates
(\ref{sumr}) guarantees that this branching process has finite mean and in
particular does not explode, so that $\St_t$ is a.s.\ a finite subtree of $\St$.

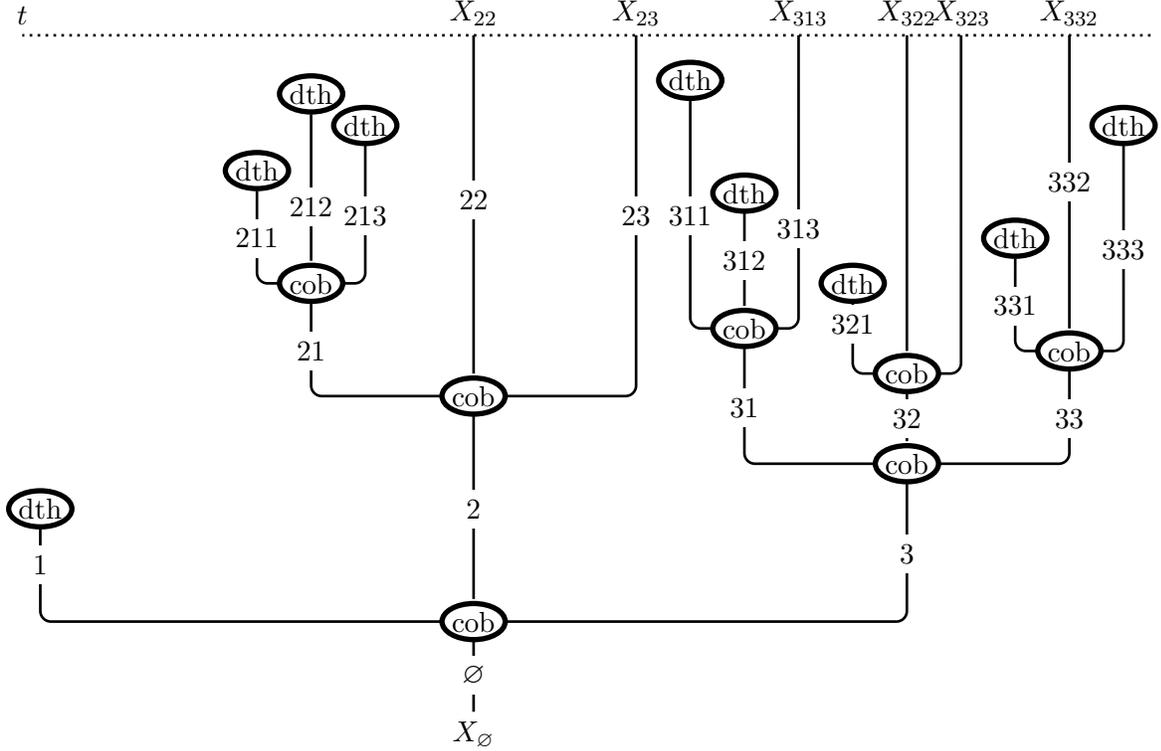
\begin{figure}
\centering

\begin{tikzpicture}[scale=0.6]

\node (wurz) at (0,2) {$\rm{cob}$};
\node (1) at (-9.6,4.5) {$\rm{dth}$};
\node (2) at (0,7) {$\rm{cob}$};
\node (3) at (9.6,5.5) {$\rm{cob}$};
\node (21) at (-3.6,9.5) {$\rm{cob}$};
\node[above] (22) at (0,15) {$X_{22}$};
\node[above] (23) at (3.6,15) {$X_{23}$};
\node (31) at (6,8.5) {$\rm{cob}$};
\node (32) at (9.6,7.5) {$\rm{cob}$};
\node (33) at (13.2,8) {$\rm{cob}$};
\node (211) at (-4.8,12) {$\rm{dth}$};
\node (212) at (-3.6,13.7) {$\rm{dth}$};
\node (213) at (-2.4,13) {$\rm{dth}$};
\node (311) at (4.8,14) {$\rm{dth}$};
\node (312) at (6,11.5) {$\rm{dth}$};
\node[above] (313) at (7.2,15) {$X_{313}$};
\node (321) at (8.4,9.5) {$\rm{dth}$};
\node[above] (322) at (9.6,15) {$X_{322}$};
\node[above] (323) at (10.8,15) {$X_{323}$};
\node (331) at (12,10.5) {$\rm{dth}$};
\node[above] (332) at (13.2,15) {$X_{332}$};
\node (333) at (14.4,13) {$\rm{dth}$};

\draw[line width=1pt] (0,0)--(wurz);
\draw (0,0.8) node[fill=white] {$\wurz$};
\draw[line width=2pt] (wurz) ellipse (0.7cm and 0.4cm);
\node[below] at (0,0) {$X_\wurz$};

\draw[line width=1pt,rounded corners] (-0.7,2)--(-9.6,2)--(1);
\draw (-9.6,3.25) node[fill=white] {$1$};
\draw[line width=2pt] (1) ellipse (0.7cm and 0.4cm);

\draw[line width=1pt] (0,2.5)--(2);
\draw (0,4.5) node[fill=white] {$2$};
\draw[line width=2pt] (2) ellipse (0.7cm and 0.4cm);

\draw[line width=1pt,rounded corners] (0.7,2)--(9.6,2)--(3);
\draw (9.6,3.5) node[fill=white] {$3$};
\draw[line width=2pt] (3) ellipse (0.7cm and 0.4cm);

\draw[line width=1pt,rounded corners] (-0.7,7)--(-3.6,7)--(21);
\draw (-3.6,8) node[fill=white] {$21$};
\draw[line width=2pt] (21) ellipse (0.7cm and 0.4cm);

\draw[line width=1pt] (0,7.5)--(22);
\draw (0,11.35) node[fill=white] {$22$};

\draw[line width=1pt,rounded corners] (0.7,7)-- (3.6,7)--(23);
\draw (3.6,11) node[fill=white] {$23$};

\draw[line width=1pt,rounded corners] (8.9,5.5)--(6,5.5)--(31);
\draw (6,6.75) node[fill=white] {$31$};
\draw[line width=2pt] (31) ellipse (0.7cm and 0.4cm);

\draw[line width=1pt] (9.6,6)--(32);
\draw (9.6,6.5) node[fill=white] {$32$};
\draw[line width=2pt] (32) ellipse (0.7cm and 0.4cm);

\draw[line width=1pt,rounded corners] (10.3,5.5)--(13.2,5.5)--(33);
\draw (13.2,6.5) node[fill=white] {$33$};
\draw[line width=2pt] (33) ellipse (0.7cm and 0.4cm);

\draw[line width=1pt,rounded corners] (-4.3,9.5)--(-4.8,9.5)--(211);
\draw (-4.8,10.5) node[fill=white] {$211$};
\draw[line width=2pt] (211) ellipse (0.7cm and 0.4cm);

\draw[line width=1pt,rounded corners] (-3.6,10)--(212);
\draw (-3.6,11.2) node[fill=white] {$212$};
\draw[line width=2pt] (212) ellipse (0.7cm and 0.4cm);

\draw[line width=1pt,rounded corners] (-2.9,9.5)--(-2.4,9.5)--(213);
\draw (-2.4,11) node[fill=white] {$213$};
\draw[line width=2pt] (213) ellipse (0.7cm and 0.4cm);

\draw[line width=1pt,rounded corners] (5.3,8.5)--(4.8,8.5)--(311);
\draw (4.8,11) node[fill=white] {$311$};
\draw[line width=2pt] (311) ellipse (0.7cm and 0.4cm);

\draw[line width=1pt,rounded corners] (6,9)--(312);
\draw (6,10) node[fill=white] {$312$};
\draw[line width=2pt] (312) ellipse (0.7cm and 0.4cm);

\draw[line width=1pt,rounded corners] (6.7,8.5)--(7.2,8.5)--(313);
\draw (7.2,10.7) node[fill=white] {$313$};

\draw[line width=1pt,rounded corners] (8.9,7.5)--(8.4,7.5)--(321);
\draw (8.4,8.6) node[fill=white] {$321$};
\draw[line width=2pt] (321) ellipse (0.7cm and 0.4cm);

\draw[line width=1pt] (9.6,8)--(322);

\draw[line width=1pt,rounded corners] (10.3,7.5)--(10.8,7.5)--(323);

\draw[line width=1pt,rounded corners] (12.5,8)--(12,8)--(331);
\draw (12,9) node[fill=white] {$331$};
\draw[line width=2pt] (331) ellipse (0.7cm and 0.4cm);

\draw[line width=1pt,rounded corners] (13.2,8.5)--(332);
\draw (13.2,11.75) node[fill=white] {$332$};

\draw[line width=1pt,rounded corners] (13.9,8)--(14.4,8)--(333);
\draw (14.4,10.25) node[fill=white] {$333$};
\draw[line width=2pt] (333) ellipse (0.7cm and 0.4cm);

\draw[line width=1pt, dotted](15,15)--(-10,15) node[above]{$t$};

\end{tikzpicture}
\caption{A particular realization of the branching process $\nab\St_t$ for a
  system with cooperative branching and deaths. The random map $G_{\St_t}$ is
  the concatenation of random maps attached to the vertices of the family
  tree $\St_t$ of the individuals alive at time $t$. In this example,
  $\nab\St_t=\{22,23,313,322,323,332\}$ and the maps $\cob$ and $\dth$ are
  as defined in (\ref{cobdth}).}
\label{fig:treemap}
\end{figure}

Let $(\omb_\ibf)_{\ibf\in\T}$ be i.i.d.\ with common law $|\rb|^{-1}\rb$,
independent of the lifetimes $(\sig_\ibf)_{\ibf\in\T}$.
For any finite rooted subtree $\U\sub\St$ and for each
$(x_\ibf)_{\ibf\in\nab\U}=x\in S^{\nab\U}$, we can inductively define $x_\ibf$
for $\ibf\in\U$ by
\be\label{xinduc}
x_\ibf:=\ga[\omb_\ibf](x_{\ibf 1},\ldots,x_{\ibf\kappa(\omb_\ibf)})
\qquad(\ibf\in\U).
\ee
Then the value $x_\wurz$ we obtain at the root is a
function of $(x_\ibf)_{\ibf\in\nab\U}$. Let us denote this function by
$G_\U:S^{\nab\U}\to S$, i.e.,
\be\label{GUdef}
G_\U\big((x_\ibf)_{\ibf\in\nab\U}\big):=x_{\wurz},
\quad\mbox{with $(x_\ibf)_{\ibf\in\U}$ defined as in (\ref{xinduc}).}
\ee
We can think of $G_\U$ as the ``concatenation'' of the maps
$(\ga[\omb_\ibf])_{\ibf\in\U}$. We will in particular be interested in the
random maps
\be\label{Gtdef}
G_t:=G_{\St_t}\qquad(t\geq 0)
\ee
with $\St_t$ as in (\ref{Stdef}). For our running example of a system with
cooperative branching and deaths, these definitions are illustrated in
Figure~\ref{fig:treemap}.

Let $(\Fi_t)_{t\geq 0}$, defined as
\be\label{Fidef}
\Fi_t:=
\sig\big(\nab\St_t,(\omb_\ibf,\sig_\ibf)_{\ibf\in\St_t}\big)\qquad(t\geq 0).
\ee
be the natural filtration associated with our evolving marked tree, that
contains information about which individuals are alive at time $t$, as well as
the random elements $\omb_\ibf$ and lifetimes $\sig_\ibf$ associated with all
individuals that have died by time $t$. In particular, $G_t$ is measurable
w.r.t.\ $\Fi_t$. The following theorem is a precise formulation of the random
mapping representation of solutions of the mean-field equation (\ref{mean}),
anticipated in (\ref{Ttrep}).

\bt[Recursive tree representation]
Let\label{T:meanrep} $S$ and $\Om$ be Polish spaces, let $\rb$ be a nonzero
finite measure on $\Om$, and let $\ga:\Om\times S^{\N_+}\to S$ and
$\kappa:\Om\to\N$ be measurable functions satisfying (\ref{kappa})
and (\ref{sumr}). Let $(\omb_\ibf)_{\ibf\in\St}$ be i.i.d.\ with
common law $|\rb|^{-1}\rb$ and let $(\sig_\ibf)_{\ibf\in\St}$ be an
independent i.i.d.\ collection of exponentially distributed random variables
with mean $|\rb|^{-1}$. Fix $t\geq 0$ and let $G_t$ and $\Fi_t$ be defined as in 
(\ref{Gtdef}) and (\ref{Fidef}). Conditional on $\Fi_t$, let
$(X_\ibf)_{\ibf\in\nab\St_t}$ be i.i.d.\ $S$-valued random variables with
common law $\mu$. Then
\be\label{meanrep}
T_t(\mu)=\mbox{ the law of }G_t\big((X_\ibf)_{\ibf\in\nab\St_t}\big),
\ee
where $T_t$ is defined in (\ref{Tt}).
\et

Recalling the definition of $G_t$, we can also formulate
Theorem~\ref{T:meanrep} as follows. With $(\omb_\ibf,\sig_\ibf)_{\ibf\in\St}$
as above, fix $t>0$ and let $(X_\ibf)_{\ibf\in\St_t\cup\nab\St_t}$ be random
variables such that
\be\ba{rl}\label{finrec}
{\rm(i)}&\mbox{Conditional on $\Fi_t$, the r.v.'s
  $(X_\ibf)_{\ibf\in\nab\St_t}$ are i.i.d.\ with common law $\mu$,}\\[5pt]
{\rm(ii)}&\dis X_\ibf=\ga[\omb_\ibf]
\big(X_{\ibf1},\ldots,X_{\ibf\kappa(\omb_\ibf)}\big)
\qquad(\ibf\in\St_t).
\ec
Then (\ref{meanrep}) says that the state at the root $X_\wurz$ has law
$T_t(\mu)$. This is a continuous-time analogue of the FRTP (\ref{FRTP}).

In our proofs, we will first prove Theorem~\ref{T:meanrep} and then use this
to prove Theorem~\ref{T:meanlim} about the mean-field limit of interacting
particle systems. Recall that these particle systems are constructed from a
stochastic flow $(\Xb_{s,t})_{s\leq t}$ as in (\ref{XNdef}). To find the
empirical measure of $X(t)=\Xb_{0,t}(X(0))$, we pick a site $i\in[N]$ at
random and ask for its type $X_i(t)$ which via $\Xb_{0,t}$ is a function of
the initial state $X(0)$. When $N$ is large, $X_i(t)$ does not depend on all
coordinates $(X_j(0))_{j\in[N]}$ but only on a random subset of them, and
indeed one can show that the map that gives $X_i(t)$ as a function of these
coordinates approximates the map $G_t$ from Theorem~\ref{T:meanrep}, in an
appropriate sense. The heuristics behind this are explained in some more
detail in Subsection~\ref{S:proofline} below.\med

\begin{remark}
Another way to write (\ref{meanrep}) is
\be\label{oldmeanrep}
\mu_t=\E\big[T_{G_t}(\mu_0)\big]\qquad(t\geq 0),
\ee
where $T_{G_t}$ is defined as in (\ref{Tgdef}) for the random map $G_t$
and $(\mu_t)_{t\geq 0}$ is a solution to (\ref{mean}).
One can check that $(\nab\St_t,G_t)_{t\geq 0}$ is a Markov process.
Let us informally denote this process by $(G_t)_{t\geq 0}$ and its state space
by $\Gi$. Then equation (\ref{meanrep}) can be understood as a (generalized)
duality relationship between $(G_t)_{t\geq 0}$ and $(\mu_t)_{t\geq 0}$ with
(generalized) duality function $H:\Gi\times\Pc(S)\to\Pc(S)$ given by
\be\label{Hdef}
H(g,\mu):=T_{g}(\mu).
\ee
With this definition, using the fact that $G_0$ is the identity map,
(\ref{oldmeanrep}) reads
\be\label{Hdual}
\mu_t=H(G_0,\mu_t)=\E[H(G_t,\mu_0)]
\ee
and we can obtain a family of usual (real-valued) dualities by integrating
against a test function~$\phi$.
\end{remark}

\subsection{Recursive tree processes}

Recall the definition of the operator $T$ in (\ref{Tdef}) and the semigroup
$(T_t)_{t\geq 0}$ in (\ref{Tt}). It is clear from (\ref{mean}) that for a
measure $\nu\in\Pc(S)$, the following two conditions are equivalent:
\be\label{RDE}
{\rm(i)}\quad T_t(\nu)=\nu\quad(t\geq 0)
\qquad\mbox{and}\qquad
{\rm(ii)}\quad T(\nu)=\nu.
\ee
We call such a measure $\nu$ a \emph{fixed point} of the mean-field equation
(\ref{mean}). Condition~(ii) is equivalent to saying that a random variable $X$
with law $\nu$ satisfies
\be\label{isd}
X\isd\ga[\om](X_1,\ldots,X_{\kappa(\om)}),
\ee
where $\isd$ denotes equality in distribution, $X_1,X_2,\ldots$ are
i.i.d.\ copies of $X$, and $\om$ is an independent $\Om$-valued random variable
with law $|\rb|^{-1}\rb$. Equations of this type are called \emph{Recursive
  Distributional Equations} (RDEs).

FRTPs as in (\ref{FRTP}) are consistent in the sense that if
$(X_\ibf)_{\ibf\in{\pa\T_{(n)}}}$ are as in (\ref{FRTP}), then for any $1\leq
m\leq n$,
\be\ba{rl}\label{FRTPconsist}
{\rm(i)}&\mbox{the $(X_\ibf)_{\ibf\in{\pa\T_{(m)}}}$ are i.i.d.\ with common
 law $T^{n-m}(\mu)$ and independent of $(\omb_\ibf)_{\ibf\in\T_{(m)}}$,}\\[5pt]
{\rm(ii)}&\dis X_\ibf:=\ga[\omb_\ibf](X_{\ibf 1},\ldots,X_{\ibf\kappa(\omb_\ibf)})
\qquad(\ibf\in\T_{(m)}).
\ec
The following lemma states a similar consistency property in the
continuous-time setting.

\bl[Consistency]
Fix\label{L:consis} $t>0$ and let $(X_\ibf)_{\ibf\in\St_t\cup\nab\St_t}$ be as in
(\ref{finrec}). \\
Then, for each $s\in[0,t]$,
\be\ba{rl}\label{consis}
{\rm(i)}&\mbox{conditional on $\Fi_s$, the r.v.'s $(X_\ibf)_{\ibf\in\nab\St_s}$
 are i.i.d.\ with common law $T_{t-s}(\mu)$,}\\[5pt]
{\rm(ii)}&\dis X_\ibf=\ga[\omb_\ibf]
\big(X_{\ibf1},\ldots,X_{\ibf\kappa(\omb_\ibf)}\big)
\qquad(\ibf\in\St_s),
\ec
where $(T_t)_{t\geq 0}$ is defined in (\ref{Tt}).
\el

Using the consistency relation (\ref{FRTPconsist}) and Kolmogorov's extension
theorem, it is not hard to see that if $\nu$ solves the RDE (\ref{RDE}), then
it is possible to define a stationary recursive process on an infinite tree
such that each vertex has law $\nu$. This was already observed in
\cite{AB05}. The following lemma is a slight reformulation of their
observation.

\bl[Recursive Tree Process]
Let\label{L:RTP} $\nu$ be a solution to the RDE (\ref{RDE}). Then there exists
a collection $(\omb_\ibf,X_\ibf)_{\ibf\in\ovT}$ of random variables whose
joint law is uniquely characterized by the following requirements.
\be\ba{rl}\label{RTP}
{\rm(i)}&\mbox{$(\omb_\ibf)_{\ibf\in\ovT}$ is an i.i.d.\ collection of
  $\Om$-valued r.v.'s with common law }|\rb|^{-1}\rb.\\[5pt]
{\rm(ii)}&\mbox{For each finite rooted subtree $\U\sub\ovT$, the
  r.v.'s $(X_\ibf)_{\ibf\in\poa\U}$ are i.i.d.\ with common}\\
&\mbox{law $\nu$ and independent of $(\omb_\ibf)_{\ibf\in\U}$.}\\[5pt]
{\rm(iii)}&\dis X_\ibf=\ga[\omb_\ibf]
\big(X_{\ibf1},\ldots,X_{\ibf\kappa(\omb_\ibf)}\big)
\qquad(\ibf\in\ovT).
\ec
\el

We call a collection of random variables $(\omb_\ibf,X_\ibf)_{\ibf\in\ovT}$ as
in Lemma~\ref{L:RTP} the \emph{Recursive Tree Process} (RTP) corresponding to
the map $\ga$ and the solution $\nu$ of the RDE (\ref{RDE}). We can view such
an RTP as a generalization of a stationary backward Markov chain. For most
purposes, we will only need the random variables $\omb_\ibf,X_\ibf$ with
$\ibf\in\St$, the random subtree defined in (\ref{Sdef}). The following
proposition shows that by adding independent exponential lifetimes to an RTP,
we obtain a stationary version of (\ref{consis}).

\bp[Continuous-time RTP]
Let\label{P:contRTP} $(\omb_\ibf,X_\ibf)_{\ibf\in\ovT}$ be an RTP corresponding
to a solution $\nu$ of the RDE (\ref{RDE}), and let
$(\sig_\ibf)_{\ibf\in\ovT}$ be an independent i.i.d.\ collection of exponentially
distributed random variables with mean $|\rb|^{-1}$. Then, for each $t\geq 0$,
\be\ba{rl}\label{contRTP}
{\rm(i)}&\mbox{conditional on $\Fi_t$, the r.v.'s $(X_\ibf)_{\ibf\in\nab\St_t}$ are
i.i.d.\ with common law $\nu$,}\\[5pt]
{\rm(ii)}&\dis X_\ibf=\ga[\omb_\ibf]\big(X_{\ibf1},\ldots,X_{\ibf\kappa(\omb_\ibf)}\big)
\qquad(\ibf\in\St).
\ec
\ep

At the end of Subsection~\ref{S:flow} we have seen that in our example of a
system with cooperative branching, the RDE (\ref{RDE}) has three solutions
when the branching rate satisfies $\al>4$, two solutions for $\al=4$, and only
one solution for $\al<4$. For $\al>4$, the solutions to the RDE are $\nu_{\rm
  low},\nu_{\rm mid}$, and $\nu_{\rm upp}$, where we let $\nu_{\rm\ldots}$
denote the probability measure on $\{0,1\}$ with mean
$\nu_{\rm\ldots}(\{1\})=z_{\rm\ldots}$ $(\ldots={\rm low},{\rm mid},{\rm
  upp})$ as defined around (\ref{fixpts}). By Lemma~\ref{L:RTP}, each of these
solutions to the RDE defines an RTP.

\subsection{Endogeny and bivariate uniqueness}\label{S:end}

In \cite[Def~7]{AB05}, an RTP $(\omb_\ibf,X_\ibf)_{\ibf\in\ovT}$ corresponding
to a solution $\nu$ of the RDE (\ref{RDE}) is called \emph{endogenous} if
$X_\wurz$ is a.s.\ measurable w.r.t.\ the \si-field generated by the random
variables $(\omb_\ibf)_{\ibf\in\ovT}$. In Lemma~\ref{L:endT} below, we will
show that this is equivalent to $X_\wurz$ being a.s.\ measurable w.r.t.\ the
\si-field generated by the random variables $\St$ and
$(\omb_\ibf)_{\ibf\in\St}$, where $\St$ is the random tree defined in
(\ref{Sdef}). Aldous and Bandyopadhyay have shown that endogeny is equivalent
to \emph{bivariate uniqueness}, which we now explain.

Let $\Pc_{\rm sym}(S^n)$ denote the space of probability measures on $S^n$
that are symmetric with respect to permutations of the coordinates. Let
$\pi_m:S^n\to S$ denote the projection on the $m$-th coordinate, i.e.,
$\pi_m(x^1,\ldots,x^n):=x_m$, and let $\mu^{(n)}\circ\pi_m^{-1}$ denote the
$m$-th marginal of a measure $\mu^{(n)}\in\Pc(S^n)$. For any $\mu\in\Pc(S)$,
we define
\be
\Pc(S^n)_\mu:=\big\{\mu^{(n)}\in\Pc(S^n):
\mu^{(n)}\circ\pi_m^{-1}=\mu\ \forall 1\leq m\leq n\big\}
\ee
to be the set of probability measures on $S^n$ whose one-dimensional marginals
are all equal to $\mu$, and we denote
$\Pc_{\rm{sym}}(S^n)_\mu:=\Pc_{\rm{sym}}(S^n)\cap\Pc(S^n)_\mu$.
Finally, we define a ``diagonal'' set
\be\label{diag}
S^n_{\rm diag}:=\big\{(x_1,\ldots,x_n)\in S^n:x_1=\cdots=x_n\big\}
\ee
and given a measure $\mu\in\Pc(S)$, we let $\ov\mu^{(n)}$ denote the unique
element of $\Pc(S^n)_{\mu}\cap\Pc(S^n_{\rm diag})$, i.e.,
\be\label{nmu}
\ov\mu^{(n)}:=\P\big[(X,\ldots,X)\in\,\cdot\,\big],
\textrm{ where $X$ has law $\mu$}.
\ee
Recall the definition of the $n$-variate map $T^{(n)}$ in (\ref{Tndef}).
The following theorem has been proved in \cite[Thm~1]{MSS18}, and in a
slightly weaker form in \cite[Thm~11]{AB05}. Below, $\Rightarrow$ denotes weak
convergence of probability measures.

\bt[Endogeny and $n$-variate uniqueness]
Let\label{T:bivar} $\nu$ be a solution of the RDE (\ref{RDE}). Then the
following statements are equivalent.
\begin{enumerate}
\item The RTP corresponding to $\nu$ is endogenous.
\item $\dis(T^{(n)})^m(\mu)\Asto{m}\ov\nu^{(n)}$ for all
  $\mu\in\Pc(S^n)_\nu$ and $n\geq 1$.
\item $\ov\nu^{(2)}$ is the only fixed point of $T^{(2)}$ in the space
  $\Pc_{\rm{sym}}(S^2)_\nu$.
\end{enumerate}
\et
We remark that bivariate uniqueness as introduced in \cite{AB05} refers to
$\ov\nu^{(2)}$ being the only fixed point of $T^{(2)}$ in the space
$\Pc(S^2)_\nu$. The equivalences in the above theorem tells us that bivariate
uniqueness already follows from the weaker condition (iii) since it implies
(ii), which implies n-variate uniqueness for any $n\geq 1$.

We will prove a continuous-time extension of
Theorem~\ref{T:bivar}, relating endogeny to solutions of the 
\emph{$n$-variate mean-field equation}
\be\label{nvariate}
\dif{t}\mu^{(n)}_t=|\rb|\big\{T^{(n)}(\mu^{(n)}_t)-\mu^{(n)}_t\big\}
\qquad(t\geq 0),
\ee
where we have replaced $T$ in (\ref{mean}) by $T^{(n)}$ and we write
$\mu^{(n)}_t$ to remind ourselves that this is a measure on $S^n$, rather than
on $S$.

This equation has the following interpretation. As
in Subsection~\ref{S:flow}, let $(\Xb_{s,u})_{s\leq u}$ be a stochastic flow
on $S^N$ constructed from a Poisson point set $\Pi$. Let
$(X^1(0),\ldots,X^n(0))$ be a random variable with values in $S^n$,
independent of $(\Xb_{s,u})_{s\leq u}$. Then setting
\be
\big(X^1(t),\ldots,X^n(t)\big)
:=\big(\Xb_{0,t}(X^1(0)),\ldots,\Xb_{0,t}(X^n(0))\big)
\qquad(t\geq 0)
\ee
defines a Markov process $(X^1(t),\ldots,X^n(t))_{t\geq 0}$ that consists of
$n$ Markov processes with initial states $X^1(0),\ldots,X^n(0)$ that are
coupled in such a way that they are constructed using the same stochastic
flow. Applying Theorem~\ref{T:meanlim} to this \emph{$n$-variate Markov
  process}, we see that the mean-field equation for the $n$-variate process
takes the form (\ref{nvariate}).

We note that if $\mu^{(n)}_t$ solves the $n$-variate mean-field equation, then
any $m$-dimensional marginal of $\mu^{(n)}_t$ solves the $m$-variate
mean-field equation.  Also, solutions to (\ref{nvariate}) started in an
initial condition $\mu^{(n)}_0\in\Pc_{\rm sym}(S^n)$ satisfy
$\mu^{(n)}_t\in\Pc_{\rm sym}(S^n)$ for all $t\geq 0$. Finally, it is easy to
see that $\mu^{(n)}_0\in\Pc(S^n_{\rm diag})$ implies
$\mu^{(n)}_t\in\Pc(S^n_{\rm diag})$ for all $t\geq 0$.

We now formulate a continuous-time extension of
Theorem~\ref{T:bivar}. Note that in view of (\ref{RDE}), a measure $\nu^{(2)}$
is a fixed point of the \emph{bivariate mean-field equation} (i.e.,
(\ref{nvariate}) with $n=2$) if and only if it is a fixed point of
$T^{(2)}$. Therefore, the equivalence of points~(i) and (iii) from
Theorem~\ref{T:bivar} immediately implies an analogue statement in the
continuous-time setting.

\bt[Endogeny and the n-variate mean-field equation]
Under\label{T:bivar2} the assumptions of Theorem~\ref{T:bivar},
the following conditions are equivalent.
\begin{enumerate}
\item The RTP corresponding to $\nu$ is endogenous.
\item For any $\mu^{(n)}_0\in\Pc(S^n)_\nu$ and $n\geq 1$, the solution
  $(\mu^{(n)}_t)_{t\geq 0}$ to the $n$-variate equation (\ref{nvariate})
  started in $\mu^{(n)}_0$ satisfies $\dis\mu^{(n)}_t\Asto{t}\ov{\nu}^{(n)}$.
\end{enumerate}
\et

Theorem~\ref{T:bivar2} motivates us to study the bivariate mean-field equation
in our example of a particle system with cooperative branching. Recall
  that in this example, $\Gi:=\{\cob,\dth\}$ with $\cob$ and $\dth$ as in
\eqref{cobdth}, and $\pi$ is defined in (\ref{picob}). In line with
(\ref{coopmu11}) we write the bivariate mean-field equation as
\be\label{pibivar}
\dif{t}\mu^{(2)}_t=\al\big\{T_{\cob^{(2)}}(\mu^{(2)}_t)-\mu^{(2)}_t\big\}
+\big\{T_{\dth^{(2)}}(\mu^{(2)}_t)-\mu^{(2)}_t\big\}.
\ee
For simplicity, we restrict ourselves to symmetric solutions, i.e., solutions
that take values in $\Pc_{\rm sym}(\{0,1\}^2)$. For any probability measure
$\mu^{(2)}\in\Pc_{\rm sym}(\{0,1\}^2)$, we let $\mu^{(1)}$ denote its
one-dimensional marginals, which are equal by symmetry. We let
$\nu_{\rm low},\nu_{\rm mid},\nu_{\rm upp}$ denote the probability measures on
$\{0,1\}$ with mean $\nu_{\rm\ldots}(\{1\})=z_{\rm\ldots}$
$(\ldots={\rm low},{\rm mid},{\rm upp})$ as defined around (\ref{fixpts}).

\bp[Bivariate equation for cooperative branching]
For\label{P:copbiv} $\al>4$, the bivariate mean-field equation (\ref{pibivar})
has precisely four fixed points in $\Pc_{\rm sym}(\{0,1\}^2)$, namely
\be\label{bifix}
\ov\nu^{(2)}_{\rm low},\quad\un\nu^{(2)}_{\rm mid},
\quad\ov\nu^{(2)}_{\rm mid},\quand\ov\nu^{(2)}_{\rm upp}.
\ee
which are uniquely characterized by their respective marginals
$\nu_{\rm low},\nu_{\rm mid},\nu_{\rm mid},\nu_{\rm upp}$,
as well as the fact that $\ov\nu^{(2)}_{\rm low},\ov\nu^{(2)}_{\rm mid}$, and
$\ov\nu^{(2)}_{\rm upp}$ are concentrated on
$\{0,1\}^2_{\rm diag}=\{(0,0),(1,1)\}$, but $\un\nu^{(2)}_{\rm mid}$ is not.

For any $\mu^{(2)}_0\in\Pc_{\rm sym}(\{0,1\}^2)$, the solution
to (\ref{pibivar}) started in $\mu^{(2)}_0$ 
converges as $t\to\infty$ to one of the fixed points in (\ref{bifix}), the
respective domains of attraction being
\be\ba{l}
\dis\big\{\mu^{(2)}_0:\mu^{(1)}_0(\{1\})<z_{\rm mid}\big\},\quad
\big\{\mu^{(2)}_0:\mu^{(1)}_0(\{1\})=z_{\rm mid},
\ \mu^{(2)}_0\neq\ov\nu^{(2)}_{\rm mid}\big\},\\[5pt]
\dis\big\{\ov\nu^{(2)}_{\rm mid}\big\},\quand
\big\{\mu^{(2)}_0:\mu^{(1)}_0(\{1\})>z_{\rm mid}\big\}.
\ec
For $\al=4$, there are two fixed points $\ov\nu^{(2)}_{\rm low}$ and
$\ov\nu^{(2)}_{\rm upp}$ with respective domains of attraction
\be
\big\{\mu^{(2)}_0:\mu^{(1)}_0(\{1\})<z_{\rm mid}\big\}\quand
\big\{\mu^{(2)}_0:\mu^{(1)}_0(\{1\})\geq z_{\rm mid}\big\},
\ee
while for $\al<4$ all solutions converge to $\ov\nu^{(2)}_{\rm low}$.
\ep

Combining Proposition~\ref{P:copbiv} with Theorem~\ref{T:bivar2}, we see
that the RTPs corresponding to $\nu_{\rm low}$ and $\nu_{\rm upp}$ are
endogenous, but for $\al>4$, the RTP corresponding to $\nu_{\rm mid}$ is
not. As is clear from \cite[Table~1]{AB05}, few examples of nonendogenous RTPs
were known at the time. Contrary to what is stated in \cite[Table~1]{AB05},
frozen percolation is now generally conjectured to be nonendogenous, but until
recently few ``natural'' examples of nonendogenous RTPs have appeared in the
literature. In fact, the RTP corresponding to $\nu_{\rm mid}$ seems
to be one of the simplest nontrivial examples of a nonendogenous RTP
discovered so far. Another nice class of nonendogenous RTPs has recently been
described in \cite{MS18}.

\subsection{The higher-level mean-field equation}\label{S:high}

Following \cite[formula~(1.1)]{MSS18}, if $S$ is a Polish space and
$g:S^k\to S$ is a measurable map, then we define a measurable map
$\ch g:\Pc(S)^k\to\Pc(S)$ by
\be\bac
\dis\ch g&:=&\dis\mbox{ the law of }g(X_1,\ldots,X_k),\\
&&\dis\mbox{where }(X_1,\ldots,X_k)\mbox{ are independent with laws }
\mu_1,\ldots,\mu_k.
\ec
Note that in this notation, the map $T_g:\Pc(S)\to\Pc(S)$ from (\ref{Tgdef})
is given by $T_g(\mu)=\ch g(\mu,\ldots,\mu)$. As in
\cite[formula~(4.2)]{MSS18}, we define a \emph{higher-level map}
$\ch T :\Pc(\Pc(S))\to\Pc(\Pc(S))$ by
\be\label{chTdef}
\ch T(\rho):=\mbox{ the law of }\ch\ga[\omb](\xi_1,\xi_2,\ldots),
\ee
where $\omb$ is an $\Om$-valued random variable with law $|\rb|^{-1}\rb$ and
$(\xi_i)_{i\geq 1}$ are i.i.d.\ $\Pc(S)$-valued random variables with law $\rho$.
Iterates of the map $\ch T$ have been studied in \cite[Section~4]{MSS18}.
We will be interested in the \emph{higher-level mean-field equation}
\be\label{levelmean}
\dif{t}\rho_t=|\rb|\big\{\ch T(\rho_t)-\rho_t\big\}\qquad(t\geq 0).
\ee
A measure $\rho\in\Pc(\Pc(S))$ is the law of a random probability measure
$\xi$ on $S$. We denote the \emph{$n$-th moment measure} of such a random
measure $\xi$ by
\be
\rho^{(n)}:=\E\big[\xi\otimes\cdots\otimes\xi]
\quad\mbox{where }\xi\mbox{ has law }\rho.
\ee
(Here $\E[\,\cdot\,]$ denotes the expectation of a random measure; see the
remark above Theorem~\ref{T:meanlim}.) Our notation for moment measures is on
purpose similar to our earlier notation for solutions to the $n$-variate
equation, because of the following proposition.

\bp[Moment measures]
If\label{P:levmom} $(\rho_t)_{t\geq 0}$ solves the higher-level mean-field
equation (\ref{levelmean}), then its $n$-th moment measures
$(\rho^{(n)}_t)_{t\geq 0}$ solve the $n$-variate equation (\ref{nvariate}).
\ep

Similarly to Proposition~\ref{P:levmom}, it has been shown in
\cite[Lemma~2]{MSS18} that $\ch T(\rho)^{(n)}=T^{(n)}(\rho^{(n)})$, and this
formula holds even for $n=\infty$. In view of this, as discussed in
Subsection~\ref{S:firstintro}, the higher-level map $\ch T$ is effectively
equivalent to the $\infty$-variate map $\T^{(\infty)}:\Pc_{\rm
  sym}(S^\infty)\to\Pc_{\rm sym}(S^\infty)$. It follows from
Proposition~\ref{P:levmom} that if $\rho$ solves the \emph{higher-level RDE}
\be\label{RDElevel}
\ch T(\rho)=\rho,
\ee
then its $n$-th moment measures solve the $n$-variate RDE
$T^{(n)}(\rho^{(n)})=\rho^{(n)}$, with $T^{(n)}$ as in (\ref{Tndef}).

If $X$ is an $S$-valued random variable defined on some probability space
$(\Om,\Fi,\P)$ and $\Hi\sub\Fi$ is a sub-\si-field, then
$\P[X\in\,\cdot\,|\Hi]$ is a random probability measure\footnote{Here we use
  that since $S$ is Polish, regular versions of conditional expectations
  exist.} on $S$. As a consequence, the law of $\P[X\in\,\cdot\,|\Hi]$ is an
element of $\Pc(\Pc(S))$. In the following theorem, which is based on
\cite[Thm~2]{Str65} and which in its present form we cite from
\cite[Thm~13]{MSS18}, we use the fact that each Polish space $S$ has a
metrizable compactification $\ov S$ \cite[\Parag 6 No.~1,
  Theorem~1]{Bou58}. Moreover, we naturally identify $\Pc(S)$ with the space
of all probability measures on $\ov S$ that are concentrated on $S$.

\bt[The convex order for laws of random probability measures]
Let\label{T:Stras} $S$ be a Polish space, let $\ov S$ be a metrizable
compactification of $S$, and let $\Ci_{\rm cv}\big(\Pc(\ov S)\big)$ denote the
space of all convex continuous functions $\phi:\Pc(\ov S)\to\R$.
Then, for $\rho_1,\rho_2\in\Pc(\Pc(S))$, the following
statements are equivalent.
\begin{enumerate}
\item $\dis\int\phi\,\di\rho_1\leq\int\phi\,\di\rho_2$ for all
  $\phi\in\Ci_{\rm cv}\big(\Pc(\ov S)\big)$.
\item There exists an $S$-valued random variable $X$ defined on some
  probability space $(\Om,\Fi,\P)$ and sub-\si-fields $\Hi_1\sub\Hi_2\sub\Fi$
  such that $\dis\rho_i=\P\big[\P[X\in\,\cdot\,|\Hi_i]\in\,\cdot\,\big]$
  $(i=1,2)$.
\end{enumerate}
\et

If $\rho_1,\rho_2\in\Pc(\Pc(S))$ satisfy the equivalent conditions of
Theorem~\ref{T:Stras}, then we say that they are ordered in the \emph{convex
  order} and denote this as $\rho_1\leq_{\rm cv}\rho_2$. It follows from
\cite[Lemma~15]{MSS18} that $\leq_{\rm cv}$ is a partial order; in particular,
$\rho_1\leq_{\rm cv}\rho_2$ and $\rho_2\leq_{\rm cv}\rho_1$ imply $\rho_1=\rho_2$.

Recall that in Subsection~\ref{S:firstintro}, we defined $\Pc(\Pc(S))_\mu$,
which is $\big\{\rho\in\Pc(\Pc(S)):\rho^{(1)}=\mu\big\}$. We define
$\ov\mu\in\Pc(\Pc(S))_\mu$ by $\ov\mu:=\P[\de_X\in\,\cdot\,]$, where $X$ has
law $\mu$. It is easy to see that the $n$-th moment measures of $\ov\mu$ are
given by (\ref{nmu}), so our present notation is consistent with
earlier notation introduced there. By \cite[formula (4.7)]{MSS18}, the
measures $\de_\mu,\ov\mu$ are the extremal elements of $\Pc(\Pc(S))_\mu$
w.r.t.\ the convex order, i.e.,
\be\label{cvextr}
\de_\mu\leq_{\rm cv}\rho\leq_{\rm cv}\ov\mu\qquad\big(\rho\in\Pc(\Pc(S))_\mu\big).
\ee
The following proposition is a continuous-time version of
\cite[Prop~3]{MSS18}.

\bp[Extremal solutions in the convex order]
If\label{P:minconv} $(\rho^i_t)_{t\geq 0}$ $(i=1,2)$ are solutions to the
higher-level mean-field equation (\ref{levelmean}) such that
$\rho^1_0\leq_{\rm cv}\rho^2_0$, then $\rho^1_t\leq_{\rm cv}\rho^2_t$ for all
$t\geq 0$. If $\nu$ solves the RDE (\ref{RDE}), then $\ov\nu$ solves the
higher-level RDE (\ref{RDElevel}) and there exists a solution $\un\nu$ of
(\ref{RDElevel}) such that
\be\label{minconv}
\rho_t\Asto{t}\un\nu
\qquad\mbox{where $(\rho_t)_{t\geq 0}$ solves (\ref{levelmean}) with }
\rho_0=\de_\nu.
\ee
Here $\Rightarrow$ denotes weak convergence of measures on $\Pc(S)$, equipped
with the topology of weak convergence. Any solution $\rho\in\Pc(\Pc(S))_\nu$
to the higher-level RDE (\ref{RDElevel}) satisfies
\be\label{sandwich}
\un\nu\leq_{\rm cv}\rho\leq_{\rm cv}\ov\nu.
\ee
\ep

The following result, which we cite from \cite[Prop.~4]{MSS18}, describes the
higher-level RTPs associated with the solutions $\un\nu$ and $\ov\nu$ of the
higher-level RDE.

\bp[Higher-level RTPs]
Let\label{P:hlRTP} $\nu$ be a solution of the RDE (\ref{RDE}) and let $\un\nu$
and $\ov\nu$ as in (\ref{sandwich}) be the corresponding minimal and maximal
solutions to the higher-level RDE, with respect to the convex order. Let
$(\omb_\ibf,X_\ibf)_{\ibf\in\ovT}$ be an RTP corresponding to $\ga$ and $\nu$
and set
\be\label{xidef}
\xi_\ibf:=\P\big[X_\ibf\in\,\cdot\,|\,(\omb_{\ibf\jbf})_{\jbf\in\ovT}\big]
\qquad(\ibf\in\T).
\ee
Then $(\omb_\ibf,\xi_\ibf)_{\ibf\in\ovT}$ is an RTP corresponding to $\ch\ga$ and
$\un\nu$. Also, $(\omb_\ibf,\de_{X_\ibf})_{\ibf\in\ovT}$ is an RTP
corresponding to $\ch\ga$ and $\ov\nu$.
\ep

Proposition~\ref{P:hlRTP} gives a more concrete interpretation of the
solutions $\un\nu$ and $\ov\nu$ to the higher-level RDE from (\ref{sandwich}).
Indeed, if $(\omb_\ibf,X_\ibf)_{\ibf\in\ovT}$ is an RTP corresponding to $\nu$, then
\be\label{ovnu}
\ov\nu=\P\big[\de_{X_\wurz}\in\,\cdot\,\big],
\ee
which corresponds to ``perfect knowledge'' about the state $X_\wurz$ of the
root, while
\be\label{unnu}
\un\nu=\P\big[\P\big[X_\wurz\in\,\cdot\,|\,(\omb_\ibf)_{\ibf\in\ovT}\big]
\in\,\cdot\,\big]
\ee
corresponds to the knowledge about $X_\wurz$ that is contained in the random
variables $(\omb_\ibf)_{\ibf\in\ovT}$. Since $X_\wurz$ is a measurable
function of $(\om_\ibf)_{\ibf\in\T}$ if and only if its conditional law given
$(\omb_\ibf)_{\ibf\in\ovT}$ equals $\de_{X_\wurz}$, it follows from
(\ref{ovnu}) and (\ref{unnu}) that the RTP corresponding to $\nu$ is
endogenous if and only if $\un\nu=\ov\nu$.

It is instructive to demonstrate the general theory on our concrete example of
a system with cooperative branching and deaths. Recall that for $\al>4$, the
mean-field equation (\ref{coopmu11}) has three fixed points $\nu_{\rm
  low},\nu_{\rm mid},\nu_{\rm upp}$. We denote the corresponding minimal and
maximal solutions to the higher-level RDE in the sense of (\ref{sandwich}) by
$\un\nu_{\rm\ldots}$ and $\ov\nu_{\rm\ldots}$ $(\ldots={\rm low},{\rm
  mid},{\rm upp})$. The following theorem lifts the results from
Proposition~\ref{P:copbiv} about the bivariate equation to a higher level.
Indeed, using the theorem below, it is easy to see that the measures
$\ov\nu^{(2)}_{\rm low},\un\nu^{(2)}_{\rm mid},\ov\nu^{(2)}_{\rm mid}$ and
$\ov\nu^{(2)}_{\rm upp}$ from Proposition~\ref{P:copbiv} are in fact the
second moment measures of the measures $\ov\nu_{\rm low},\un\nu_{\rm
  mid},\ov\nu_{\rm mid}$ and $\ov\nu_{\rm upp}$.

\bt[Higher-level equation for cooperative branching]
Let\label{T:coblev} $\nu_{\rm low},\nu_{\rm mid}$, and $\nu_{\rm upp}$ denote
the fixed points of the mean-field equation (\ref{coopmu11}) defined above
Proposition~\ref{P:copbiv}.
Then we have for the corresponding minimal and maximal solutions to the higher-level RDE that
\be
\un\nu_{\rm low}=\ov\nu_{\rm low},\quad
\un\nu_{\rm upp}=\ov\nu_{\rm upp},\quad\mbox{but}\quad
\un\nu_{\rm mid}\neq\ov\nu_{\rm mid}\quad(\al>4).
\ee
For $\al>4$, the higher-level RDE (\ref{RDElevel}) has four solutions, namely
\be\label{levfix}
\ov\nu_{\rm low},\quad\un\nu_{\rm mid},
\quad\ov\nu_{\rm mid},\quand\ov\nu_{\rm upp}.
\ee
Any solution $(\rho_t)_{t\geq 0}$ to the higher-level mean-field equation
(\ref{levelmean}) converges as $t\to\infty$ to one of the fixed points in
(\ref{levfix}), the respective domains of attraction being
\be\ba{l}\label{domains}
\dis\big\{\rho_0:\rho^{(1)}_0(\{1\})<z_{\rm mid}\big\},\quad
\big\{\rho_0:\rho^{(1)}_0(\{1\})=z_{\rm mid},
\ \rho_0\neq\ov\nu_{\rm mid}\big\},\\[5pt]
\dis\big\{\ov\nu_{\rm mid}\big\},\quand
\big\{\rho_0:\rho^{(1)}_0(\{1\})>z_{\rm mid}\big\}.
\ec
For $\al=4$, there are two fixed points $\ov\nu_{\rm low}$ and
$\ov\nu_{\rm upp}$ with respective domains of attraction
\be
\big\{\rho_0:\rho^{(1)}_0(\{1\})<z_{\rm mid}\big\}\quand
\big\{\rho_0:\rho^{(1)}_0(\{1\})\geq z_{\rm mid}\big\},
\ee
while for $\al<4$ all solutions converge to $\ov\nu_{\rm low}$.
\et

Since a probability measure $\mu\in\Pc(\{0,1\})$ is uniquely characterized by
$\mu(\{1\})\in[0,1]$, there is a natural identification
$\Pc(\{0,1\})\cong[0,1]$. Let $\widehat\cob$ and $\widehat\dth$ denote the
higher-level maps $\ch g$ corresponding to $g=\cob,\dth$, which using the
identification $\Pc(\{0,1\})\cong[0,1]$ we view as maps
$\widehat\cob:[0,1]^3\to[0,1]$ and $\widehat\dth:[0,1]^0\to[0,1]$. One can
check that
\be\label{coblev}
\widehat\dth(\wurz)=0\quand
\widehat\cob(\eta_1,\eta_2,\eta_3)=\eta_1+(1-\eta_1)\eta_2\eta_3
\qquad\big(\eta_1,\eta_2,\eta_3\in[0,1]\big).
\ee
Identifying $\Pc(\Pc(\{0,1\}))\cong\Pc[0,1]$, we can identify the measures
$\ov\nu_{\rm low},\un\nu_{\rm mid},\ov\nu_{\rm mid}$, and $\ov\nu_{\rm upp}$
with probability laws on $[0,1]$. Letting $\eta$ denote a random variable with
law $\nu\in\Pc[0,1]$, the higher-level RDE, written in the form (\ref{isd}),
then reads
\be\label{ourRDE}
\eta\isd\chi\cdot\big(\eta_1+(1-\eta_1)\eta_2\eta_3\big),
\ee
where $\eta_1,\eta_2,\eta_3$ are independent copies of $\eta$ and $\chi$ is an
independent Bernoulli random variable with $\P[\chi=1]=\al/(\al+1)$.
Theorem~\ref{T:coblev} says that for $\al>4$, this equation has four
solutions. Three ``trivial'' solutions $\ov\nu_{\rm low},\ov\nu_{\rm
  mid},\ov\nu_{\rm upp}$ that correspond to Bernoulli $\eta$ with parameters
\be
\P[\eta=1]=z,\quad\P[\eta=0]=1-z\quad\mbox{with}\quad
z=z_{\rm low},z_{\rm mid},z_{\rm low},
\ee
and a ``nontrivial'' solution $\un\nu_{\rm mid}$ for which $\P[0<\eta<1]>0$.
In view of Proposition~\ref{P:hlRTP}, we can interpret this nontrivial
solution (viewed as a probability law on $[0,1]$) as
\be
\un\nu_{\rm mid}
=\P\big[\P[X_\wurz=1\,|\,(\omb_\ibf)_{\ibf\in\ovT}\big]\in\,\cdot\,\big],
\ee
where $(\omb_\ibf,X_\ibf)_{\ibf\in\ovT}$ is the RTP corresponding to $\nu_{\rm
  mid}$. The following lemma summarizes some elementary facts about the law
$\un\nu_{\rm mid}$. We note that by solving the $n$-variate RDE for $n\geq 3$,
one should in principle be able to calculate higher moments of $\un\nu_{\rm
  mid}$, although the formulas quickly become unwieldy.

\bl[Nontrivial solution of the higher-level RDE]
Let\label{L:nontriv} $\al>4$ and let $\eta$ be a random variable with law
$\un\nu_{\rm mid}$. Then
\bc
\dis\E[\eta]&=&\dis z_{\rm mid}
=\ffrac{1}{2}-\sqrt{\ffrac{1}{4}-\ffrac{1}{\al}},\\[5pt]
\dis\E[\eta^2]&=&\dis\ha z_{\rm mid}-\ffrac{1}{2}
+\ha\sqrt{13z_{\rm mid}^2-6z_{\rm mid}+1+\ffrac{4}{\al}}.
\ec
Moreover,
\be
\dis\P[\eta=0]=z_{\rm mid}\quand\P[\eta=1]=0.
\ee
\el

\begin{figure}
\centering
\begin{tikzpicture}[>=triangle 45]

\begin{scope}[scale=3.3]
\draw[->] (-0.1,0) -- (1.2,0);
\draw[->] (0,-0.1) -- (0,1.2);
\draw (1.1,0) node[above] {$\eta$};
\draw (0,1) node[right] {$F(\eta)$};
\foreach \x in {0.2,0.4,0.6,0.8,1}
 \draw[thick] (\x,-0.02) -- (\x,0.02);
\foreach \x in {0.2,0.4,0.6,0.8,1}
 \draw (\x,-0.02) node[below] {\x};
\foreach \y in {0.2,0.4,0.6,0.8,1}
 \draw[thick] (-0.02,\y) -- (0.02,\y);
\foreach \y in {0.2,0.4,0.6,0.8,1}
 \draw (-0.02,\y) node[left] {\y};

\draw[very thick] plot file {plot41.dat};
\draw (0.65,0.35) node {$\al=4.1$};
\end{scope}

\begin{scope}[scale=3.3,xshift=1.5cm]
\draw[->] (-0.1,0) -- (1.2,0);
\draw[->] (0,-0.1) -- (0,1.2);
\draw (1.1,0) node[above] {$\eta$};
\draw (0,1) node[right] {$F(\eta)$};
\foreach \x in {0.2,0.4,0.6,0.8,1}
 \draw[thick] (\x,-0.02) -- (\x,0.02);
\foreach \x in {0.2,0.4,0.6,0.8,1}
 \draw (\x,-0.02) node[below] {\x};
\foreach \y in {0.2,0.4,0.6,0.8,1}
 \draw[thick] (-0.02,\y) -- (0.02,\y);
\foreach \y in {0.2,0.4,0.6,0.8,1}
 \draw (-0.02,\y) node[left] {\y};

\draw[very thick] plot file {unif100.dat};
\draw (0.65,0.35) node {$\al=4.5$};
\end{scope}

\begin{scope}[scale=3.3,xshift=3cm]
\draw[->] (-0.1,0) -- (1.2,0);
\draw[->] (0,-0.1) -- (0,1.2);
\draw (1.1,0) node[above] {$\eta$};
\draw (0,1) node[right] {$F(\eta)$};
\foreach \x in {0.2,0.4,0.6,0.8,1}
 \draw[thick] (\x,-0.02) -- (\x,0.02);
\foreach \x in {0.2,0.4,0.6,0.8,1}
 \draw (\x,-0.02) node[below] {\x};
\foreach \y in {0.2,0.4,0.6,0.8,1}
 \draw[thick] (-0.02,\y) -- (0.02,\y);
\foreach \y in {0.2,0.4,0.6,0.8,1}
 \draw (-0.02,\y) node[left] {\y};

\draw[very thick] plot file {plot100.dat};
\draw (0.65,0.35) node {$\al=10$};
\end{scope}


\end{tikzpicture}

\caption{Numerical data for the distribution function $F(\eta):=\un\nu_{\rm
    mid}([0,\eta])$ of the nontrivial solution $\un\nu_{\rm mid}$ to the RDE
  (\ref{ourRDE}).}
\label{fig:numid}
\end{figure}
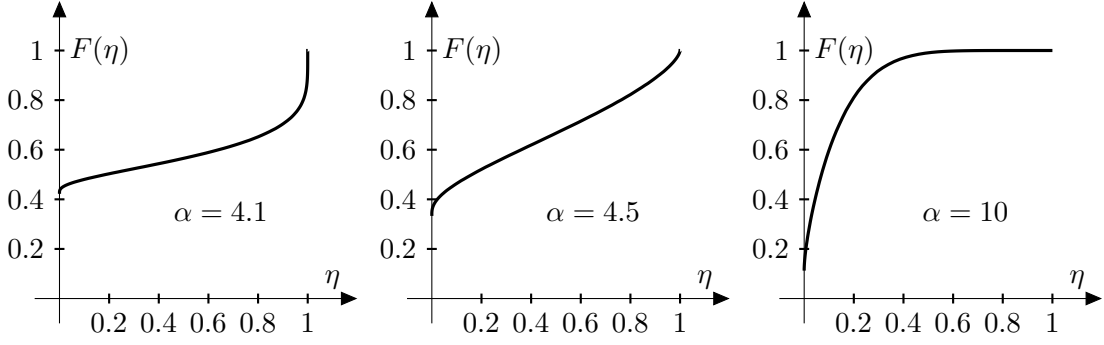

It is not too hard to obtain numerical data for $\un\nu_{\rm mid}$, see
Figure~\ref{fig:numid}. These data suggest that apart from the atom in $0$,
the measure $\un\nu_{\rm mid}$ has a smooth density with respect to the
Lebesgue measure, but we have no proof for this. We have tried to
find an explicit formula for the density but have not been successful.

\subsection{Lower and upper solutions}

In this and the next subsection we collect a few further results on endogeny
and the uniqueness of solutions to RDEs. In the present subsection, we show
that the endogeny of the RTPs corresponding to $\nu_{\rm low}$ and $\nu_{\rm
  upp}$ follows from a general principle, discovered in \cite{AB05}, that says
that RDEs that are defined by monotone maps always have a minimal and maximal
solution with respect to the stochastic order, and that the RTPs corresponding
to these solutions are always endogenous.

Let $S$ be a compact metrizable space that is equipped with a partial
order $\leq$ that is \emph{closed} in the sense that
\be\label{compat}
\big\{(x,y)\in S^2:x\leq y\big\}
\ee
is a closed subset of $S^2$, equipped with the product topology. Recall
that a function $f$ from one partially ordered space into another is monotone
if $x\leq y$ implies $f(x)\leq f(y)$, and a subset $A$ of a partially ordered
space is increasing if $A\ni x\leq y$ implies $y\in A$. It is known
that for two probability measures $\mu_1,\mu_2\in\Pc(S)$, the following
statements are equivalent:
\begin{enumerate}
\item $\mu_1(A)\leq\mu_2(A)$ for all closed increasing $A\sub S$.
\item $\dis\int f\,\di\mu_1\leq\int f\,\di\mu_2$ for all bounded continuous
  monotone $f:S\to\R$.
\item Two random variables $X_1,X_2$ with laws $\mu_1,\mu_2$ can be coupled
  such that $X_1\leq X_2$ a.s.
\end{enumerate}
The equivalence of (ii) and (iii) is proved in \cite[Thm~II.2.4]{Lig85}.
The equivalence of (i) and (iii) holds more generally for Polish spaces, see
\cite[Thm~1 (ii) and (vi)]{KKO77}. In the general setting of Polish spaces,
the implications (iii)$\volgt$(i) and (iii)$\volgt$(ii) are trivial, but
the implication (ii)$\volgt$(i) needs the additional assumption of
\emph{monotone normality}, see \cite[Prop.~3.6 and 3.11]{HLL18}.

If $\mu_1,\mu_2$ satisfy the above conditions, then one says that they are
\emph{stochastically ordered}, denoted as $\mu_1\leq\mu_2$. This
defines a partial order on $\Pc(S)$; in particular, by
  Lemma~\ref{L:compar} below, $\mu_1\leq\mu_2\leq\mu_1$ implies $\mu_1=\mu_2$.

The proposition below is a variant of \cite[Lemma~15]{AB05}. As in our usual
setting, we assume that $S$ and $\Om$ are Polish spaces, $\rb$ is a nonzero
finite measure on $\Om$, and $\ga:\Om\times S^{\N_+}\to S$ and
$\kappa:\Om\to\N$ are measurable functions such that (\ref{kappa}) and
(\ref{sumr}) hold. If $S$ is equipped with a partial order, then we equip
$S^k$ with the product partial order. Recall that Proposition~\ref{P:gacon}
gives sufficient conditions for $T$ to be continuous w.r.t.\ the topology of
weak convergence.

\bp[Lower and upper solutions to RDE]
Assume\label{P:lowup} that $S$ is compact and equipped with a closed
partial order. Assume that $S$ has minimal and maximal
elements, denoted by $0$ and $1$. Assume $\ga[\om]$ is monotone for each
$\om\in\Om$ and that the operator $T$ in (\ref{Tdef}) is continuous
w.r.t.\ the topology of weak convergence. Then there exists solutions
$\nu_{\rm low},\nu_{\rm upp}$ to the RDE (\ref{RDE}) that are minimal and
maximal with respect to the stochastic order, in the sense that any solution
$\nu$ to the RDE (\ref{RDE}) must satisfy
\be
\nu_{\rm low}\leq\nu\leq\nu_{\rm upp},
\ee
where $\Rightarrow$ denotes weak convergence.
Moreover, if $(\mu^{\rm low}_t)_{t\geq 0}$ and $(\mu^{\rm upp}_t)_{t\geq 0}$
denote the solutions to the mean-field equation (\ref{mean}) with
initial states $\mu^{\rm low}_0=\de_0$ and $\mu^{\rm upp}_0=\de_1$, then
\be
\mu^{\rm low}_t\Asto{t}\nu_{\rm low}
\quand
\mu^{\rm upp}_t\Asto{t}\nu_{\rm upp}.
\ee
Finally, the RTPs corresponding to $\nu_{\rm low}$ and $\nu_{\rm upp}$ are
endogenous.
\ep

We can view the solutions $\nu_{\rm low}$ and $\nu_{\rm upp}$ to the RDE
(\ref{RDE}) as mean-field versions of the lower and upper invariant laws of
monotone particle systems; compare \cite[Thm~III.2.3]{Lig85}.

In our example of a system with cooperative branching, the maps $\cob$ and
$\dth$ are monotone, so Proposition~\ref{P:lowup} is applicable. Since the
measures we called $\nu_{\rm low}$ and $\nu_{\rm upp}$ before are the $t\to\infty$
limits of the solutions of the mean-field equation started in $\de_0$ and
$\de_1$, our earlier notation agrees with the more general notation of
Proposition~\ref{P:lowup}. The endogeny of the RTPs corresponding to
$\nu_{\rm low}$ and $\nu_{\rm upp}$, which before we proved based on an
analysis of the bivariate equation, using Proposition~\ref{P:copbiv} and
Theorem~\ref{T:bivar}, alternatively follows from Proposition~\ref{P:lowup}.

\subsection{Conditions for uniqueness}

In the present subsection, we prove some results of varying generality that
allow one to conclude that a given RDE has a unique solution. In our example
of a system with cooperative branching and deaths, this happens if and only if
$\al<4$. We will see that there are some general results that can be applied
to prove uniqueness in the whole regime $\al<4$. We also make a connection
with a general duality for monotone particle systems described in \cite{SS18}.
Although duality plays only a minor role in our paper, the original motivation
for the work that led to it was to understand this duality in the mean-field
limit.

We return to our usual set-up   from 
Subsection \ref{S:firstintro} with $S$ and $\Om$ Polish spaces 
and $\ga, \kappa$ and $\rb$ satisfying (\ref{kappa}) and (\ref{sumr}).
We also recall the random subtrees  $\St_t\sub\St\sub\ovT$  
defined in (\ref{Sdef}) as well as the fact that $\St_t$ for any $t \geq 0$
are a.s.\ finite by (\ref{sumr}). The tree
$\St$ is the family tree of the branching process $(\nab\St_t)_{t\geq 0}$. In
view of this, by well-known facts about branching processes, $\St$ is
a.s.\ finite if and only if
\be\label{Tfin}
\int_\Om\!\rb(\di\om)\,(\kappa(\om)-1)<0
\quad\mbox{or}\quad
\int_\Om\!\rb(\di\om)\,(\kappa(\om)-1)=0
\mbox{ while }\rb\big(\{\om:\kappa(\om)\neq 1\}\big)>0.
\ee
Recall that $G_t=G_{\St_t}$, where for any finite subtree $\U\sub\St$ that
contains the root, $G_\U:S^{\nab\U}\to S$ is the map defined in (\ref{GUdef}).
If $\St$ is a.s.\ finite, then $\nab\St_t=\emptyset$ for $t$ sufficiently large
and hence $G_t:S^{\nab\St_t}\to S$ is eventually constant.

More generally, if $\U$ is finite subtree of $\St$ that contains the root
$\wurz$, then we say that $\U$ is a \emph{root determining subtree} if the map
$G_\U:S^{\nab\U}\to S$ is constant. Note that this can happen even if
$\nab\U\neq\emptyset$. It is easy to see that if $\V\sub\U$ and $\V$ is root
determining, then the same is true for $\U$. We say that $\U$ is a
\emph{minimal root determining subtree} if $\U$ is root determining but there
exists no $\V\sub\U$ with $\V\neq\U$ that is root determining.
By our previous remark, it suffices to check this for such $\V$ that differ
from $\U$ by a single element.

\bl[Root determining subtrees]
The\label{L:rootdet} following conditions are equivalent:
\begin{enumerate}
\item There a.s.\ exists a $t<\infty$ such that $G_s$ is constant for all
  $s\geq t$.
\item $\St$ a.s.\ contains a root determining subtree.
\item $\St$ a.s.\ contains a minimal root determining subtree.
\end{enumerate}
\el

If $\U$ is a subtree of $\St$, then we denote by $\Xi_\U$ the set of all
$x=(x_\ibf)_{\ibf\in\U\cup\nab\U}$ that satisfy (\ref{xinduc}). We say that
$\U$ is \emph{uniquely determined} if $x,y\in\Xi_\U$ imply $x_\ibf=y_\ibf$
$(\ibf\in\U)$. The following lemma is inspired by \cite[Lemma~14]{AB05} who
showed that (i) implies that the RDE (\ref{RDE}) has a unique solution and the
corresponding RTP is endogenous.

\bl[Uniquely determined subtrees]
Between\label{L:detree} the following four conditions, one has the
implications (i)$\volgt$(ii)$\volgt$(iii)$\volgt$(iv) and (ii)$\volgt$(v). If
$S$ is finite, then moreover (iii)$\volgt$(ii), and if $S=\{0,1\}$, then
(ii)$\volgt$(i).
\begin{enumerate}
\item $\St$ a.s.\ contains a finite, uniquely determined subtree that contains
  the root $\wurz$.
\item The equivalent conditions of Lemma~\ref{L:rootdet} are satisfied.
\item $\St$ is a.s.\ uniquely determined.
\item The RDE (\ref{RDE}) has at most one solution and any corresponding RTP is
  endogenous.
\item The RDE (\ref{RDE}) has a solution $\nu$ that is globally attractive in
  the sense that any solution $(\mu_t)_{t\geq 0}$ to (\ref{mean})
  satisfies $\|\mu_t-\nu\|\asto{t}0$, where $\|\,\cdot\,\|$ denotes the total
  variation norm.
\end{enumerate}
\el

\begin{figure}
\centering

\begin{tikzpicture}[scale=0.6]

\node (wurz) at (0,2) {$\rm{cob}$};
\node (1) at (-9.6,4.5) {$\rm{dth}$};
\node (2) at (0,7) {$\rm{cob}$};
\node (3) at (9.6,5.5) {$\rm{cob}$};
\node (21) at (-3.6,9.5) {$\rm{cob}$};
\node[above] (22) at (0,15) {$X_{22}$};
\node[above] (23) at (3.6,15) {$X_{23}$};
\node (31) at (6,8.5) {$\rm{cob}$};
\node (32) at (9.6,7.5) {$\rm{cob}$};
\node (33) at (13.2,8) {$\rm{cob}$};
\node (211) at (-4.8,12) {$\rm{dth}$};
\node (212) at (-3.6,13.7) {$\rm{dth}$};
\node (213) at (-2.4,13) {$\rm{dth}$};
\node (311) at (4.8,14) {$\rm{dth}$};
\node (312) at (6,11.5) {$\rm{dth}$};
\node[above] (313) at (7.2,15) {$X_{313}$};
\node (321) at (8.4,9.5) {$\rm{dth}$};
\node[above] (322) at (9.6,15) {$X_{322}$};
\node[above] (323) at (10.8,15) {$X_{323}$};
\node (331) at (12,10.5) {$\rm{dth}$};
\node[above] (332) at (13.2,15) {$X_{332}$};
\node (333) at (14.4,13) {$\rm{dth}$};

\draw[line width=6pt] (0,0)--(wurz);
\draw[line width=2pt] (wurz) ellipse (0.7cm and 0.4cm);
\node[below] at (0,0) {$X_\wurz$};

\draw[line width=6pt,rounded corners] (-0.7,2)--(-9.6,2)--(1);
\draw[line width=2pt] (1) ellipse (0.7cm and 0.4cm);

\draw[line width=1pt] (0,2.5)--(2);
\draw[line width=2pt] (2) ellipse (0.7cm and 0.4cm);

\draw[line width=6pt,rounded corners] (0.7,2)--(9.6,2)--(3);
\draw[line width=2pt] (3) ellipse (0.7cm and 0.4cm);

\draw[line width=1pt,rounded corners] (-0.7,7)--(-3.6,7)--(21);
\draw[line width=2pt] (21) ellipse (0.7cm and 0.4cm);

\draw[line width=1pt] (0,7.5)--(22);

\draw[line width=1pt,rounded corners] (0.7,7)-- (3.6,7)--(23);

\draw[line width=6pt,rounded corners] (8.9,5.5)--(6,5.5)--(31);
\draw[line width=2pt] (31) ellipse (0.7cm and 0.4cm);

\draw[line width=1pt] (9.6,6)--(32);
\draw[line width=2pt] (32) ellipse (0.7cm and 0.4cm);

\draw[line width=6pt,rounded corners] (10.3,5.5)--(13.2,5.5)--(33);
\draw[line width=2pt] (33) ellipse (0.7cm and 0.4cm);

\draw[line width=1pt,rounded corners] (-4.3,9.5)--(-4.8,9.5)--(211);
\draw[line width=2pt] (211) ellipse (0.7cm and 0.4cm);

\draw[line width=1pt,rounded corners] (-3.6,10)--(212);
\draw[line width=2pt] (212) ellipse (0.7cm and 0.4cm);

\draw[line width=1pt,rounded corners] (-2.9,9.5)--(-2.4,9.5)--(213);
\draw[line width=2pt] (213) ellipse (0.7cm and 0.4cm);

\draw[line width=6pt,rounded corners] (5.3,8.5)--(4.8,8.5)--(311);
\draw[line width=2pt] (311) ellipse (0.7cm and 0.4cm);

\draw[line width=6pt,rounded corners] (6,9)--(312);
\draw[line width=2pt] (312) ellipse (0.7cm and 0.4cm);

\draw[line width=1pt,rounded corners] (6.7,8.5)--(7.2,8.5)--(313);

\draw[line width=1pt,rounded corners] (8.9,7.5)--(8.4,7.5)--(321);
\draw[line width=2pt] (321) ellipse (0.7cm and 0.4cm);

\draw[line width=1pt] (9.6,8)--(322);

\draw[line width=1pt,rounded corners] (10.3,7.5)--(10.8,7.5)--(323);

\draw[line width=6pt,rounded corners] (12.5,8)--(12,8)--(331);
\draw[line width=2pt] (331) ellipse (0.7cm and 0.4cm);

\draw[line width=1pt,rounded corners] (13.2,8.5)--(332);

\draw[line width=6pt,rounded corners] (13.9,8)--(14.4,8)--(333);
\draw[line width=2pt] (333) ellipse (0.7cm and 0.4cm);

\draw[line width=1pt, dotted](15,15)--(-10,15) node[above]{$t$};

\end{tikzpicture}
\caption{A minimal root determining subtree. In this example, $X_\wurz=0$
  regardless of the values of $X_{22},X_{23},X_{313},X_{322},X_{323},X_{332}$.}
\label{fig:rootdet}
\end{figure}
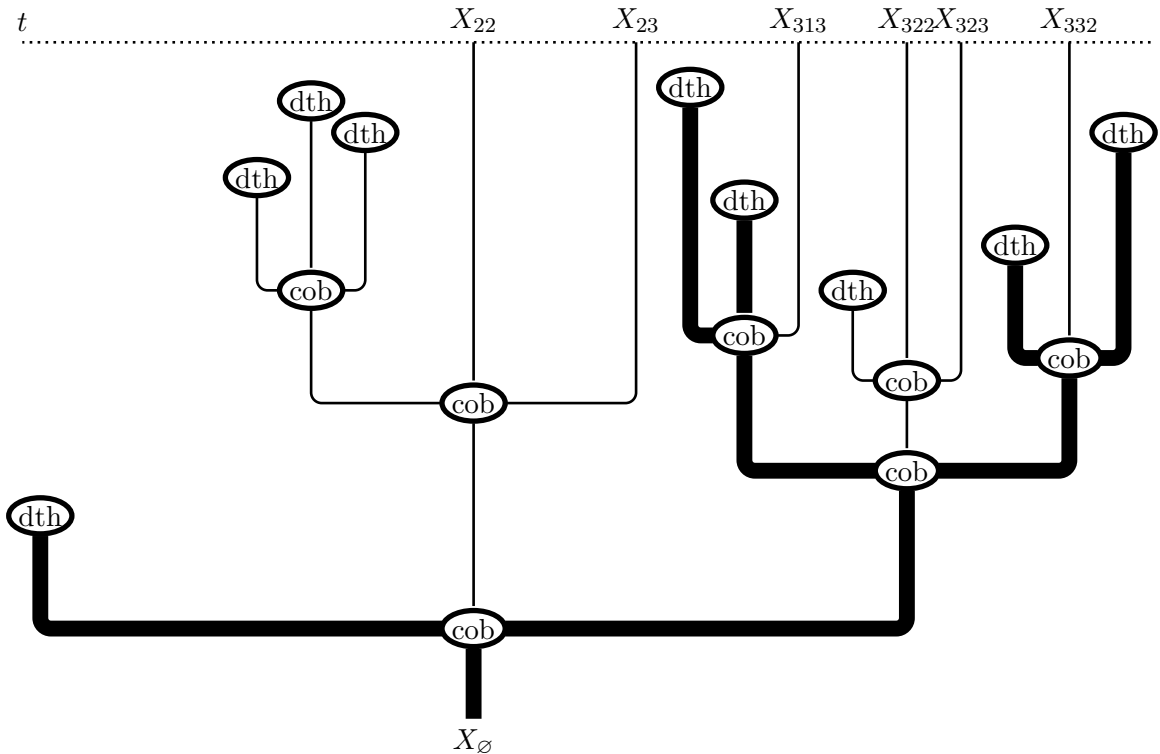

The following lemma illustrates these ideas on our example of a system with
cooperative branching and deaths. Below, $|\U\cap\{\ibf 2,\ibf 3\}|$ denotes
the cardinality of $\U\cap\{\ibf 2,\ibf 3\}$. See Figure~\ref{fig:rootdet} for
an example.

\bl[The uniqueness regime]
Let\label{L:cobU} $S=\{0,1\}$ and $\Gi:=\{\cob,\dth\}$, and let $\pi$ be as in
(\ref{picob}). Then (\ref{Tfin}) is satisfied if and only if $\al\leq\ha$,
while conditions (i)--(iii) of Lemma~\ref{L:detree} are satisfied if and only
if $\al<4$. Moreover, a finite subtree $\U\sub\St$ is a minimal root
determining subtree if and only if
\be\label{minroot}
\ibf 1\in\U\quand\big|\U\cap\{\ibf 2,\ibf 3\}\big|=1
\quad\mbox{for each $\ibf\in\U$ with }\ga[\omb_\ibf]=\cob.
\ee
\el

Lemma~\ref{L:cobU} shows that in our example of a system with
cooperative branching and deaths, the conditions of Lemma~\ref{L:rootdet} are
in fact equivalent to uniqueness of solutions to the RDE. As the next lemma
shows, this is a consequence of monotonicity.

\bl[Uniqueness for monotone systems]
Assume\label{L:monuni} that $S$ is a finite partially ordered set that
contains a minimal and maximal element, and assume that $\ga[\om]$ is monotone
for each $\om\in\Om$. Then the RDE (\ref{RDE}) has a unique solution if and
only if the equivalent conditions of Lemma~\ref{L:rootdet} are satisfied.
\el

In the remainder of this subsection, we focus on the case that $S=\{0,1\}$ and
$\ga[\om]$ is monotone for all $\om\in\Om$, which allows us to make a
connection to a general duality for monotone particle systems described in
\cite{SS18}. Recall that a set $A\sub\{0,1\}^k$ is increasing if $A\ni x\leq
y$ implies $y\in A$. A \emph{minimal element} of $A$ is an $y\in A$ such that
$A\ni x\leq y$ implies $x=y$. If $K$ is a nonempty finite set and
$G:\{0,1\}^K\to\{0,1\}$ is a monotone map, then the inverse image
$G^{-1}(\{1\})$ is an increasing set. We set
\be
Y_G:=\big\{y\in\{0,1\}^K:y\mbox{ is a minimal element of }G^{-1}(\{1\})\big\}.
\ee
Then
\be\label{GY}
G(x)=1\quad\mbox{if and only if}\quad x\geq y\mbox{ for some }y\in Y_G.
\ee
These formulas remain true when $K=\emptyset$, provided we define
$\{0,1\}^\emptyset:=\{\wurz\}$ and we let $Y_G:=\{\wurz\}$ if $G(\wurz)=1$ and 
$Y_G:=\emptyset$ if $G(\wurz)=0$.

Recall from Section~\ref{S:treerep} that $(\nab\St_t,G_t)_{t\geq 0}$ is a Markov
process. If $S=\{0,1\}$ and $\ga[\om]$ is monotone for all $\om\in\Om$, then
the random map $G_t:\{0,1\}^{\nab\St_t}\to\{0,1\}$ is monotone for each $t\geq
0$. In view of this, by (\ref{GY}), $G_t$ is uniquely characterized by
$Y_{G_t}$ and hence $(\nab\St_t,Y_{G_t})_{t\geq 0}$ is a Markov process too. For
a system with cooperative branching and deaths, this process has been defined
before in \cite[Section~I.2.1.2]{Mac17}. As explained in more detail there, it
can be seen as the mean-field limit of a general dual for monotone particle
systems described in \cite[Section~5.2]{SS18}.

Let $S=\{0,1\}$, let $\ga[\om]$ be monotone for all $\om\in\Om$, and let $\U$
be a subtree of $\St$ that contains the root $\wurz$. Borrowing
  terminology from percolation theory, we say
that $\Op$ is a \emph{open subtree of $\U$} if $\wurz\in\Op\sub\U\cup\nab\U$
and
\be\label{goodef}
A_\ibf:=\big\{j\in[\kappa(\omb_\ibf)]:\ibf j\in\Op\big\}
\quad\mbox{satisfies}\quad
1_{A_\ibf}\in Y_{\ga[\omb_\ibf]}\quad\forall\ibf\in\Op\cap\U,
\ee
where we use the convention that $1_{A_\ibf}:=\wurz$ if $\kappa(\omb_\ibf)=0$.

\bl[Open subtrees]
Assume\label{L:good} that $S=\{0,1\}$ and $\ga[\om]$ is monotone for all
$\om\in\Om$. Then
\bc\label{dualow}
\dis\nu_{\rm upp}(\{1\})
&=&\dis\P\big[\mbox{there exists an open subtree of }\St\big],\\[5pt]
\dis\nu_{\rm low}(\{1\})
&=&\dis\P\big[\mbox{there exists a finite open subtree of }\St\big].
\ec
If moreover $\ga[\om](0,\ldots,0)=0$ for each $\om\in\Om$, then
\be\label{YGS}
Y_{G_t}=\big\{y\in\{0,1\}^{\nab\St_t}:\exists\mbox{ open subtree $\Op$ of }\St_t
\mbox{ s.t.\ }y=1_{\Op\cap\nab\St_t}\big\}.
\ee
\el

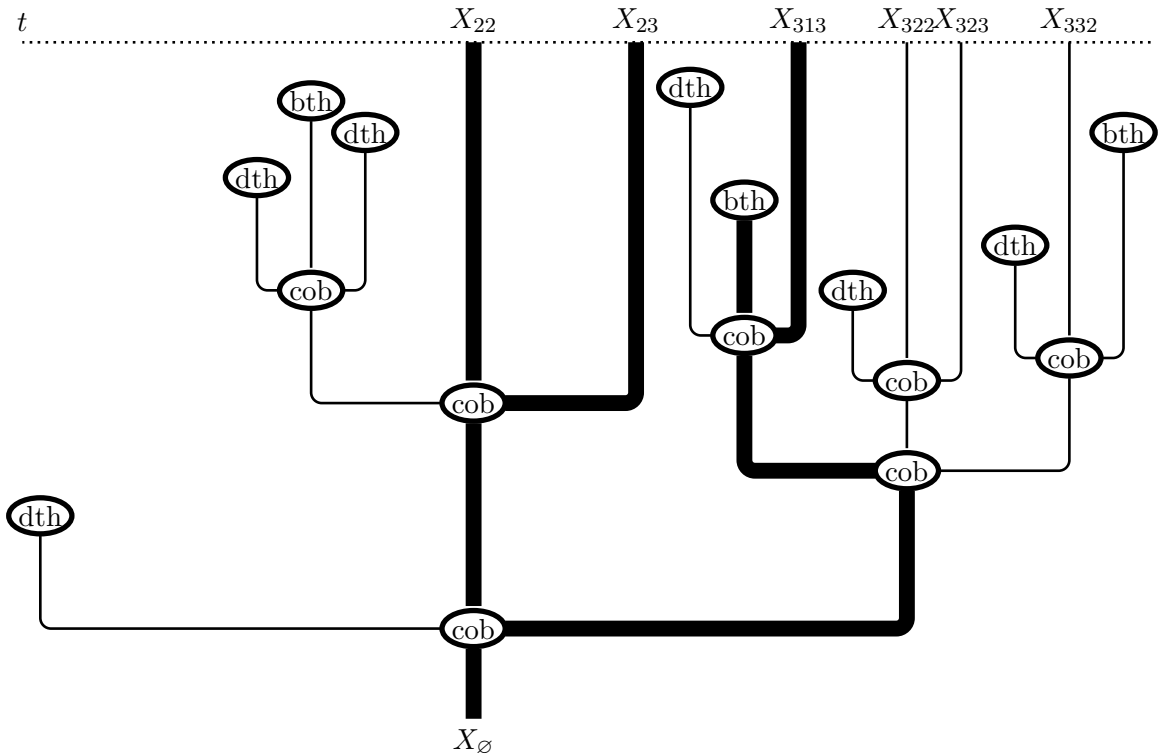
\begin{figure}
\centering

\begin{tikzpicture}[scale=0.6]

\node (wurz) at (0,2) {$\rm{cob}$};
\node (1) at (-9.6,4.5) {$\rm{dth}$};
\node (2) at (0,7) {$\rm{cob}$};
\node (3) at (9.6,5.5) {$\rm{cob}$};
\node (21) at (-3.6,9.5) {$\rm{cob}$};
\node[above] (22) at (0,15) {$X_{22}$};
\node[above] (23) at (3.6,15) {$X_{23}$};
\node (31) at (6,8.5) {$\rm{cob}$};
\node (32) at (9.6,7.5) {$\rm{cob}$};
\node (33) at (13.2,8) {$\rm{cob}$};
\node (211) at (-4.8,12) {$\rm{dth}$};
\node (212) at (-3.6,13.7) {$\rm{bth}$};
\node (213) at (-2.4,13) {$\rm{dth}$};
\node (311) at (4.8,14) {$\rm{dth}$};
\node (312) at (6,11.5) {$\rm{bth}$};
\node[above] (313) at (7.2,15) {$X_{313}$};
\node (321) at (8.4,9.5) {$\rm{dth}$};
\node[above] (322) at (9.6,15) {$X_{322}$};
\node[above] (323) at (10.8,15) {$X_{323}$};
\node (331) at (12,10.5) {$\rm{dth}$};
\node[above] (332) at (13.2,15) {$X_{332}$};
\node (333) at (14.4,13) {$\rm{bth}$};

\draw[line width=6pt] (0,0)--(wurz);
\draw[line width=2pt] (wurz) ellipse (0.7cm and 0.4cm);
\node[below] at (0,0) {$X_\wurz$};

\draw[line width=1pt,rounded corners] (-0.7,2)--(-9.6,2)--(1);
\draw[line width=2pt] (1) ellipse (0.7cm and 0.4cm);

\draw[line width=6pt] (0,2.5)--(2);
\draw[line width=2pt] (2) ellipse (0.7cm and 0.4cm);

\draw[line width=6pt,rounded corners] (0.7,2)--(9.6,2)--(3);
\draw[line width=2pt] (3) ellipse (0.7cm and 0.4cm);

\draw[line width=1pt,rounded corners] (-0.7,7)--(-3.6,7)--(21);
\draw[line width=2pt] (21) ellipse (0.7cm and 0.4cm);

\draw[line width=6pt] (0,7.5)--(22);

\draw[line width=6pt,rounded corners] (0.7,7)-- (3.6,7)--(23);

\draw[line width=6pt,rounded corners] (8.9,5.5)--(6,5.5)--(31);
\draw[line width=2pt] (31) ellipse (0.7cm and 0.4cm);

\draw[line width=1pt] (9.6,6)--(32);
\draw[line width=2pt] (32) ellipse (0.7cm and 0.4cm);

\draw[line width=1pt,rounded corners] (10.3,5.5)--(13.2,5.5)--(33);
\draw[line width=2pt] (33) ellipse (0.7cm and 0.4cm);

\draw[line width=1pt,rounded corners] (-4.3,9.5)--(-4.8,9.5)--(211);
\draw[line width=2pt] (211) ellipse (0.7cm and 0.4cm);

\draw[line width=1pt,rounded corners] (-3.6,10)--(212);
\draw[line width=2pt] (212) ellipse (0.7cm and 0.4cm);

\draw[line width=1pt,rounded corners] (-2.9,9.5)--(-2.4,9.5)--(213);
\draw[line width=2pt] (213) ellipse (0.7cm and 0.4cm);

\draw[line width=1pt,rounded corners] (5.3,8.5)--(4.8,8.5)--(311);
\draw[line width=2pt] (311) ellipse (0.7cm and 0.4cm);

\draw[line width=6pt,rounded corners] (6,9)--(312);
\draw[line width=2pt] (312) ellipse (0.7cm and 0.4cm);

\draw[line width=6pt,rounded corners] (6.7,8.5)--(7.2,8.5)--(313);

\draw[line width=1pt,rounded corners] (8.9,7.5)--(8.4,7.5)--(321);
\draw[line width=2pt] (321) ellipse (0.7cm and 0.4cm);

\draw[line width=1pt] (9.6,8)--(322);

\draw[line width=1pt,rounded corners] (10.3,7.5)--(10.8,7.5)--(323);

\draw[line width=1pt,rounded corners] (12.5,8)--(12,8)--(331);
\draw[line width=2pt] (331) ellipse (0.7cm and 0.4cm);

\draw[line width=1pt,rounded corners] (13.2,8.5)--(332);

\draw[line width=1pt,rounded corners] (13.9,8)--(14.4,8)--(333);
\draw[line width=2pt] (333) ellipse (0.7cm and 0.4cm);

\draw[line width=1pt, dotted](15,15)--(-10,15) node[above]{$t$};

\end{tikzpicture}
\caption{An open subtree of $\S_t$ for a system with cooperative branching,
  deaths, and births, represented by the maps $\cob$ and $\dth$ from
    (\ref{cobdth}) and $\bth$ from (\ref{bth}). In this example, if
  $X_{22},X_{23}$, and $X_{313}$ are 1, then the state at the root $X_\wurz$
  is also 1.}
\label{fig:good}
\end{figure}

We note that formula (\ref{dualow}) can be generalized to more general finite
partially ordered sets $S$, see Lemma~\ref{L:lowdu} below.
Again, it will be useful to illustrate our definitions on the concrete example
of a system with cooperative branching and death. To make the example more
interesting, we add a birth map $\bth:S^0\to S$, which is defined similarly to
the death map as
\be\label{bth}
\bth(\wurz):=1.
\ee
The following lemma describes open subtrees for a system described by the
maps $\cob,\dth,\bth$; see Figure~\ref{fig:good} for an illustration.

\bl[Systems with cooperative branching, deaths, and births]
Let\label{L:birth} $S=\{0,1\}$, $\Gi:=\{\cob,\dth,\bth\}$, with
\be\label{picobbra}
\pi\big(\{\cob\}\big):=\al\geq 0,\quad\pi\big(\{\dth\}\big):=1,
\quand\pi\big(\{\bth\}\big):=\bet\geq 0.
\ee
Let $\U$ be a subtree of $\St$ that contains the root and let
$\Op\sub\U\cup\nab\U$ satisfy $\wurz\in\Op$. Then $\Op$ is an open subtree
of $\U$ if and only if for all $\ibf\in\Op\cap\U$,
\be\label{cobgood}
\ga[\omb_\ibf]\neq\dth
\quand
\{\ibf 1,\ibf 2,\ibf 3\}\cap\Op=\{\ibf 1\}\mbox{ or }\{\ibf 2,\ibf 3\}
\mbox{ if }\ga[\omb_\ibf]=\cob.
\ee
\el

We can think of open subtrees as a generalization of the open paths from
oriented percolation. Outside of a mean-field setting, using ideas from
\cite[Section~5.2]{SS18}, one can characterize the upper invariant law of
quite general monotone particle system in terms of ``open structures'' that in
general are neither paths nor trees.

\section{Discussion}

This section is divided into four subsections. In Subsection~\ref{S:Moran}, we
discuss the relation of our work to \cite{BCH18}, who in parallel to our work
have studied Moran models that generalize our running example of a system with
cooperative branching and deaths. In Subsection~\ref{S:meandiscus}, we compare
our results and methods with the existing literature on mean-field limits. In
Subsection~\ref{S:Problems}, we state open problems and we conclude in
Subsection~\ref{S:Outline} with an outline of the proofs.

\subsection{A Moran model with frequency-dependent selection}\label{S:Moran}

Let $\bra:\{0,1\}^2\to\bra$ be the branching map defined as
\be
\bra(x_1,x_2):=x_1\vee x_2\qquad\big(x_1,x_2\in\{0,1\}\big).
\ee
Consider a system with $S=\{0,1\}$, $\Gi:=\{\cob,\bra,\dth,\bth\}$, with rates
\be\label{piBCH}
\pi\big(\{\cob\}\big):=\ga,\quad\pi\big(\{\bra\}\big):=s,
\quad\pi\big(\{\dth\}\big):=u\nu_0,\quand\pi\big(\{\bth\}\big):=u\nu_1.
\ee
$\ga\geq 0$, $s>0$, $\nu_0,\nu_1\geq 0$ with $\nu_0+\nu_1=1$, and $u>0$.
If $(\mu_t)_{t\geq 0}$ solves the corresponding mean-field equation (\ref{pimean}),
then $p_t:=\mu_t(\{1\})$ solves the ODE (compare (\ref{coopmean}))
\be\label{Moran}
\dif{t}p_t=\ga p_t^2(1-p_t)+sp_t(1-p_t)-u\nu_0p_t+u\nu_1(1-p_t)
\qquad(t\geq 0).
\ee
This equation has an interpretation in terms of a Moran model describing a 
fixed population of $N$ individuals which can be of two types, 0 and 1, where
type 1 is fitter than type 0. The parameter $\ga$ is the \emph{frequency
  dependent selection rate}, $s$ is the \emph{selection rate}, $u$ is the
\emph{mutation rate}, and $\nu_0,\nu_1$ are mutation probabilities. The
frequency dependent selection is of a type that is especially appropriate to
describe an advantageous, (partially) recessive gene in a diploid population.

In parallel to our work, Moran models of this form have been studied by Ellen
Baake, Fernando Cordero, and Sebastian Hummel in \cite{BCH18}. A notational
difference between their work and the discussion here is that they denote the
fitter type by 0, so their \cite[formula (2.1)]{BCH18} is our (\ref{Moran})
rewritten in terms of $y(t)=1-p_t$ and with the roles of $\nu_0$ and $\nu_1$
reversed. They prove that (\ref{Moran}) describes the mean-field limit of a
class of Moran models \cite[Prop.~4.1]{BCH18} and that in the limit
$N\to\infty$, the genealogy of a single individual is described by an
\emph{Ancestral Selection Graph} (\emph{ASG}) $\Ai_t$, which in our notation
corresponds to
\be
\Ai_t:=\big(\nab\St_t,\S_t,(\omb_\ibf,\sig_\ibf)_{\ibf\in\St_t}\big),
\ee
i.e., this is the random tree with maps attached to its branch points depicted
in Figure~\ref{fig:treemap}.

The authors of \cite{BCH18} define a duality function $H(\Ai_t,p)$ which
corresponds to the duality function in (\ref{Hdef}) after the identification
$\mu(\{1\})=p$. (Here we have slightly rephrased things compared to the
different conventions in \cite{BCH18}, where 0 denotes the fitter type and $y$
is the frequency of the unfit type.) In \cite[Lemma~4.4]{BCH18}, they
show that $H(\Ai_t,p)$ can be calculated by concatenating the higher-level
maps $\ch\ga[\omb_\ibf]$ with $\ibf\in\S_t$. For example, the equation
$y=y_1[y_2+y_3-y_2y_3]$ in \cite[Lemma~4.4 (4)]{BCH18} can be rewritten in
terms of $p=1-y$ as $p=\widehat\cob(p_1,p_2,p_3)$ with $\widehat\cob$ as in
(\ref{coblev}).

In \cite[Section~5]{BCH18}, it is shown that the ASG $\Ai_t$ can be simplified
a lot, while retaining all information necessary to calculate the duality
function $H(\Ai_t,p)$. This is done in three steps, I, IIa, and IIb.

In the step I, the ASG is \emph{pruned}. This is a process in which parts
of the tree that are irrelevant for the map $G_t$ are cut off. In particular,
if the function $G_t$ is constant, then the pruned $G_t$ consists of a single
edge ending in one of the maps $\dth$ or $\bth$. In the remaining case,
the pruned ASG is a finite tree where each branch point is marked with one
of the maps $\cob$ and $\bra$.

In steps IIa and IIb, the pruned ASG is \emph{stratified}. In step IIa, the
tree structure is changed in such a way that starting at the root, one first
sees a ternary tree containing only the map $\cob$, and then at the leaves of
this ternary tree, there are attached binary trees containing only the map
$\bra$. In step IIb, each binary tree is replaced by an integer $n\geq 0$
which records the number of leaves of the binary tree.

The result of this is a simplified process, the \emph{stratified ASG} $\Ti_t$,
which contains all necessary information about the ASG $\Ai_t$ in the sense
that there exists a function $\Hi(\Ti_t,p)$ such that
$\Hi(\Ti_t,p)=H(\Ai_t,p)$ \cite[Thm~5.13]{BCH18}. In particular, solutions of
(\ref{Moran}) can be represented as $p_t=\E[\Hi(\Ti_t,p_0)]$
(\cite[Thm~6.2]{BCH18} (compare (\ref{oldmeanrep})).

One can now check (compare Lemma~\ref{L:condroot} below) that
$\rho_t:=\P[\Hi(\Ti_t,p)\in\,\cdot\,]$ solves the higher-level mean-field
equation with initial state $\rho_0=\de_p$, where we use the identification
$\Pc(\{0,1\})\cong[0,1]$. In \cite[Thm~6.5]{BCH18}, it is observed that
$M_t:=\Hi(\Ti_t,p_0)$ is a bounded sub- or supermartingale for each
$p_0\in[0,1]$ and hence converges to an a.s.\ limit $\Hi_\infty(p_0)$. In
\cite[Prop.~6.6]{BCH18}, it is proved that if $p_0$ is not an unstable fixed
point of (\ref{Moran}), then $\Hi_\infty(p_0)$ is a Bernoulli random variable
with parameter $\lim_{t\to\infty}p_t$.

Our Propositions \ref{P:minconv} and \ref{P:hlRTP} imply
 that if $p_0$ is a fixed point of
(\ref{Moran}), then $\Hi_\infty(p_0)$ is a Bernoulli random variable if and
only if the RTP corresponding to $p_0$ is endogenous. Thus,
\cite[Prop.~6.6]{BCH18} implies that for the model in (\ref{piBCH}), RTPs
corresponding to stable fixed points are always endogenous. Since all stable fixed
points of (\ref{Moran}) are in fact lower or upper solutions, this
alternatively also follows from our Proposition~\ref{P:lowup}.

In the special case $s=0$ and $\nu_1=0$, \cite[Prop.~6.6]{BCH18} follows
alternatively from our Theorem~\ref{T:coblev}, which completely describes the
long-time behaviour of solutions to the higher-level mean-field
equation not just for initial states of the form $\rho_0=\de_p$, but for
general initial states.

\subsection{Mean-field limits}\label{S:meandiscus}

If $N$ Markov processes interact in a way that is symmetric under permutations
of the $N$ coordinates, then it is frequently possible to obtain a nontrivial
limit as $N\to\infty$. Such limits are generally called \emph{mean-field
  limits}. In the mean-field limit, the individual processes behave
asymptotically independently, but with transition probabilities that depend on
the average behavior of all processes. For systems of interacting diffusions,
this principle was demonstrated by McKean in his analysis of the Vlasov
equation \cite{McK66}. Consequently, mean-field limits are also called
McKean-Vlasov limits. There exists an extensive literature on the topic. Most
work has focused on interacting diffusions, but jump processes have also been
studied \cite{ST85,ADF18}. An elementary introduction to mean-field limits for
interacting particle systems is given in \cite[Chapter~3]{Swa17}.

In a biological setting, well-mixing populations converge in the mean-field
limit to the solution of a deterministic ODE. Similarly, spatial populations
with strong local mixing can be expected to converge, after an appropriate
rescaling, to the solution of a determinstic PDE. For interacting particle
systems whose dynamics have an exclusion process component with a large rate,
this intuition was made rigorous by De Masi, Ferrari and Lebowitz
\cite[Thm~2]{DFL86}. They state their theorem only for processes whose state
space $S$ consists of two points, and only prove the theorem for one
particular one-dimensional example, but sketch how the proof should be adapted
to the general case. In \cite[Thm~1]{DN94}, a version of the theorem is stated
where $S$ can be any finite set; it is claimed that the proof is again the
same.

In our running example of a particle system with cooperative branching
and deaths, the limiting PDE takes the form
\be\label{coopPDE}
\dif{t}p_t(x)=\ha\diff{x}p_t(x)+\al p_t(x)^2\big(1-p_t(x)\big)-p_t(x)
\qquad(t\geq 0,\ x\in\R).
\ee
This PDE was used in \cite{Nob92} to derive asymptotic properties of the
associated spatial particle system with strong mixing. We can view
(\ref{coopPDE}) as a spatial version of the ODE (\ref{coopmean}); in
particular, if $p_0(x)=p_0$ does not depend on $x$, then $p_t=p_t(x)$ solves
(\ref{coopmean}).

The intuition behind (\ref{coopPDE}), and more general PDEs of this type, is
easily explained. In the strong mixing limit, the genealogy of a single site
should be described by a branching process as in Figure~\ref{fig:treemap}
where in addition, each particle has a position in $\R$, which moves according
to an independent Brownian motion. Convergence to the PDE should then follow
from, on the one hand, convergence of the genealogies to a system of branching
Brownian motions with random maps attached to their branching events, and, on
the other hand, a representation in the spirit of Theorem~\ref{T:meanrep} of
solutions of the PDE (\ref{coopPDE}) in terms of such a system of branching
Brownian motions.

The proof of \cite[Thm~2]{DFL86} is indeed based on this sort of dual
approach, although one would wish that they had given a more explicit statement
of the stochastic representation of solutions of their general PDE.  Our proof
of Theorem~\ref{T:meanlim} follows the same strategy, i.e., we first prove the
stochastic representation of solutions to the mean-field equation
(Theorem~\ref{T:meanrep}) and then use this to prove our convergence result
(Theorem~\ref{T:meanlim}).

\subsection{Open problems}\label{S:Problems}

In the present paper, we have adapted results from \cite{AB05,MSS18} about
discrete-time Recursive Tree Processes and endogeny to the continuous-time
setting, and applied our general results on a concrete system with cooperative
branching and deaths. Among other things, we proved that for $\al>4$, the RTPs
corresponding to $\nu_{\rm low}$ and $\nu_{\rm upp}$ are endogenous but the
RTP corresponding to $\nu_{\rm mid}$ is not. The proof was based on an
analysis of the bivariate mean-field equation. Here, it was convenient to be
able to analyse a differential equation, as an analysis of the associated
discrete-time bivariate evolution would have been possible, but more messy.

Our work leaves a number of questions unanswered, both in the general setting
and more specifically for our running example with $\Gi:=\{\cob,\dth\}$  and
$\pi$ as in (\ref{picob}). Concerning the latter, we pose the following
questions.
\begin{itemize}
\item[] \textbf{Open Problem~1} Not every measure $\mu^{(n)}\in\Pc_{\rm
  sym}\big(\{0,1\}^n\big)$ is the $n$-th moment measure of a measure
  $\rho\in\Pc\big(\Pc(\{0,1\})\big)$. Determine all symmmetric solutions of
  the $n$-variate RDE, for general $n\geq 3$, and their domains of attraction.
\item[] \textbf{Open Problem~2} Same as Open Problem~1 but without the
  symmetry assumption and for general $n\geq 2$.
\item[] \textbf{Open Problem~3} Prove that apart from the atom at zero, the
  law $\un\nu_{\rm mid}$, viewed as a probability law on $[0,1]$, has a smooth
  density with respect to the Lebesgue measure.
\item[] \textbf{Open Problem~4} Determine the aymptotics of the distribution
  function $F$ of $\un\nu_{\rm mid}$ near $0$ and $1$.
\item[] \textbf{Open Problem~5} For the more general model in (\ref{piBCH}),
  is it true that unstable fixed points of the mean-field equation that
    separate the domains of attraction of two stable fixed points correspond
  to nonendogenous RTPs? Is the picture for the higher-level RDE
  the same?
\end{itemize}
Partly inspired by our concrete example, we ask the following problems in the
general setting.
\begin{itemize}
\item[] \textbf{Open Problem~6} Can (\ref{sumr}) be relaxed to allow for
  branching processes $(\nab\St_t)_{t\geq 0}$ that are nonexplosive but have
  infinite mean?
\item[] \textbf{Question~7} Are there general results linking the
  (in)stability of fixed points of the mean-field equation to (non)endogeny of
  the related RTP?
\item[] \textbf{Question~8} In our example, the higher-level RDE has two
  solutions $\un\nu_{\rm mid}$ and $\ov\nu_{\rm mid}$ with mean $\nu_{\rm
    mid}$, of which the former is stable and the latter is unstable. Is this a
  general phenomenon in the nonendogenous case? Can one prove nonendogeny of
  an RTP corresponding to a solution $\nu$ of the RDE by showing that $\ov\nu$
  is unstable?
\item[] \textbf{Question~9} Are there examples of higher-level RDEs that have
  solutions $\rho\not\in\{\un\nu,\ov\nu\}$?
\item[] \textbf{Open Problem~10} Is the higher-level RTP
  $(\omb_\ibf,\xi_\ibf)_{\ibf\in\St}$ from Proposition~\ref{P:hlRTP} always
  endogenous?
\end{itemize}
Finally, we mention the problem of proving nonendogeny for the frozen
percolation of \cite{Ald00}, which to our knowledge is still open. Although we
did not attempt to solve this problem here, one might hope that the methods of
the present paper can provide a useful new point of view on this old problem.

\subsection{Outline of the proofs}\label{S:Outline}

In the remainder of the paper, we prove all results stated so far, except for
Theorems~\ref{T:bivar} and \ref{T:Stras} as well as Proposition~\ref{P:hlRTP},
which we cite from \cite[Thm~1, Thm~13, and Prop~4]{MSS18}.

In Section~\ref{S:meanfield} we prove Theorem~\ref{T:mean}, Propositions
  \ref{P:TVcont} and \ref{P:gacon}, and Lemma~\ref{L:simpmean}, which state
  elementary properties of solutions of the mean-field equations
  (\ref{vecmean}) and (\ref{mean}), as well as Theorem~\ref{T:meanrep}, which
gives a stochastic representation of solutions of the mean-field equation in
terms of finite recursive tree processes. In Section~\ref{S:meanfield3}, we
use this stochastic representation to prove Theorem~\ref{T:meanlim} about
convergence of finite systems to a solution of the mean-field equation.

In Section~\ref{S:RTP}, we prove our main results about RTPs with continuous
time, which are largely analogous to known results from the discrete-time
setting. Basic results are Lemma~\ref{L:consis} and
Proposition~\ref{P:contRTP}, as well as Lemma~\ref{L:RTP} which deals with
discrete time and is a slight reformulation of known results. Following
\cite{AB05}, Theorem~\ref{T:bivar2} links the $n$-variate equation to
endogeny, while Propositions~\ref{P:levmom} and \ref{P:minconv} are
concerned with the higher-level equation, and closely follow ideas from
\cite{MSS18}.

In Section~\ref{S:fur} we prove some additional results about RTPs, first
Proposition~\ref{P:lowup}, which generalizes \cite[Lemma~15]{AB05} and shows that
upper and lower solutions of a monotonous RDE are always endogenous, and then
Lemmas~\ref{L:rootdet}, \ref{L:detree}, \ref{L:monuni}, and \ref{L:good} which
give conditions for uniqueness in a general setting and then more specifically
for monotone systems.

In Section~\ref{S:coop}, finally, we have collected all proofs that deal
specifically with our running example of a system with cooperative branching
and deaths. The first such result is Proposition~\ref{P:copbiv} about the
bivariate equation, which is a two dimensional ODE for which by elementary
means we find all fixed points and their domains of attraction. By combining
Proposition~\ref{P:copbiv} with ideas involving the convex order we then prove
the much stronger Theorem~\ref{T:coblev} which gives all fixed points and
domains of attraction for the higher-level equation. The picture is then
completed by the proofs of Lemma~\ref{L:nontriv}, which gives some properties
of the nontrivial fixed point of the higher-level equation, as well as
Lemmas~\ref{L:cobU} and \ref{L:birth} which illustrate ideas from
Section~\ref{S:fur} in the concrete set-up of our example.


\section{The mean-field equation}\label{S:meanfield}

In this section, we prove Theorems~\ref{T:mean} and \ref{T:meanrep}, which
state that the mean-field equation (\ref{mean}) has a unique solution and can
be represented in terms of a random tree generated by a branching process,
with random maps attached to its vertices. In addition, we also prove
Propositions \ref{P:TVcont} and \ref{P:gacon}, as well as
Lemma~\ref{L:simpmean}.

In Subsection~\ref{S:formul}, we start with some preliminaries, showing, in
particular, that the integral in (\ref{weakmean}) is well-defined, and
Lemma~\ref{L:simpmean}, which says that mean-field equations of the form
(\ref{vecmean}) can be rewritten in the simpler form (\ref{mean}).

Next, in Subsection~\ref{S:meanfield1}, we prove uniqueness of solutions of
(\ref{mean}), which yields the uniqueness parts of
Theorem~\ref{T:mean}. To prove existence, in
Subsection~\ref{S:meanfield2}, we show that the right-hand side of
(\ref{meanrep}) solves (\ref{mean}), which not only completes the proof
of Theorem~\ref{T:mean} but also yields the stochastic representation that
is Theorem~\ref{T:meanrep}.

The proofs of Propositions \ref{P:TVcont} and \ref{P:gacon}, finally, can be
found in Subsection~\ref{S:contP}.

\subsection{Preliminaries}\label{S:formul}

Recall that we interpret the mean-field equation (\ref{mean}) as in
(\ref{weakmean}), where, by (\ref{Tgdef}),
\be\label{Tgom}
\li T_{\ga[\om]}(\mu),\phi\re=\int_S\!\mu(\di x_1)\cdots\int_S\!\mu(\di x_k)\,
\phi\big(\ga[\om](x_1,\ldots,x_k)\big)\qquad(\om\in\Om_k,\ k\geq 1).
\ee
Since by assumption, $\ga[\om](x_1,\ldots,x_k)$ is jointly measurable in $\om$
and $x_1,\ldots,x_k$, the right-hand side of (\ref{Tgom}) is measurable as a
function of $\om$ and hence the integral in (\ref{weakmean}) is
well-defined.\med

\bpro[of Lemma~\ref{L:simpmean}]
Recall from Subsection~\ref{S:flow} that the basic ingredients that go into
the equation (\ref{vecmean}) are the measure space $(\Omm,\qb)$ and function
$\la$, as well as, for each $\om\in\Omm$ and $1\leq i\leq\la(\om)$, the
function $\ga_i[\om]$ and set $K_i(\om)$. Also, $\kappa_i(\om):=|K_i(\om)|$.
In terms of these basic ingredients we need to define $\Om,\rb,\kappa$, and
$\ga$ as in Subsection~\ref{S:firstintro} so that (\ref{vecmean}) takes the
simpler form (\ref{mean}).

Since we want to replace the integral and sum in (\ref{vecmean}) by a single
integral, we put
\be\label{Omqu}
\Om:=\bigcup_{l=1}^\infty\Omm_l\times[l],
\ee
where as before $\Omm_l:=\{\om\in\Omm:\la(\om)=l\}$ and $[l]:=\{1,\ldots,l\}$,
and we equip $\Om$ with the measure
\be\label{rq}
\rb\big(A\times\{k\}\big):=\qb(A)
\qquad(A\sub\Omm_l\mbox{ measurable, }1\leq k\leq l).
\ee
In general, $\Om$ need not be a Polish space, as required in
Subsection~\ref{S:firstintro}. We will fix this problem at the end of our
proof, but for the sake of the presentation we neglect it for the moment being.
We define $\kappa:\Om\to\N$ as in Subsection~\ref{S:firstintro} by
$\kappa(\om,i):=\kappa_i(\om)$, where the right-hand side is the function from
Subsection~\ref{S:flow}. We write
\be\label{Kord}
K_i(\om)=\big\{j_1,\ldots,j_{\kappa_i(\om)}\big\}
\quad\mbox{with}\quad
j_1<\cdots<j_{\kappa_i(\om)}.
\ee
Since $\ga_i[\om](x_1,\ldots,x_{\la(\om)})$ depends only on coordinates in
$K_i(\om)$, there exists a function $\ga[\om,i]:S^{\kappa(\om,i)}\to S$ such
that
\be
\ga_i[\om](x_1,\ldots,x_{\la(\om)})=\ga[\om,i](x_{j_1},\ldots,x_{j_{\kappa(\om,i)}})
\qquad(x\in S^{\la(\om)}).
\ee
Note that $T_{\ga_i[\om]}=T_{\ga[\om,i]}$ by (\ref{Tgdef}).
As in (\ref{kappa}), we can associate $\ga[\om,i]$ with a function that is
defined on $S^{\N_+}$ but depends only on the first $\kappa(\om,i)$
coordinates. We take this as our definition of the function $\ga:\Om\times
S^{\N_+}\to S$ from Subsection~\ref{S:firstintro}. It follows from
(\ref{vecmeas}) that $\ga[\om,i](x)$ is jointly measurabe as a function of
$(\om,i)$ and $x$.

Replacing the integral and sum in (\ref{vecmean}) by a single integral over
$\rb$ as defined in (\ref{rq}), using the fact that $T_{\ga_i[\om]}=T_{\ga[\om,i]}$
we see that (\ref{vecmean}) can be rewritten as
\be\label{intmean}
\dif{t}\mu_t=\int_\Om\!\rb(\di\om)\big\{T_{\ga[\om,i]}(\mu_t)-\mu_t\big\},
\ee
which coincides with (\ref{mean}). The condition that $\rb$ should be a finite
measure translates to (\ref{rmom2})~(i), while the condition (\ref{sumr}),
written in terms of $\qb$, becomes (\ref{rmom2})~(ii). Moreover, if $\qb$
satisfies (\ref{qasco}), then $\rb$ satisfies (\ref{asco}).

We still have to fix the problem that $\Om$, as defined in (\ref{Omqu}), is in
general not a Polish space. There are several possible ways to fix
this.\footnote{For example, we can strengthen our assumptions on $\la$ in the
  sense that $\{\om:\la(\om)=l\}$ is a $G_\de$-set for each $l\in\N_+$, or we
  can relax our assumptions on $\Om$ allowing it to be a Lusin space, instead
  of just a Polish space, throughout.} The solution we will choose is to
replace $\Om$ by the Polish space
\be\label{ovOmqu}
\ov\Om:=\bigcup_{l=1}^\infty\ov{\Omm_l}\times[l],
\ee
were $\ov{\Omm_l}$ denotes the closure of $\Omm_l$ in $\Omm$. We view $\rb$ as
a measure on $\ov\Om$ that is concentrated on $\Om$ and extend $\kappa$ and
$\ga$ in a measurable way to the larger space, which is possible since $\Om$
is a measurable subset of $\ov\Om$. Since $\rb$ is concentrated on $\Om$, it
does not matter how we extend $\kappa$ and $\ga$ as this has no effect on
(\ref{intmean}).
\epro

\subsection{Uniqueness}\label{S:meanfield1}

In the present section, we prove that under the assumption (\ref{sumr}),
solutions to (\ref{mean}) are unique, which settles the uniqueness part
of Theorem~\ref{T:mean}.

Below, we let $\Mi(S)$ denote the space of all finite signed measures on $S$.
The total variation norm has already been mentioned several times. There are
two conventional definitions, which differ by a factor 2. We will use the
definition
\be\label{TVdef}
\|\mu\|:=\ha\sup_{|f|\leq 1}\Big|\int\! f\,\di\mu\Big|
\qquad\big(\mu\in\Mi(S)\big),
\ee
where the supremum runs over all measurable functions $f:S\to[-1,1]$.
If $X,Y$ are $S$-valued random variables, then it is easy to see that
$\|\mu-\nu\|\leq\P[X\neq Y]$. Conversely, it is well-known
\cite[page~19]{Lin92} that if $\mu,\nu\in\Pc(S)$, then it is possible to
couple $S$-valued random variables $X,Y$ in such a way that
\be\label{TVcoup}
\|\mu-\nu\|=\P[X\neq Y].
\ee 

\bl[Lipschitz continuity]
Let\label{L:TLip} $g:S^k\to S$ be measurable and let $T_g$ be defined as
in (\ref{Tgdef}). Then
\be\label{TgLip}
\big\|T_g(\mu)-T_g(\nu)\big\|\leq k\|\mu-\nu\|
\qquad\big(\mu,\nu\in\Pc(S)\big).
\ee
Moreover, if $T$ is defined as in (\ref{Tdef}), then
\be\label{TLip}
\big\|T(\mu)-T(\nu)\big\|
\leq\Big(|\rb|^{-1}\!\int_\Om\rb(\di\om)\,\kappa(\om)\Big)\,\|\mu-\nu\|
\qquad\big(\mu,\nu\in\Pc(S)\big).
\ee
\el
\bpro
By (\ref{TVcoup}) we can find an $S^2$-valued random variable $(X,Y)$ such that
$\|\mu-\nu\|=\P[X\neq Y]$. Let $(X_1,Y_1),\ldots,(X_k,Y_k)$ be i.i.d.\ copies
of $(X,Y)$. Then, by (\ref{Tgdef}),
\bc
\dis\big\|T_g(\mu)-T_g(\nu)\big\|
&=&\dis\big\|\P\big[g(X_1,\ldots,X_k)\in\,\cdot\,\big]
-\P\big[g(Y_1,\ldots,Y_k)\in\,\cdot\,\big]\big\|\\[5pt]
&\leq&\dis\P\big[g(X_1,\ldots,X_k)\neq g(Y_1,\ldots,Y_k)\big]\\[5pt]
&\leq&\dis\sum_{i=1}^k\P[X_i\neq Y_i]=k\|\mu-\nu\|.
\ec
This proves (\ref{TgLip}). Formula (\ref{TLip}) follows by integrating over
$\om$.
\epro

Our next lemma gives equivalent formulations of the mean-field equation
(\ref{mean}), that will also be useful in the next subsection where we
prove existence of solutions. Below, we interpret an integral of a
measure-valued integrand in the usual way, i.e., $\int_0^t\!\di s\,\mu_s$
denotes the measure defined by
\be
\big\li\int_0^t\!\di s\,\mu_s,\phi\big\re:=
\int_0^t\!\di s\,\li\mu_s,\phi\re
\ee
for any bounded measurable $\phi:S\to\R$.

\bl[Equivalent formulations of the mean-field equation]
Assume\label{L:meanequ} (\ref{sumr}). Let $\half\ni t\mapsto\mu_t\in\Pc(S)$ be
measurable. Then of the following conditions, (i) implies (ii) and (iii).
If $\half\ni t\mapsto\mu_t\in\Pc(S)$ is continuous with
respect to the total variation norm, then all three conditions are equivalent.
\begin{enumerate}
\item For each bounded measurable $\phi:S\to\R$, the function
  $t\mapsto\li\mu_t,\phi\re$ is continuously differentiable and
  $\dis\dif{t}\mu_t=|\rb|\big\{T(\mu_t)-\mu_t\}$ $(t\geq 0)$.
\item $\dis\mu_t=\mu_0+|\rb|\int_0^t\!\di s\,\big\{T(\mu_s)-\mu_s\}$
 $(t\geq 0)$.
\item $\dis\mu_t=e^{-|\rb|t}\mu_0+|\rb|\int_0^t\!\di s\,e^{-|\rb|s}\,T(\mu_{t-s})$
 $(t\geq 0)$.
\end{enumerate}
\el
\bpro
Integrating the equation in (i) from time 0 until time $t$, we see that (i)
implies (ii). Also, we can equivalently write the equation in (i) as
\be
\dif{t}\big(e^{|\rb|t}\mu_t\big)=|\rb|e^{|\rb|t}\,T(\mu_t)
\qquad(t\geq 0).
\ee
Integrating from time 0 until time $t$ now yields
\be\label{print}
e^{|\rb|t}\mu_t
=\mu_0+|\rb|\int_0^t\!\di s\,e^{|\rb|s}\,T(\mu_s).
\ee
Multiplying by $e^{-|\rb|t}$ and substituting $s\mapsto t-s$ in the integral
then yields the equation in (iii).

If $t\mapsto\mu_t\in\Pc(S)$ is continuous with respect to the total variation
norm, then Lemma~\ref{L:TLip} together with (\ref{sumr}) imply that also
$t\mapsto T(\mu_t)\in\Pc(S)$ is continuous with respect to the total variation
norm. It follows that $t\mapsto\li\mu_t,\phi\re$ and $t\mapsto\li
T(\mu_t),\phi\re$ are continuous for each bounded measurable $\phi:S\to\R$. As
a result, the right-hand side of (ii), integrated against any bounded
measurable $\phi$, is continuously differentiable as a function of $t$, and
(ii) implies (i). By the same argument, rewriting (iii) as (\ref{print}) and
differentiating, we see that (iii) implies (i).
\epro

We now prove the promised uniqueness of solutions to
(\ref{mean}). Proposition~\ref{P:TVcont}, which will be proved in
  Subsection~\ref{S:contP} below, shows that the constant $L$ from
  (\ref{Ldef}) is not optimal and can be replaced by the constant $K$ from
  (\ref{Kdef}).

\bl[Uniqueness]
Let\label{L:meanuni} $(\mu_t)_{t\geq 0}$ and $(\nu_t)_{t\geq 0}$ be
solutions of the mean-field equation (\ref{mean}). Then
\be
\|\mu_t-\nu_t\|\leq\ex{Lt}\|\mu_0-\nu_0\|\qquad(t\geq 0),
\ee
where
\be\label{Ldef}
L:=\int_\Om\rb(\di\om)\,\big(\kappa(\om)+1\big).
\ee
\el
\bpro
Equation~(ii) of Lemma~\ref{L:meanequ} implies that
\bc
\dis\|\mu_t-\nu_t\|
&\leq&\dis\|\mu_0-\nu_0\|
+|\rb|\,\Big\|\int_0^t\!\di s\,\big\{T(\mu_s)-\mu_s\}
-\int_0^t\!\di s\,\big\{T(\nu_s)-\nu_s\}\Big\|\\[7pt]
&\leq&\dis\|\mu_0-\nu_0\|
+|\rb|\int_0^t\!\di s\,\big\|T(\mu_s)-T(\nu_s)\big\|
+|\rb|\int_0^t\!\di s\,\|\mu_s-\nu_s\|\\[7pt]
&\leq&\dis\|\mu_0-\nu_0\|+L\int_0^t\!\di s\,\|\mu_s-\nu_s\|,
\ec
where $L=|\rb|+\int_\Om\rb(\di\om)\,\kappa(\om)$ using (\ref{TLip}) of
Lemma~\ref{L:TLip}. The claim now follows from Gronwall's lemma
\cite[Thm~A.5.1]{EK86}.
\epro

\subsection{The stochastic representation}\label{S:meanfield2}

In this section, we prove the following proposition, that settles the
existence part of Theorem~\ref{T:mean}. Together with
Lemma~\ref{L:meanuni}, this completes the proof of Theorem~\ref{T:mean}
and at the same time also proves Theorem~\ref{T:meanrep}.

We work in our usual set-up where $S$ and $\Om$ are Polish spaces,
$\kappa:\Om\to\N$ is measurable, $\ga$ is as in Subsection~\ref{S:firstintro},
and $\rb$ is a nonzero finite measure on $\Om$ satisfying (\ref{sumr}). We fix
$\ovT$ as in Section~\ref{S:treerep} and let $(\omb_\ibf)_{\ibf\in\ovT}$ be
i.i.d.\ with common law $|\rb|^{-1}\rb$. We let $(\sig_\ibf)_{\ibf\in\ovT}$ be
an independent i.i.d.\ collection of exponentially distributed random
variables with mean $|\rb|^{-1}$ and define $\St$, $\St_t$, $\nab\St_t$, and
$G_t$ as in (\ref{Sdef}), (\ref{Stdef}), and (\ref{Gtdef}).

\bp[Recursive tree representation]
For\label{P:meanrep} any $\mu_0\in\Pc(S)$, setting 
\be\label{meanrepdef}
\mu_t:=\E\big[T_{G_t}(\mu_0)\big]\qquad(t\geq 0)
\ee
defines a solution $(\mu_t)_{t\geq 0}$ to the mean-field equation
(\ref{mean}). Moreover, $\half\ni t\mapsto\mu_t$ is continuous with
respect to the total variation norm.
\ep

To prepare for the proof of Proposition~\ref{P:meanrep}, we need one lemma.
Recall that $|\ibf|$ denotes the length of a word $\ibf$, i.e., $|i_1\cdots
i_n|:=n$. Let
\be
\St_{t,(n)}:=\St_t\cap\St_{(n)}
\quad\mbox{with}\quad
\St_{(n)}:=\{\ibf\in\St:|\ibf|<n\}\qquad(n\geq 1).
\ee
Fix $\mu_0\in\Pc(S)$ and using notation as in (\ref{GUdef}), set
\be\label{munt}
\mu_{t,(n)}:=\E\big[T_{G_{t,(n)}}(\mu_0)\big]
\quad\mbox{with}\quad
G_{t,(n)}:=G_{\St_{t,(n)}}\qquad(t\geq 0, n\geq 1),
\ee
and set $\mu_{t,(0)}:=\mu_0$ $(t\geq 0)$. The following lemma is a ``cut-off''
version of Proposition~\ref{P:meanrep}.

\bl[Representation with cut-off]
The\label{L:mun} measures $(\mu_{t,(n)})_{t\geq 0}$ defined in (\ref{munt})
satisfy
\be\label{mun}
\mu_{t,(n)}=e^{-|\rb|t}\mu_0+|\rb|\int_0^t\!\di s\,e^{-|\rb|s}\,T(\mu_{t-s,(n-1)})
\qquad(n\geq 1,\ t\geq 0).
\ee
\el
\bpro
Let $(Y_\ibf)_{\ibf\in\ovT}$ be i.i.d.\ with common law $\mu_0$, independent
of $(\omb_\ibf,\sig_\ibf)_{\ibf\in\ovT}$. Set
\be
X_\ibf:=Y_\ibf\qquad(\ibf\in\nab\St_{t,(n)}),
\ee
and define $(X_\ibf)_{\ibf\in\St_{t,(n)}}$ inductively by
\be
X_\ibf:=\ga[\omb_\ibf](X_{\ibf 1},\ldots,X_{\ibf k_\ibf})\qquad(\ibf\in\St_{t,(n)}).
\ee
Then $X_\wurz=G_{t,(n)}\big((X_\ibf)_{\ibf\in\nab\St_{t,(n)}}\big)$
and hence, in the same way as (\ref{meanrep}) is equivalent to
(\ref{oldmeanrep}),
\be
\mu_{t,(n)}=\P[X_\wurz\in\,\cdot\,].
\ee
Conditioning on $\sig_\wurz$ and then also on $\omb_\wurz$, we see that
\bc
\dis\P[X_\wurz\in\,\cdot\,]
&=&\dis\int_0^\infty\!\P\big[\sig_\wurz\in\di s\big]
\,\P\big[X_\wurz\in\,\cdot\,\big|\,\sig_\wurz=s\big]\\[5pt]
&=&\dis\int_t^\infty\!|\rb|e^{-|\rb|s}\,\di s\,\mu_0
+\int_0^t\!|\rb|e^{-|\rb|s}\,\di s
\,P\big[X_\wurz\in\,\cdot\,\big|\,\sig_\wurz=s\big]\\[5pt]
&=&\dis e^{-|\rb|t}\mu_0
+\int_0^t\!|\rb|e^{-|\rb|s}\,\di s
\;|\rb|^{-1}\!\int_\Om\rb(\di\om)
\,P\big[X_\wurz\in\,\cdot\,\big|\,\sig_\wurz=s,\ \omb_\wurz=\om\big],
\ec
where we have used that $\omb_\wurz$ is independent of $\sig_\wurz$ with law
$|\rb|^{-1}\rb$. We see from this that
\bc
\dis\P[X_\wurz\in\,\cdot\,]
&=&\dis e^{-|\rb|t}\mu_0
+\int_0^t\!e^{-|\rb|s}\,\di s\int_\Om\rb(\di\om)
\,P\big[\ga[\om](X_1,\ldots,X_{\kappa(\om)})\in\,\cdot
\,\big|\,\sig_\wurz=s,\ \omb_\wurz=\om\big]\\[5pt]
&=&\dis e^{-|\rb|t}\mu_0
+\int_0^t\!e^{-|\rb|s}\,\di s\;|\rb|\,T(\mu_{t-s,(n-1)}),
\ec
where we have used (\ref{Tdef}) and the observation that conditional on
$\sig_\wurz=s$ and $\omb_\wurz=\om$, the random variables
$X_1,\ldots,X_{\kappa(\om)}$ are i.i.d.\ with common law $\mu_{t-s,(n-1)}$.
\epro

\bpro[of Proposition~\ref{P:meanrep}]
The condition (\ref{sumr}) guarantees that $(\nab\St_t)_{t\geq 0}$ is a finite
mean branching process; more precisely, by standard theory,
\be\label{EpaTt}
\E\big[|\nab\St_t|\big]=\ex{Kt}
\quad\mbox{with}\quad K:=\int_\Om\!\rb(\di\om)\,(\kappa(\om)-1).
\ee
Fix $\mu_0\in\Pc(S)$ and define $\mu_t$ and $\mu_{t,(n)}$ as in
(\ref{meanrepdef}) and (\ref{mun}). Then the total variation norm distance
between these measures can be bounded by
\be\label{muntot}
\big\|\mu_{t,(n)}-\mu_t\big\|\leq\P\big[\St_{t,(n)}\neq\St_t\big],
\ee
which tends to zero as $n\to\infty$ since $\St_t$ is a.s.\ finite by
(\ref{EpaTt}). In fact, since
$\P[\St_{s,(n)}\neq\St_s]\leq\P[\St_{t,(n)}\neq\St_t]$ for all
$s\leq t$, we have that
\be\label{munconv}
\sup_{0\leq s\leq t}\big\|\mu_{s,(n)}-\mu_{s}\big\|\asto{n}0
\qquad\forall t<\infty.
\ee
Using this and the Lipschitz continuity of $T$ with respect to the total
variation norm (Lemma \ref{L:TLip}), we can let $n\to\infty$ in (\ref{mun}) to
obtain
\be\label{mun2}
\mu_t=e^{-|\rb|t}\mu_0+|\rb|\int_0^t\!\di s\,e^{-|\rb|t}\,T(\mu_{t-s})
\qquad(t\geq 0).
\ee
Since
\be
\big\|\mu_{t+\eps}-\mu_t\big\|\leq\P\big[\St_{t+\eps}\neq\St_t\big],
\ee
using the fact that the branching process $(\nab\St_t)_{t\geq 0}$ a.s.\ does not
jump at deterministic times, we see that $\half\ni t\mapsto\mu_t$ is
continuous with respect to the total variation norm. Using this and
(\ref{mun2}), we see from Lemma~\ref{L:meanequ} that $(\mu_t)_{t\geq 0}$ solves
the mean-field equation (\ref{mean}).
\epro

\subsection{Continuity in the initial state}\label{S:contP}

In this subsection, we prove Propositions \ref{P:TVcont} and \ref{P:gacon}.\med

\bpro[of Proposition~\ref{P:TVcont}]
It follows from Theorem~\ref{T:meanrep} and Lemma~\ref{L:TLip} that
\be\ba{l}
\big\|T_t(\mu)-T_t(\nu)\big\|=\big\|\E[T_{G_t}(\mu)-T_{G_t}(\nu)]\big\|\\[5pt]
\dis\quad\leq
\E\big[\big\|\E[T_{G_t}(\mu)-T_{G_t}(\nu)\,|\,\Fi_t]\big\|\big]
\leq\E\big[|\nab\S_t|\,\|\mu-\nu\|\big]
=e^{Kt}\|\mu-\nu\|,
\ec
where $\Fi_t$ is the filtration defined in (\ref{Fidef}).
\epro

Proposition~\ref{P:gacon} follows from the following two lemmas.

\bl[Continuity of $T$]
Under\label{L:gacon} the condition (\ref{asco}), the operator $T$ in
(\ref{Tdef}) is continuous w.r.t.\ the topology of weak convergence.
\el
\bpro
If $\mu_n\in\Pc(S)$ converge weakly to a limit $\mu_\infty$, then by
Skorohod's representation theorem there exists random variables $X^n$ with
laws $\mu_n$ that converge a.s.\ to a limit $X^\infty$ with law $\mu_\infty$.
Let $\big((X^n_i)_{n\in\N\cup\{\infty\}}\Big)_{i\geq 1}$ be i.i.d.\ copies of such
a sequence $(X^n)_{n\in\N\cup\{\infty\}}$ and let $\omb$ be an independent
random variable with law $|\rb|^{-1}\rb$. Then by (\ref{asco}),
\be
\ga[\omb](X^n_1,\ldots,X^n_{\kappa(\omb)})
\asto{n}
\ga[\omb](X^\infty_1,\ldots,X^\infty_{\kappa(\omb)})\quad{\rm a.s.}
\ee
and hence $T(\mu_n)$ converges weakly to $T(\mu_\infty)$ by (\ref{Tdef}).
\epro

\bl[Continuity in the initial state]
Assume\label{L:mucon} that the operator $T$ in (\ref{Tdef}) is continuous
w.r.t.\ the topology of weak convergence. Then the same is true for the
operators $T_t$ $(t\geq 0)$ defined in
(\ref{Tt}).
\el

\bpro
We need to show that solutions of the mean-field equation (\ref{mean})
are continuous in their initial state, in the sense that if $(\mu^k_t)_{t\geq
  0}$ $(k\in\N\cup\{\infty\})$ are started in initial states such that
$\mu^k_0\Rightarrow\mu^\infty_0$, then $\mu^k_t\Rightarrow\mu^\infty_t$ for
all $t\geq 0$.

To see this, inductively define $\mu^k_{t,(n)}$ as in (\ref{mun}) with
$\mu_0$ replaced by $\mu^k_0$. Using the continuity
of $T$, by induction, we see that $\mu^k_{t,(n)}\Rightarrow\mu^\infty_{t,(n)}$
as $k\to\infty$ for all $n\geq 1$ and $t\geq 0$. By (\ref{muntot}), for each
bounded continuous $\phi:S\to\R$, the quantity $\li\mu^k_{t,(n)},\phi\re$
converges to $\li\mu^k_t,\phi\re$ uniformly in $k\in\N\cup\{\infty\}$, which
allows us to conclude that $\li\mu^k_t,\phi\re\to\li\mu^\infty_t,\phi\re$ as
$k\to\infty$ for all $t\geq 0$.
\epro

\section{Approximation by finite systems}\label{S:meanfield3}

\subsection{Main line of the proof}\label{S:proofline}

In this section, we prove Theorem~\ref{T:meanlim}.  The basic idea, which
already goes back to \cite{DFL86}, is that in the mean-field limit, the
genealogy of a site converges to a branching process, and sites are
independent in the limit. More precisely, consider $n$ sites, sampled
uniformly at random from $[N]$. To find out what their states are at time $t$,
we follow the sites back until the last time when a random map is applied that
has the potential to change the state of one of our sites. At this point, we
stop following that given site but replace it by the sites that are relevant
for the outcome of the map at the given site, and we continue in this
way. When $N$ is large, the new sites that are added in each step are with
high probability sites we have not been following before, so that  in the limit 
we obtain a branching process with random maps attached to its
branch points. Making this idea precise yields the following proposition,
that will be proved in Subsection~\ref{S:samp} below.

\bp[State at sampled sites]
For\label{P:samp} each $N\in\N_+$ let $(X^{(N)}(t))_{t\geq 0}$ be a process
as in Theorem~\ref{T:meanlim} started in a deterministic initial state
$X^{(N)}(0)$. Fix $t\geq 0$ and let $T_t$ be defined as in (\ref{Tt})
but with the mean-field equation (\ref{mean}) replaced by (\ref{vecmean}).
Fix $n\geq 1$ and let $I_1,\ldots,I_n$ be i.i.d.\ uniformly
distributed on $[N]$ and independent of $X^N(t)$. Then
\be\label{samp}
\Big\|\P\big[\big(X^{(N)}_{I_1}(Nt),\ldots,X^{(N)}_{I_n}(Nt)\big)\in\,\cdot\,\big]
-\underbrace{T_t(\mu^N_0)\otimes\cdots\otimes
T_t(\mu^N_0)}_{\mbox{$n$ times}}\Big\|\asto{N}0,
\ee
where $\|\,\cdot\,\|$ denotes the total variation norm, and the convergence in
(\ref{samp}) is uniform w.r.t.\ the initial state $X^{(N)}(0)$.
\ep

Proposition~\ref{P:samp} allows us to control the mean and variance of
$\mu^N_{Nt}$, which is enough to prove the convergence of $\mu^N_{Nt}$ to
$\mu_t$ for fixed times $t$. To boost this up to pathwise convergence, we use
the following lemma, that will be proved in Subsection~\ref{S:TVtight} below.

\bl[Tightness in total variation]
For\label{L:TVtight} each $N\in\N_+$ let $(X^{(N)}(t))_{t\geq 0}$ be a process
as in Theorem~\ref{T:meanlim} started in a deterministic initial state
$X^{(N)}(0)$, and let $\mu^N_t:=\mu\big\{X^{(N)}(t)\big\}$ denote the
empirical measure of $X^{(N)}(t)$. Then there exist random processes
$(\tau^N_t)_{t\geq 0}$ such that $\tau^N:\R\to\R$ is a.s.\ nondecreasing with
$\tau^N_0=0$ and
\begin{enumerate}
\item $\dis\P\big[\sup_{0\leq t\leq T}
       |\tau^N_t-t|\geq\eps\big]\asto{N}0\quad(\eps>0,\ T<\infty)$,
\item $\dis\|\mu^N_{Nt}-\mu^N_{Ns}\|\leq L(\tau^N_t-\tau^N_s)
       \quad(0\leq s\leq t)$ a.s.,
\end{enumerate}
where $\|\,\cdot\,\|$ denotes the total variation norm and
$\dis L:=\int_\Om\!\qb(\di\om)\,\la(\om)$.
\el

In Subsection~\ref{S:meanproof}, we will derive Theorem~\ref{T:meanlim}
from Proposition~\ref{P:samp}, Lemma~\ref{L:TVtight}, and some abstract
considerations.

\subsection{The state at sampled sites}\label{S:samp}

In this subsection we prove Proposition~\ref{P:samp}. We start with two
preparatory lemmas.

Let $(\Xb^N_{s,t})_{s\leq t}$ be the stochastic flow defined in
(\ref{flowdef}), where we have made the dependence on $N$ explicit.
Let $I$ be uniformly distributed on $[N]$ and independent of
$(\Xb^N_{s,t})_{s\leq t}$. For each $t\geq 0$, let $\ti M^N_t:S^{[N]}\to S$ be
defined as (note the factor $N$ rescaling the speed of time):
\be\label{MNti}
\ti M^N_t(x):=\Xb^N_{-Nt,0}(x)_I.
\ee
Let $G_t:S^{\nab\S_t}\to S$ be the random map defined in (\ref{Gtdef}),
where $\Om,\rb,\ga$ and $\kappa$ from Subsection~\ref{S:firstintro}
are defined in terms of the ``ingredients'' $\Omm,\qb$, $\ga_i[\om]$ and
$K_i(\om)$ from Subsection~\ref{S:flow}, see the proof of
Lemma~\ref{L:simpmean} in Subsection~\ref{S:formul}. Let
$(I_\ibf)_{\ibf\in\T}$ be i.i.d.\ uniformly distributed on $[N]$ and
independent of $(\nab\S_t,G_t)_{t\geq 0}$. For each $t\geq 0$, let
$M^N_t:S^{[N]}\to S$ be defined as
\be\label{MN}
M^N_t(x):=G_t\big((x_{I_\ibf})_{\ibf\in\nab\S_t}\big).
\ee
The following lemma says that for large $N$, the map in (\ref{MNti}) can be
approximated by the map in (\ref{MN}).

\bl[Coupling of maps]
For\label{L:mapcoup} each $t\geq 0$, it is possible to couple the random maps
$\ti M^N_t$ and $M^N_t$ with $N\in\N_+$ in such a way that
\be\label{mapcoup}
\P\big[\ti M^N_t\neq M^N_t\big]\asto{N}0.
\ee
\el
\bpro
The essence of the proof can be summarized as follows: since for large $N$,
sampling with or without replacement from $[N]$ is almost the same, the
genealogy of a given site is approximately given by a branching process. In
spite of this simple idea, the proof is quite long, mainly because we have to
take care of a lot of definitions, such as the way $\Om,\rb,\ga$ and $\kappa$
are defined in terms of $\Omm,\qb$, $\ga_i[\om]$ and $K_i(\om)$ in the proof
of Lemma~\ref{L:simpmean}.

We start by recalling that the random map $G_t$ from (\ref{Gtdef}) can be seen
as the concatenation of random maps assigned to the branch points of a
branching process. We then embed this branching process in the set $[N]$ and
prove that what we obtain is a good approximation for the genealogy of a given
site.

We observe that in order to construct the map $G_t:S^{\nab\S_t}\to S$
from (\ref{Gtdef}), it suffices to know
\be\label{SomS}
\big(\S_t,(\omb_\ibf)_{\ibf\in\S_t}\big),
\ee
where $\S_t$ is defined in (\ref{Stdef}). Indeed, from the information in
(\ref{SomS}) we can determine $\nab\S_t$, since
\be
\nab\S_t=\big\{\ibf k:
\ibf\in\S_t,\ 1\leq k\leq\kappa(\omb_\ibf),\ \ibf k\not\in\S_t\big\},
\ee
and the map $G_t:S^{\nab\S_t}\to S$ is obtained by concatenating the maps
$\ga[\omb_\ibf]$ with $\ibf\in\S_t$ according to the tree structure of $\S_t$.

The object in (\ref{SomS}) is in fact a Markov chain as a function of $t$.
Starting from the initial state $\S_0=\emptyset$ and $\nab\S_0=\{\wurz\}$, 
its evolution is as follows: Independently for each $\ibf\in\nab\S_t$, with
rate $|\rb|$, we add $\ibf$ to $\S_t$ and assign to it a value $\omb_\ibf$
chosen according to the probability law $|\rb|^{-1}\rb$.

We will be interested in the process in (\ref{SomS}) in the special case when
$\Om,\rb,\kappa$, and $\ga$ are defined in terms of $\Omm,\qb,\la$, $K_i(\om)$,
and $\ga_i[\om]$ as in the proof of Lemma~\ref{L:simpmean}. In this case,
elements of $\Om$ are pairs $(\om,n)$ where $\om\in\Omm$ and $1\leq
n\leq\la(\om)$, so we denote the process in (\ref{SomS}) as
\be\label{altSomS}
\big(\S_t,(\omb_\ibf,n_\ibf)_{\ibf\in\S_t}\big),
\ee
where $\omb_\ibf\in\Omm$ and $1\leq n_\ibf\leq\la(\omb_\ibf)$.
The set $\nab\S_t$ is now given by
\be\label{nabS}
\nab\S_t=\big\{\ibf k:
\ibf\in\S_t,\ 1\leq k\leq\kappa_{n_\ibf}(\omb_\ibf),\ \ibf k\not\in\S_t\big\}.
\ee
Defining $\rb$ as in (\ref{rq}), the process in (\ref{SomS}) now evolves in
such a way that independently for each $\ibf\in\nab\S_t$, with
rate $|\rb|$, we add $\ibf$ to $\S_t$ and assign values $(\omb_\ibf,n_\ibf)$
to it that are chosen according to the probability law $|\rb|^{-1}\rb$.

Let $\al\in[N]$ be fixed. Our next aim is to ``embed'' the process from
(\ref{altSomS}) in the set $[N]$, in such a way that it approximates the
genealogy of the site $\al$. To this aim, we define, for each time, a random
function $\psi^N_t:\S_t\cup\nab\S_t\to[N]$. Initially, we set
$\psi^N_0(\wurz):=\al$. We let the function $\psi^N_t$ evolve in a Markovian
way together with the process in (\ref{altSomS}) in the following way. Recall
that when we add an element $\ibf$ to $\S_t$ and assign values
$(\omb_\ibf,n_\ibf)$ to it, this element is at the same time removed from
$\nab\S_t$ and replaced by new elements $\ibf
1,\ldots,\ibf\kappa_{n_\ibf}(\omb_\ibf)$. We assign labels $\psi^N_t(\ibf k)$
$(k=1,\ldots,\kappa_{n_\ibf}(\omb_\ibf))$ to these new elements as
follows. First, we choose $(I_l)_{l=1,\ldots,\la(\omb_\ibf)}$ in such a way that
$I_{n_\ibf}:=\psi^N_t(\ibf)$ and
\be
(I_l)_{l\neq n_\ibf}\quad\mbox{are i.i.d.\ uniformly chosen from }[N],
\ee
and next, we set $\psi^N_t(\ibf k):=I_{j_k}$, where as in
(\ref{Kord}), we order the elements of
$K_{n_\ibf}(\omb_\ibf)\sub\{1,\ldots,\la(\omb_\ibf)\}$ as
\be\label{Kordn}
K_{n_\ibf}(\omb_\ibf)=\big\{j_1,\ldots,j_{\kappa_{n_\ibf}(\omb_\ibf)}\big\}
\quad\mbox{with}\quad
j_1<\cdots<j_{\kappa_{n_\ibf}(\omb_\ibf)}.
\ee
Note that this has the effect that if $n_\ibf$ is an element of
$K_{n_\ibf}(\omb_\ibf)$, say $n_\ibf=j_k$, then the corresponding element
$\ibf k$ gets the same label as $\ibf$, i.e., $\psi^N_t(\ibf
k)=\psi^N_t(\ibf)$. Otherwise, we assign new i.i.d.\ labels to all new
elements of $\nab\S_t$.

Using the function $\psi^N_t$ that embeds the process in (\ref{altSomS}) in
the set $[N]$, we define a function $\phi^N_t:S^N\to S^{\nab\S_t}$ by
\be\label{phiNt}
\phi^N_t(x)_\ibf:=x_{\psi^N_t(\ibf)}\qquad(\ibf\in\nab\S_t).
\ee
We now consider the maps
\be
\ti M^N_t(x):=\Xb_{-Nt,0}(x)_\al
\quand
M^N_t(x):=G_t\circ\phi^N_t(x)\qquad(x\in S^N),
\ee
where $\al=\psi^N_0(\wurz)\in[N]$ is the label initially assigned to the root.
We claim that
\be\label{mapcoupTV}
\big\|\P\big[\ti M^N_t\in\,\cdot\,]-\P\big[M^N_t\in\,\cdot\,]\big\|\asto{N}0,
\ee
where $\|\,\cdot\,\|$ denotes the total variation norm. In
particular, if $\al$ is chosen uniformly distributed in $[N]$ and independent
of everything else, then $(\psi^N_t(\ibf))_{\ibf\in\nab\S_t}$ are
i.i.d.\ uniformly distributed in $[N]$ and independent of the map $G_t$,
so (\ref{mapcoupTV}) implies (\ref{mapcoup}).

To prove (\ref{mapcoupTV}), we construct a process similar to the process in
(\ref{altSomS}), together with an embedding in $[N]$, that describes the true
genalogy of the site $\al$, and show that the error we make by replacing this
true genealogy by the process we had before is small. We denote this process
as
\be\label{Spsi}
\big(\ti\S_t,(\omb_\ibf,n_\ibf)_{\ibf\in\ti\S_t},\ti\psi^N_t\big).
\ee
At each time, $\nab\ti\S_t$ is defined in terms of this process in the same way
as $\nab\S_t$ is defined in (\ref{nabS}). We also define $\ti G_t$ and
$\ti\phi^N_t:S^N\to S^{\nab\ti\S_t}$ as before, i.e., $\ti G_t$ is the
concatenation of the random maps $\ga_{n_\ibf}[\om_\ibf]$ with $\ibf\in\ti\S_t$
according to the tree structure of $\ti\S_t$, and $\ti\phi^N_t:S^N\to
S^{\nab\ti\S_t}$ is defined in terms of $\ti\psi^N_t$ as in (\ref{phiNt}).

Recall that $\ti M^N_t(x)=\Xb^N_{-Nt,0}(x)_\al$. As for our previous process
we start with $\ti\S_0=\emptyset$, $\nab\ti\S_0=\{\wurz\}$, and
$\ti\psi^N(\wurz)=\al$. In Subsection~\ref{S:flow}, the
stochastic flow $(\Xb^N_{s,t})_{s\leq t}$ is constructed from a Poisson point
set $\Pi$. We will construct the process in (\ref{Spsi}) in terms of $\Pi$ in
such a way that 
\be\label{geneal}
\ti G_t\circ\ti\phi^N_t(x)=\ti M^N_t(x)\qquad(t\geq 0),
\ee
which expresses the fact that the process in (\ref{Spsi}) describes the ``true
genealogy'' of the site~$\al$.

The Poisson set $\Pi$ consists of triples $(\om,\ibf,t)$ which express the
fact that at time $t$ the random map $\vec\ga[\om]$ should be applied to the
coordinates $\ibf=(i_1,\ldots,i_{\la(\om)})$.
Note that we are interested in $\Xb^N_{-Nt,0}(x)_\al$, which means that we look at
negative times and need to rescale time by a factor $N$.
For each $(\om,\ibf,-Nt)\in\Pi$ and $\jbf\in\nab\ti\S_t$ such that
$\ti\psi^N_t(\jbf)=i_l$ for some $1\leq l\leq\la(\om)$,
we update the process in (\ref{Spsi}) as follows:
\begin{enumerate}
\item We remove $\jbf$ from $\nab\ti\S_t$ and add it to $\ti\S_t$.
\item We set $\om_\jbf:=\om$ and $n_\jbf:=l$.
\item We add $\jbf 1,\ldots,\jbf\kappa_l(\om)$ to $\ti\S_t$.
\item We define $\ti\psi^N_t(\jbf k):=i_{j_k}$ $(k=1,\ldots,\kappa_l(\om))$,
  where $K_l(\om)=\{j_1,\ldots,j_{\kappa_l(\om)}\}$ as in (\ref{Kord}).
\end{enumerate}

It is straightforward to check that these rules guarantee that (\ref{geneal})
holds and hence the process in (\ref{Spsi}) describes the true genealogy of
the site $\al$. As some more explanation, we can add the following: we follow
a site $\bet$ back in time till the first time when a map is applied that has
the possibility to change the value of $\bet$. From that moment on, we follow
back all sites that are relevant for the outcome of the map at $\bet$, and we
number them according to the convention in (\ref{Kord}). This defines a family
structure, i.e., $\ibf=i_1i_2i_3$ is the $i_3$-th child of the $i_2$-th
child of the $i_1$-th child of the original site $\al$. The map $\ti\psi^N_t$
applied to $\ibf$ tells us where this ancestor lives in the set $[N]$. There
may be some overlap, i.e., it is possible that
$\ti\psi^N_t(\ibf)=\ti\psi^N_t(\jbf)$ for some
$\ibf,\jbf\in\ti\S_t\cup\nab\ti\S_t$. For $\ibf,\jbf\in\nab\ti\S_t$, however,
the probability that two ancestors live at the same site in $[N]$ tends to
zero as $N\to\infty$, as we will see in a moment.

In view of (\ref{geneal}), to prove (\ref{mapcoupTV}), it suffices to prove
that the Markov process in (\ref{Spsi}) is close in total variation distance
to the process with $\ti\S_t$ and $\ti\psi^N_t$ replaced by $\S_t$ and
$\psi^N_t$. Since the latter process is nonexplosive by (\ref{EpaTt}), it
suffices to prove convergence for the processes stopped at the first time when
the cardinality of $\ti\nab\S_t$ resp.\ $\nab\S_t$ exceeds a certain value,
and then at the end send this value to infinity. We will prove convergence of
the stopped processes in a number of steps, by making small changes in the
jump rates. Here we use the fact that if the transition kernels of two
continuous-time Markov chains are close in total variation norm, uniformly in
the starting point, then by standard arguments the two processes can be
coupled so that their laws at fixed time are close in total variation norm.

Let $\ti\Psi^N_t:=\ti\psi^N_t(\nab\ti\S_t)$ denote the image of $\nab\ti\S_t$
under the map $\ti\psi^N_t$. As a first step, we change the dynamics of the
(stopped) process from (\ref{Spsi}) in such a way that elements
$(\om,\ibf,-Nt)\in\Pi$ have no effect if $\{i_1,\ldots,i_{\la(\om)}\}$
intersects $\ti\Psi^N_t$ in more that one point. Then the modified process is
still Markovian; we claim the change in jump rates compared to the original
process is of order $N^{-1}$. Indeed, for fixed $l$, if $i_1,\ldots,i_l$ are
chosen uniformly without replacement from $[N]$, then the probability that
one, resp.\ two or more of them lie in a set $A$ of fixed cardinality is of
order $N^{-1}$ resp.\ $N^{-2}$ as $N\to\infty$. Taking into
account the fact that we rescale time by a factor $N$, as well as the
summability condition (\ref{rmom2})~(i), this translates into a change in 
jump rates of order $N^{-1}$ for the modified process, stopped at the first
time when the cardinality of $\nab\ti\S_t$ exceeds a fixed value.

Recall that by (\ref{geneal}), $\Xb^N_{-Nt,0}(x)_\al$ is a function only of
$(x_\bet)_{\bet\in\ti\Psi^N_t}$. The modified process we have just constructed
has the property that $\ti\psi^N_t:\nab\ti\S_t\to\ti\Psi^N_t$ is a bijection,
i.e., each element $\bet\in\ti\Psi^N_t$ corresponds only to a single place
$(\psi^N_t)^{-1}(\bet)$ in the family tree. The dynamics of the modified process
can be described as follows:
\begin{enumerate}
\item Independently for each $\bet\in\ti\Psi^N_t$, with rates described by the
  measure $\rb$ from (\ref{rq}), we choose a pair $(\om,n)$ with $1\leq
  n\leq\la(\om)$.
\item If $\la(\om)>N$, we do nothing.
\item Otherwise, we choose $(\bet'_k)_{k=1,\ldots,\la(\om)}$ such that
  $\bet'_n:=\bet$ and $(\bet'_k)_{k\neq n}$ are drawn from $[N]\beh\{\bet\}$
  without replacement.
\item If some of the $(\bet'_k)_{k\neq n}$ are elements of $\ti\Psi^N_t$, we
  do nothing.
\item Otherwise, we remove $\bet$ from $\ti\Psi^N_t$ and add $(\bet_{j_k})_{1\leq
  k\leq\kappa_n(\om)}$ to $\ti\Psi^N_t$, where
  $K_n(\om)=\{j_1,\ldots,j_{\kappa_n(\om)}\}$ with
  $j_1<\cdots<j_{\kappa_n(\om)}$.
\item If $\jbf=(\ti\psi^N_{t-})^{-1}(\bet)$ is the place of $\bet$ in the
  family tree immediately prior to time $t$, then we assign to each new element
  of $\ti\Psi^N_t$ a place in the family tree by setting
  $(\ti\psi^N_{t-})^{-1}(\bet_{j_k}):=\jbf k$.
\end{enumerate}
Note that the measure $\rb$ from (\ref{rq}) occurs naturally here, since
each $\la(\om)$-tuple of sites in $[N]$ can contain a given site $\bet$ in
$\la(\om)$ different ways, as its 1st, 2nd,\ldots, $\la(\om)$-th member.

Removing the restrictions in points (ii) and (iv) above, and performing
sampling without replacement instead of sampling with replacement in point
(iii), we only make changes in the transition rates of order $N^{-1}$, and
arrive at a process whose family tree evolves as the process in
(\ref{altSomS}) and where to new members of the family tree, sites in $[N]$ are
assigned chosen uniformly with replacement, as described by the process
$\psi^N_t$.
\epro

In the proof of Lemma~\ref{L:mapcoup}, we have seen that in the mean-field
limit $N\to\infty$, the genealogy of a single site can be approximated by a
branching process with random maps attached to its branch points. Similarly,
the genealogy of $n$ randomly chosen sites can be approximated by $n$
independent branching processes, which leads to the following extension of
Lemma~\ref{L:mapcoup}.

\bl[The genealogy of multiple sites]
Let\label{L:geneal} $(\Xb^N_{s,t})_{s\leq t}$ be the stochastic flow defined in
(\ref{flowdef}) and let $I_1,\ldots,I_n$ be i.i.d.\ uniformly distributed on
$[N]$, independent of $(\Xb^N_{s,t})_{s\leq t}$. Let $\Om,\rb,\ga$ and
$\kappa$ be defined in terms of $\Omm,\qb$, $\ga_i[\om]$ and $K_i(\om)$ as in
the proof of Lemma~\ref{L:simpmean}. Fix $t\geq 0$ and let
$(\nab\S^i_t,G^i_t)$ $(i=1,\ldots,n)$ be i.i.d.\ copies of the random set and
map defined in (\ref{Stdef}) and (\ref{Gtdef}). Conditional on
$(\nab\S^i_t,G^i_t)_{i=1,\ldots,n}$, let
$(I^i_\jbf)^{i=1,\ldots,n}_{\jbf\in\nab\S^i_t}$ be i.i.d.\ uniformly distributed on
$[N]$. Define $\ti M^N_t:S^N\to S^n$ and $M^N_t:S^N\to S^n$ by
\be\label{vecM}
\ti M^N_t(x)_i:=\Xb^N_{-Nt,0}(x)_{I_i}
\quand
M^N_t(x)_i:=G^i_t\big((x_{I^i_\jbf})_{\jbf\in\nab\S^i_t}\big).
\ee
Then $\ti M^N_t$ and $M^N_t$ can be coupled such that $\dis\P\big[\ti
  M^N_t\neq M^N_t\big]\asto{N}0$.
\el
\bpro
The proof is the same as the proof of Lemma~\ref{L:mapcoup}, except that
instead of following back the genealogy of one site, one follows the
genealogies of $n$ sites. By the same arguments as given in the proof of
Lemma~\ref{L:mapcoup}, when $N$ is large, with high probability, the
genealogies do not intersect, and hence can be approximated by independent
branching processes. Although writing down all objects involved is
notationally complicated, no new ideas are needed so we omit the details.
\epro

\bpro[of Proposition~\ref{P:samp}]
Let $x:=X^{(N)}(0)$ be the (deterministic) initial state and using notation as
in (\ref{empir0}) let $\mu^N_0=\mu\{x\}$ denote its empirical measure.
Define maps $\ti M^N_t$ and $M^N_t$ as in Lemma~\ref{L:geneal}.
Then $\big(X^{(N)}_{I_1}(Nt),\ldots,X^{(N)}_{I_n}(Nt)\big)$ 
has law $\ti M^N_t(x)$ while the coordinates of
$M^N_t(x)$ are i.i.d.\ with a law that by Theorem~\ref{T:meanrep} equals
$T_t(\mu^N_0)$. In view of this, the claim follows from Lemma~\ref{L:geneal}.
\epro

\subsection{Tightness in total variation}\label{S:TVtight}

In this subsection we prove Lemma~\ref{L:TVtight}.\med

\bpro[of Lemma~\ref{L:TVtight}]
The process $(X^{(N)}(t))_{t\geq 0}$ is defined in (\ref{XNdef}) in terms of a
stochastic flow which is in turn defined in terms of a Poisson set $\Pi$.
Elements of $\Pi$ are triples $(\om,\ibf,s)$ which tell us that at time $s$
the map $\vec\ga[\om]$ should be applied to the coordinates
$\ibf=(i_1,\ldots,i_{\la(\om)})$. We let
\be
\tau^N_t:=\frac{1}{NL}\sum_{(\om,\ibf,s)\in\Pi:\,0<s\leq Nt}\la(\om)\qquad(t\geq 0),
\ee
where $L:=\int_\Om\!\qb(\di\om)\,\la(\om)$, which is finite by (\ref{rmom2}).
Then (i) follows from a functional law of large numbers. Since for any $s\leq
t$, the fraction of sites in $[N]$ that changes its type is bounded from above
by $L(\tau^N_t-\tau^N_s)$, in view of (\ref{TVcoup}), we obtain also (ii).
\epro

\subsection{Convergence to the mean-field equation}\label{S:meanproof}

In this subsection, we prove Theorem~\ref{T:meanlim}. The proof is split into a
number of lemmas. We start by proving convergence at fixed times. This part of
the proof is based on Proposition~\ref{P:samp}. At the end of the proof, we
use Lemma~\ref{L:TVtight} to obtain pathwise convergence.

\bl[Expectation of test functions]
Let\label{L:momcal} $\Omm,\qb,\la$, and $\vec\ga$ be as in
Subsection~\ref{S:flow}, and assume (\ref{rmom2}). Let $(T_t)_{t\geq 0}$
denote the semigroup defined as in (\ref{Tt}) but with the mean-field equation
(\ref{mean}) replaced by (\ref{vecmean}). For each $N\in\N_+$, let
$(X^{(N)}(t))_{t\geq 0}$ be Markov processes with state space $S^N$ as defined
in (\ref{XNdef}), and let $\mu^N_t=\mu\{X^{(N)}(t)\}$ denote their associated
empirical measures. Then
\be\label{momcal}
\sup_{|\phi|\leq 1}
\P\big[\big|\li\mu^N_{Nt},\phi\re-\li T_t(\mu^N_0),\phi\re\big|\geq\eps\big]
\asto{N}0\qquad(\eps>0,\ t\geq 0),
\ee
where the supremum runs over all measurable functions $\phi:S\to[-1,1]$.
\el
\bpro
Fix $t\geq 0$. Let $\phi:S\to[-1,1]$ be measurable. Let $I_1$ and $I_2$ be
uniformly distributed on $[N]$ and independent of each other and of $X^N(t)$.
Since
\be
\li\mu^N_{Nt},\phi\re=\frac{1}{N}\sum_{i=1}^N\phi\big(X^{(N)}_i(Nt)\big),
\ee
we see that
\bc
\dis\E\big[\li\mu^N_{Nt},\phi\re\big]
&=&\dis\E\big[\phi\big(X^{(N)}_{I_1}(Nt)\big)\big],\\[5pt]
\dis\E\big[\li\mu^N_{Nt},\phi\re^2\big]
&=&\dis\E\big[\phi\big(X^{(N)}_{I_1}(Nt)\big)\phi\big(X^{(N)}_{I_2}(Nt)\big)\big].
\ec
Assume for the moment that $X^{(N)}(0)$ is deterministic. Then applying
Proposition~\ref{P:samp} with $n=1,2$ we find that
\be\ba{r}
\dis\sup_{|\phi|\leq 1}
\Big|\E\big[\li\mu^N_{Nt},\phi\re\big]-\li T_t(\mu^N_0),\phi\re\Big|
\asto{N}0,\\[5pt]
\dis\sup_{|\phi|\leq 1}
\Big|\E\big[\li\mu^N_{Nt},\phi\re^2\big]-\li T_t(\mu^N_0),\phi\re^2\Big|
\asto{N}0,
\ec
where we take the supremum over all measurable $\phi:S\to[-1,1]$. It follows
that
\be
\sup_{|\phi|\leq 1}\var\big(\li\mu^N_{Nt},\phi\re\big)
=\sup_{|\phi|\leq 1}\Big(
\E\big[\li\mu^N_{Nt},\phi\re^2\big]-\E\big[\li\mu^N_{Nt},\phi\re\big]^2\Big)
\asto{N}0,
\ee
and hence (\ref{momcal}) follows by Chebyshev's inequality. To obtain
(\ref{momcal}) more generally when $X^{(N)}(0)$ is random, we condition on the
initial state to get, for each $\eps>0$ and measurable $\psi:S\to[-1,1]$.
\be\ba{l}
\dis\int\P[X^{(N)}(0)\in\di x]\,
\P\big[\big|\li\mu^N_{Nt},\psi\re-\li T_t(\mu^N_0),\psi\re\big|\geq\eps
\,\big|\,X^{(N)}(0)=x\big]\\[5pt]
\dis\quad\leq\int\P[X^{(N)}(0)\in\di x]\,\sup_{|\phi|\leq 1}
\P\big[\big|\li\mu^N_{Nt},\phi\re-\li T_t(\mu^N_0),\phi\re\big|\geq\eps
\,\big|\,X^{(N)}(0)=x\big].
\ec
Since the integrand on the right-hand side does not depend on $\psi$ and
tends to zero in a bounded pointwise way as a function of $x\in S^N$,
(\ref{momcal}) follows.
\epro

Our next aim is to prove that if in addition to the assumptions of
Lemma~\ref{L:momcal}, condition (i) or (ii) of Theorem~\ref{T:meanlim} is
satisfied, then
\be\label{fixconv}
\P\big[d(\mu^N_{Nt},T_t(\mu_0))\geq\eps\big]\asto{N}0
\qquad(\eps>0,\ t\geq 0),
\ee
where $d$ is any metric on $\Pc(S)$ that generates the topology of weak
convergence. Applying the following well-known fact to the Polish space
$\Pc(S)$, we see that if (\ref{fixconv}) holds for one such metric, then it
holds for all of them.

\bl[Convergence in probability]
Let\label{L:toP} $X_n$ be random variables taking values in a Polish space
$S$, let $x\in S$ be deterministic, and let $d$ be a metric generating the
topology on $S$. Then one has
\be\label{toP}
\P\big[d(X_n,x)\geq\eps\big]\asto{n}0\qquad(\eps>0),
\ee
if and only if
\be\label{towk}
\P\big[X_n\in\,\cdot\,]\Asto{n}\de_x,
\ee
where $\Rightarrow$ denotes weak convergence of probability measures on $S$.
\el
\bpro
It is easy to see that (\ref{toP}) implies $\E[\phi(X_n)]\to\phi(x)$ for all
bounded continuous $\phi:S\to\R$, so (\ref{toP}) implies
(\ref{towk}). Conversely, if (\ref{towk}) holds, then by Skorohod's
representation theorem it is possible to couple the random variables $X_n$
such that $X_n\to x$ a.s., which implies (\ref{toP}). 
\epro

The following lemma gives sufficient conditions for the type of convergence of
(\ref{fixconv}).

\bl[Convergence to a deterministic measure]
Let\label{L:measconv} $S$ be a Polish space, let $\mu\in\Pc(S)$ be
deterministic, and let $\mu^N$ be random variables with values in
$\Pc(S)$. Let $d$ be a metric on $\Pc(S)$ generating the topology of weak
convergence. Then the following conditions are equivalent.
\begin{enumerate}
\item $\dis\P\big[d(\mu^N,\mu)\geq\eps\big]\asto{N}0$ for all $\eps>0$.
\item $\dis\P\big[\big|\li\mu^N,\phi\re-\li\mu,\phi\re\big|\geq\eps\big]
       \asto{N}0$ for all $\eps>0$ and bounded continuous $\phi:S\to\R$.
\item $\dis\E\big[\prod_{i=1}^n\li\mu^N,\phi_i\re\big]\asto{N}
       \prod_{i=1}^n\li\mu,\phi_i\re$ for all bounded continuous
       functions $\phi_1,\ldots,\phi_n$ $(n\geq 1)$.
\end{enumerate} 
\el
\bpro
We equip $\Pc(S)$ with the topology of weak convergence, making it into a
Polish space. Then by Lemma~\ref{L:toP}, condition~(i) is equivalent to
\begin{itemize}
\item[{\rm(i)'}] $\dis\P[\mu^N\in\,\cdot\,]\Asto{N}\de_\mu$.
\end{itemize}
We will prove (i)'$\volgt$(ii)$\volgt$(iii)$\volgt$(i)'.

(i)'$\volgt$(ii). By Skorohod's representation theorem, (i)' implies that the
$\mu^N$ can be coupled such that $\mu^N\Asto{N}\mu$ a.s., which implies~(ii).

(ii)$\volgt$(iii). Without loss of generality we may assume that the
$\phi_i$'s take values in $[-1,1]$. Since the function
$(x_1,\ldots,x_n)\mapsto\prod_{i=1}^nx_i$ is continuous, (ii) implies that
\be
\P\big[\big|\prod_{i=1}^n\li\mu^N,\phi_i\re-\prod_{i=1}^n\li\mu,\phi_i\re\big|
\geq\eps\big]\asto{N}0
\ee
for all $\eps>0$ and bounded continuous functions $\phi_1,\ldots,\phi_n$.
Since moreover $|\prod_{i=1}^n\li\mu^N,\phi_i\re|\leq 1$, this implies (iii).

(iii)$\volgt$(i)'. Since $S$ is Polish, it has a metrizable compactification,
i.e., there exists a compact metrizable space $\ov S$ such that $S$ is a dense
subset of $\ov S$ and the topology on $S$ is the induced topology from $\ov
S$ \cite[Theorem~6.3]{Cho69}. It is known that this implies that $S$ is a
$G_\de$-subset of $\ov S$ \cite[\Parag 6 No.~1, Theorem.~1]{Bou58}. In
particular, $S$ is a Borel measurable subset of $\ov S$ and we can identify
$\Pc(S)$ with the space of probability measures on $\ov S$ that are
concentrated on $S$. If we equip $\Pc(\ov S)$ with the topology of weak
convergence, then the induced topology on $\Pc(S)$ is also the topology of weak
convergence (this follows, e.g., from \cite[Thm~3.3.1]{EK86}), and in fact
$\Pc(\ov S)$ (being compact by Prohorov's theorem) is a metrizable
compactification of $\Pc(S)$.

We view $\mu^N$ and $\mu$ as probability measures on $\ov S$. Since $\ov
S$ is compact, so are $\Pc(\ov S)$ and $\Pc(\Pc(\ov S))$, so by going to a
subsequence if necessary, we can assume that the laws $\P[\mu^N\in\,\cdot\,]$
converge weakly to some limit $\rho\in\Pc(\Pc(\ov S))$. Since the
restriction to $S$ of a continuous function $\phi:\ov S\to\R$ is a bounded
continuous function on $S$, condition~(ii) implies that 
\be\label{rhochar}
\int_{\Pc(\ov S)}\!\rho(\di\nu)\,\prod_{i=1}^n\li\nu,\phi_i\re
=\prod_{i=1}^n\li\mu,\phi_i\re
\ee
for general $n\geq 1$ and continuous functions $\phi_i:\ov S\to\R$
$(i=1,\ldots,n)$. By the Stone-Weierstrass theorem, the linear span of
functions of the form $\nu\mapsto\prod_{i=1}^n\li\mu,\phi_i\re$ is dense in
the space of continuous functions on $\Pc(\ov S)$, and hence (\ref{rhochar})
implies $\rho=\de_\mu$.
\epro

We now prove (\ref{fixconv}) under either of the conditions (i) and (ii) of
Theorem~\ref{T:meanlim}.

\bl[Continuity argument]
In\label{L:coarg} addition to the assumptions of Lemma~\ref{L:momcal}, assume
that condition~(i) of Theorem~\ref{T:meanlim} is satisfied. Then (\ref{fixconv})
holds.
\el
\bpro
Fix $t\geq 0$. In view of Lemma~\ref{L:measconv}~(ii), it suffices to show that
\be
\P\big[\big|\li\mu^N_{Nt},\phi\re-\li T_t(\mu_0),\phi\re\big|\geq\eps\big]
\asto{N}0\qquad(\eps>0)
\ee
for any bounded continuous $\phi:S\to\R$. By Lemma~\ref{L:momcal}, it suffices
to show that
\be\label{initcond}
\P\big[\big|\li T_t(\mu^N_0),\phi\re-\li T_t(\mu_0),\phi\re\big|\geq\eps\big]
\asto{N}0\qquad(\eps>0).
\ee
By the second part of condition~(i), Lemma~\ref{L:simpmean}, and
Proposition~\ref{P:gacon}, the operator $T_t$ is continuous w.r.t.\ weak
convergence. In view of this, (\ref{initcond}) is implied by the first part of
condition~(i).
\epro

\bl[Moment argument]
In\label{L:moarg} addition to the assumptions of Lemma~\ref{L:momcal}, assume
that condition~(ii) of Theorem~\ref{T:meanlim} is satisfied. Then
(\ref{fixconv}) holds.
\el
\bpro
Fix $t\geq 0$. In view of Lemma~\ref{L:measconv}~(iii), it suffices to show that
\be\label{suf}
\E\big[\prod_{i=1}^n\li\mu^N_{Nt},\phi_i\re\big]
\asto{N}\prod_{i=1}^n\li T_t(\mu_0),\phi_i\re
\ee
for all $n\geq 1$ and bounded continuous functions $\phi_i:S\to\R$,
$i=1,\ldots,n$. Without loss of generality we may assume that the $\phi_i$'s
take values in $[-1,1]$. Let $X^{(N)}(t)$ be as in Theorem~\ref{T:meanlim} and
let $I_1,\ldots,I_n$ be i.i.d.\ uniformly distributed on $[N]$ and independent
of $X^N(t)$. Then
\be\label{prodli}
\E\big[\prod_{i=1}^n\li\mu^N_{Nt},\phi_i\re\big]
=\E\big[\prod_{i=1}^n\phi_i\big(X^{(N)}_{I_i}(N t)\big)\big].
\ee
By Proposition~\ref{P:samp} applied to the process conditioned on
$X^{(N)}(0)$, there exist $\eps_N\to 0$ such that
\be
\Big\|\P\Big[\big(X^{(N)}_{I_1}(Nt),\ldots,X^{(N)}_{I_n}(Nt)\big)\in\,\cdot\,
\big|\,X^{(N)}(0)\Big]-T_t(\mu^N_0)^{\otimes n}\Big\|\leq\eps_N.
\ee
In view of (\ref{TVdef}), it follows that
\be
\Big|\E\Big[\prod_{i=1}^n\phi_i\big(X^{(N)}_{I_i}(t)\big)\Big|\,X^{(N)}(0)\Big]
-\prod_{i=1}^n\li T_t(\mu^N_0),\phi_i\re\Big|\leq\eps_N.
\ee
Combining this with (\ref{prodli}), taking the expectation, we obtain that
\be\label{erbd}
\Big|\E\big[\prod_{i=1}^n\li\mu^N_{Nt},\phi_i\re\big]
-\E\big[\prod_{i=1}^n\li T_t(\mu^N_0),\phi_i\re\big]\Big|\leq\eps_N.
\ee
In view of this, to prove (\ref{suf}), it suffices to show that
\be\label{TmuNT}
\E\Big[\prod_{i=1}^n\li T_t(\mu^N_0),\phi_i\re\Big]
\asto{N}\prod_{i=1}^n\li T_t(\mu_0),\phi_i\re.
\ee
If $\mu\in\Pc(S)$ is deterministic, then Theorem~\ref{T:meanrep} tells us that
\be\label{repconseq}
\li T_t(\mu),\phi\re
=\E\big[\phi\circ G_t\big((X_\jbf)_{\jbf\in\nab\S_t}\big)\big],
\ee
where $\nab\S_t$ and $G_t$ are as in (\ref{Stdef}) and (\ref{Gtdef}) and
$(X_\jbf)_{\jbf\in\T}$ are i.i.d.\ with law $\mu$. Conditional on $\mu^N_0$,
let $(X^i_\jbf)^{i=1,\ldots,n}_{\jbf\in\T}$ be i.i.d.\ with common law $\mu^N_0$.
Let $(\nab\S^i_t,G^i_t)$ $(i=1,\ldots,n)$ be i.i.d.\ and distributed as the
random variables in (\ref{Stdef}) and (\ref{Gtdef}), independent of
$\mu^N_0$ and  $(X^i_\jbf)^{i=1,\ldots,n}_{\jbf\in\T}$. Then (\ref{repconseq})
implies that
\be\label{prorep}
\E\Big[\prod_{i=1}^n\li T_t(\mu^N_0),\phi_i\re\Big]
=\E\Big[\prod_{i=1}^n\phi_i\circ G^i_t\big((X^i_\jbf)_{\jbf\in\nab\S^i_t}\big)\Big].
\ee
If we replace the expectation on the right-hand side by a conditional
expectation given $(\nab\S^i_t,G^i_t)_{i=1,\ldots,n}$, then this is the
integral of a measurable $[-1,1]$-valued function with respect to the
expectation of a product measure of the form $(\mu^N_0)^{\otimes m}$, where
$m=\sum_{i=1}^n|\nab\S^i_t|$. Condition~(ii) of Theorem~\ref{T:meanlim} allows
us to replace the integral w.r.t.\ $\E[(\mu^N_0)^{\otimes m}]$ by
the integral w.r.t.\ $\mu_0^{\otimes m}$ at the cost of a small error.
Thus,
\be\ba{l}
\dis\E\Big[\prod_{i=1}^n\phi_i\circ G^i_t\big((X^i_\jbf)_{\jbf\in\nab\S^i_t}\big)
\,\Big|\,(\nab\S^i_t,G^i_t)_{i=1,\ldots,n}\Big]\\[5pt]
\dis\quad=
\E\Big[\prod_{i=1}^n\phi_i\circ G^i_t\big((\ti X^i_\jbf)_{\jbf\in\nab\S^i_t}\big)
\,\Big|\,(\nab\S^i_t,G^i_t)_{i=1,\ldots,n}\Big]+R^N,
\ec
where the $(\ti X^i_\jbf)^{i=1,\ldots,n}_{\jbf\in\T}$ are i.i.d.\ with common
law $\mu_0$ and independent of $(\nab\S^i_t,G^i_t)_{i=1,\ldots,n}$,
and $R^N$ is a random error term that by condition~(ii) can be estimated as
\be
|R^N|\leq\eps^N_m\quad\mbox{with}\quad m=\sum_{i=1}^n|\nab\S^i_t|,
\ee
where $\lim_{N\to\infty}\eps^N_m=0$ for each $m$. Note that moreover
$|R^N|\leq 2$ since the $\phi_i$'s take values in $[-1,1]$. Integrating over
the randomness of $(\nab\S^i_t,G^i_t)_{i=1,\ldots,n}$, using bounded
  convergence, (\ref{repconseq}) and (\ref{prorep}), (\ref{TmuNT})
  follows.
\epro

With Lemmas~\ref{L:coarg} and \ref{L:moarg} proved, most of the work needed
for proving Theorem~\ref{T:meanlim} is done. The only remaining task is to
improve the convergence at fixed times in (\ref{fixconv}) to pathwise
convergence as in (\ref{convmu2}). Our first aim is to show that the condition
(\ref{convmu2}) does not depend on the choice of the metric $d$. This follows
from the following lemma, applied to the Polish space $\Pc(S)$.

\bl[Convergence in path space]
Let\label{L:pathconv} $S$ be a Polish space and let $d$ be a metric generating
the topology on $S$. Let $\Di_S\half$ be the space of cadlag functions
$x:\half\to S$, equipped with the Skorohod topology. Let $X_n=(X_n(t))_{t\geq
  0}$ be random variables with values in $\Di_S\half$ and let $x:\half\to S$
be a continuous function. Then one has
\be\label{locuni}
\P\big[\sup_{0\leq t\leq T}d\big(X_n(t),x(t)\big)\geq\eps\big]\asto{n}0
\qquad(\eps>0,\ T<\infty),
\ee
if and only if
\be\label{Skorlaw}
\P\big[X_n\in\,\cdot\,]\Asto{n}\de_x,
\ee
where $\Rightarrow$ denotes weak convergence of probability measures on
$\Di_S\half$.
\el
\bpro
It is well-known that $\Di_S\half$ is a Polish space \cite[Sect.~3.5]{EK86}.
Let $d_{\rm S}$ be the metric generating the topology on $\Di_S\half$ defined
in \cite[(5.2) of Chapter~3]{EK86}. Then it is easy to see that for all
$\de>0$ there exist $\eps>0$ and $T<\infty$ such that
\be
\sup_{0\leq t\leq T}d\big(y(t),x(t)\big)\leq\eps
\quad\mbox{implies}\quad
d_{\rm S}(y,x)\leq\de\qquad(x,y\in\Di_S\half).
\ee
In view of this, (\ref{locuni}) implies
\be\label{locuni2}
\P\big[d_{\rm S}\big(X_n,x\big)\geq\eps\big]\asto{n}0
\qquad(\eps>0),
\ee
which by Lemma~\ref{L:toP} implies (\ref{Skorlaw}).
Conversely, if (\ref{Skorlaw}) holds, then by Skorohod's
representation theorem it is possible to couple the random variables $X_n$
such that $d_{\rm S}(X_n,x)\to 0$ a.s. By the continuity of $x$ and
\cite[Lemma~3.10.1]{EK86}, this implies that
\be
\sup_{0\leq t\leq T}d\big(X_n(t),x(t)\big)\asto{n}0\quad{\rm a.s.}
\qquad(T<\infty),
\ee
which implies (\ref{locuni}).
\epro

Before the proof of Theorem~\ref{T:meanlim} we need one more lemma.

\bl[Weak convergence and convergence in total variation norm]
Let\label{L:weakstrong} $S$ be a Polish space. Then there exists a metric $d$
on $\Pc(S)$ such that $d$ generates the topology of weak convergence and
$d(\mu,\nu)\leq\|\mu-\nu\|$ $(\mu,\nu\in\Pc(S))$, where $\|\,\cdot\,\|$ denotes
the total variation norm.
\el
\bpro
Let $r$ be a metric generating the topology on $S$. Replacing $r(x,y)$ by
$r(x,y)\wedge 1$ if necessary we can assume without loss of generality that
$r\leq 1$. Let $\Li$ be the space of all functions $\phi:S\to\R$ such that
$|\phi(x)-\phi(y)|\leq r(x,y)$ $(x,y\in S)$, i.e., these are Lipschitz
continuous functions with Lipschitz constant $\leq 1$. Then
\be\label{Wass}
d(\mu,\nu):=\sup_{\phi\in\Li}\Big|\int\phi\,\di\mu-\int\phi\,\di\nu\Big|
\ee
is the 1-Wasserstein metric on $\Pc(S)$, which is known to generate the
topology of weak convergence. Let $\Li':=\{\phi\in\Li:\sup_{x\in
  S}|\phi(x)|\leq 1\}$. Since $r\leq 1$, each function $\phi\in\Li$ can be
written as $\phi=\ha\phi'+c$ with $\phi'\in\Li'$ and $c\in\R$.
In view of this and (\ref{TVdef}),
\be
d(\mu,\nu)=\ha\sup_{\phi\in\Li'}\Big|\int\phi\,\di\mu-\int\phi\,\di\nu\Big|
\leq\ha\sup_{|\phi|\leq 1}\Big|\int\phi\,\di\mu-\int\phi\,\di\nu\Big|
=\|\mu-\nu\|.
\ee
\epro

\bpro[of Theorem~\ref{T:meanlim}]
Lemmas~\ref{L:coarg} and \ref{L:moarg} show that either of the conditions (i)
and (ii) implies (\ref{fixconv}). We will use Lemma~\ref{L:TVtight} to improve 
(\ref{fixconv}) to pathwise convergence as in (\ref{convmu2}). By
Lemma~\ref{L:pathconv} it suffices to prove (\ref{convmu2}) for one particular
metric $d$ on $\Pc(S)$ that generates the topology of weak convergence.
We will choose a metric $d$ as in Lemma~\ref{L:weakstrong}.

Set $\mu_t:=T_t(\mu_0)$ $(t\geq 0)$ denote the solution to the mean-field
equation (\ref{vecmean}) with initial state $\mu_0$. Lemma~\ref{L:TVtight}
implies that
\be
\P\big[\|\mu^N_{Ns}-\mu^N_{Nt}\|\geq L|t-s|+\eps\big]\asto{N}0
\qquad(\eps>0,\ s,t\geq 0).
\ee
Taking the limit $N\to\infty$, using the fact that
$d(\mu,\nu)\leq\|\mu-\nu\|$ and (\ref{fixconv}), it follows that
\be\label{muLip}
d(\mu_s,\mu_t)\leq L|t-s|\qquad(s,t\geq 0).
\ee
Since for any $s,t\geq 0$,
\be
d(\mu^N_{Ns},\mu_s)
\leq\|\mu^N_{Ns}-\mu^N_{Nt}\|+d(\mu^N_{Nt},\mu_t)+d(\mu_t,\mu_s),
\ee
using Lemma~\ref{L:TVtight}, (\ref{fixconv}), and (\ref{muLip}), we see that
for each $T>0$ and $t\in[0,T]$,
\be
\P\big[\exists s\in[0,T]\mbox{ s.t.\ }d(\mu^N_{Ns},\mu_s)\geq 2L|t-s|+\eps\big]
\asto{N}0\qquad(\eps>0).
\ee
Combining this with the fact that by (\ref{fixconv}), for any $n\geq 1$,
\be
\P\big[\sup_{k=0,\ldots,n}d(\mu^N_{N(k/n)T},\mu_{(k/n)T})\geq\eps\big]\asto{N}0
\qquad(\eps>0),
\ee
we find that
\be
\P\big[\sup_{t\in[0,T]}d(\mu^N_{Nt},\mu_{t})\geq L/n+\eps\big]\asto{N}0
\qquad(\eps>0,\ n\geq 1).
\ee
Since $\eps$ and $n$ are arbitrary, this implies (\ref{convmu2}).
\epro

\section{Recursive Tree Processes}\label{S:RTP}

In this section, we prove our main results about RTPs with continuous time.
For completeness, we also prove Lemma~\ref{L:RTP} which deals with discrete
time and says that each solution to the RDE (\ref{RDE}) gives
rise to an RTP. This is done in Subsection~\ref{S:RTPconstr}

Our basic results about continuous-time RTPs are Lemma~\ref{L:consis}
and Proposition~\ref{P:contRTP}. Lemma~\ref{L:consis} describes the evolution
of the law of the process
\be\label{contreeMark}
\{X_\ibf:\ibf\in\nab\St_{t-s}\}_{0\leq s\leq t},
\ee
that is constructed by assigning independent values $X_\ibf$ to
elements $\ibf\in\nab\St_t$ and then calculating backwards.
Proposition~\ref{P:contRTP} says that adding exponential lifetimes to the
elements of an RTP yields a stationary version of the process in
(\ref{contreeMark}). These results are proved in Subsection~\ref{S:contRTP}.

In Subsection~\ref{S:endlev}, we prove continuous-time analogues of known
discrete-time results related to endogeny. Following \cite{AB05},
Theorem~\ref{T:bivar2} links the $n$-variate mean-field equation to endogeny, while
Propositions~\ref{P:levmom} and \ref{P:minconv} are concerned with the
higher-level mean-field equation, and closely follow ideas from \cite{MSS18}.

\subsection{Construction of RTPs}\label{S:RTPconstr}

\bpro[of Lemma~\ref{L:RTP}]
For each finite subtree $\U\sub\ovT$ that contains the root, we can construct
random variables $(\omb_\ibf)_{\ibf\in\U}$ and $(X_\ibf)_{\ibf\in\U\cup\poa\U}$ such
that the $(\omb_\ibf)_{\ibf\in\U}$ are independent with common law
$|\rb|^{-1}\rb$, the $(X_\ibf)_{\ibf\in\poa\U}$ are i.i.d.\ with common law
$\nu$ and independent of the $(\omb_\ibf)_{\ibf\in\U}$, and the
$(X_\ibf)_{\ibf\in\U}$ are inductively defined by
\be\label{Xind}
X_\ibf=\ga[\omb_\ibf]\big(X_{\ibf1},\ldots,X_{\ibf\kappa(\omb_\ibf)}\big)
\qquad(\ibf\in\U).
\ee
The joint law of $(\omb_\ibf)_{\ibf\in\U}$ and $(X_\ibf)_{\ibf\in\U\cup\poa\U}$
is a probability law $\P_\U$ on $\Om^\U\times S^{\U\cup\poa\U}$.
Since $\Om$ and $S$ are Polish spaces, we can apply Kolmogorov's extension
theorem. The statement of the lemma then follows provided we can show that
the laws $\P_\U$ are consistent in the sense that if
$\V\sub\U$ is another subtree that contains the root, then the projection of 
$\P_\U$ on $\Om^\V\times S^{\V\cup\poa\V}$ equals $\P_\V$. It suffices to
prove this when $\U$ and $\V$ differ by one element only, say
$\U=\V\cup\{\ibf\}$ where $\ibf\in\nab\V$. It follows from (\ref{Xind}) and the
fact that $\nu$ solves the RDE (\ref{RDE}) that $X_\ibf$ has law $\nu$ and is
independent of $(X_\jbf)_{\jbf\in\nab\V\beh\{\ibf\}}$, and from this we see
that the projection of $\P_\U$ is indeed $\P_\V$.
\epro

It will be useful in what follows to have a somewhat stronger version of 
Lemma~\ref{L:RTP} that applies also to certain random subtrees $\U\sub\ovT$.
Let $\Ti$ denote the set of all finite subtrees $\U\sub\ovT$ such
that either $\wurz\in\U$ or $\U=\emptyset$. Let us define a \emph{stopping
  tree} to be a random variable $\U$ with values in $\Ti$ such that
\be\label{stoptree}
\{\U=\V\}\mbox{ is measurable w.r.t.\ }\sig\big((\omb_\ibf)_{\ibf\in\V}\big)
\quad\mbox{for all }\V\in\Ti.
\ee
In the special case that $\kappa\equiv 1$ and $\ovT=\N$, a stopping tree is
just a stopping time w.r.t.\ the filtration generated by
$\omb_\wurz,\omb_1,\omb_{11},\ldots$.

\bl[RTPs and stopping trees]
Let\label{L:stoptree} $(\omb_\ibf,X_\ibf)_{\ibf\in\ovT}$ be an RTP
corresponding to a map $\ga$ and a solution $\nu$ to the RDE (\ref{RDE}), and
let $\U\sub\ovT$ be a stopping tree. Then conditional on
$\U$, the random variables $(X_\ibf)_{\ibf\in\poa\U}$ are i.i.d.\ with common
law $\nu$ and independent of $(\omb_\ibf)_{\ibf\in\U}$.
\el
\bpro
For each fixed $\V\in\Ti$, by Lemma~\ref{L:RTP}, conditional on
$(\omb_\ibf)_{\ibf\in\V}$, the random variables $(X_\ibf)_{\ibf\in\poa\V}$ are
i.i.d.\ with common law $\nu$. By (\ref{stoptree}), it follows that
conditional on the event $\{\U=\V\}$ and $(\omb_\ibf)_{\ibf\in\V}$, the random
variables $(X_\ibf)_{\ibf\in\poa\V}$ are i.i.d.\ with common law $\nu$.
Since this holds for all $\V\in\Ti$, and since $\U\in\Ti$ a.s., the claim
follows.
\epro

\subsection{Continuous-time RTPs}\label{S:contRTP}

In this subsection, we prove Lemma~\ref{L:consis} and
Proposition~\ref{P:contRTP}. We work in our usual set-up as described above
Proposition~\ref{P:meanrep}. We start with a preparatory lemma that says that
if we condition on the \si-field $\Fi_t$ defined in (\ref{Fidef}), then the
subtrees of $\St$ rooted at $\ibf\in\nab\St_t$ are i.i.d.\ with the same
distribution as $\St$. To formulate this properly, we need some notation.

We call the object
\be\label{MBT}
\big(\St,(\omb_\ibf,\sig_\ibf)_{\ibf\in\St}\big).
\ee
a \emph{marked branching tree}. For each $\ibf\in\St$, let $\St^\ibf$ describe
the subtree of $\St$ that is rooted at $\ibf$, i.e.,
\be\label{Tibf}
\St^\ibf:=\{\jbf\in\ovT:\ibf\jbf\in\St\}.
\ee
We set $\om^\ibf_\jbf:=\omb_{\ibf\jbf}$ $(\ibf,\jbf\in\ovT)$, so that
$\om^\ibf_\jbf$ is the random element of $\Om$ that ``belongs'' to
$\jbf\in\St^\ibf$. Fix $t\geq 0$. For each $\ibf\in\nab\St_t$, let
$\sig^{\ibf,t}_\jbf$ describe the lifetime of an individual $\jbf\in\St^\ibf$
after time $t$, i.e.,
\be
\sig^{\ibf,t}_\wurz:=\sig_\ibf-(t-\tau^\ast_\ibf)
\quand
\sig^{\ibf,t}_\jbf:=\sig_{\ibf\jbf}\quad(\wurz\neq\jbf\in\St^\ibf),
\ee
where $t-\tau^\ast_\ibf$ is the age of the individual $\ibf$ at time $t$.

\bl[Memoryless property]
For\label{L:memless} each $t\geq 0$, conditional on the \si-field $\Fi_t$, the
marked branching trees
\be\label{memless}
\big(\St^\ibf,(\om^\ibf_\jbf,\sig^{\ibf,t}_\jbf)_{\jbf\in\St^\ibf}\big)_{\ibf\in\nab\St_t}
\ee
are i.i.d.\ with the same distribution as  the marked branching tree
in (\ref{MBT}).
\el
\bpro
Let $\Ti$ be as defined above (\ref{stoptree}). Then, for each
$\V\in\Ti$, the event $\{\St_t=\V\}$ is measurable w.r.t.\ the \si-field
generated by the random variables
\be\label{TtV}
(\omb_\ibf,\sig_\ibf)_{\ibf\in\V}\quand(1_{\{\sig_\ibf>t-\tau^\ast_\ibf\}})_{\ibf\in\nab\V}.
\ee
Note that here $\nab\V=\{\ibf j:\ibf\in\V,\ j\leq\kappa(\omb_\ibf)\}$ is
measurable w.r.t.\ the \si-field generated by $(\omb_\ibf)_{\ibf\in\V}$
while for each $\ibf\in\nab\V$, the random variable $\tau^\ast_\ibf$ is
measurable w.r.t.\ the \si-field generated by $(\sig_\ibf)_{\ibf\in\V}$.

Conditional on $\{\St_t=\V\}$ and the random variables in (\ref{TtV}),
the random variables $(\omb_\ibf)_{\ibf\in\ovT\beh\V}$ are still i.i.d.\ with
their original law and independent of $(\sig_\ibf)_{\ibf\in\ovT\beh\V}$.
The latter are also still independent of each other and the
$(\sig_\ibf)_{\ibf\in\ovT\beh(\V\cup\nab\V)}$ still have their original law,
but the laws of $(\sig_\ibf)_{\ibf\in\nab\V}$ are changed since conditioning on
$\{\St_t=\V\}$ entails conditioning on $\sig_\ibf>t-\tau^\ast_\ibf$ for each
$\ibf\in\nab\V$.

Since this holds for each $\V\in\St$, we see that if we condition on
$\Fi_t$ as in (\ref{Fidef}), then under the conditional law the random variables
$\omb_\ibf$ and $\sig_\ibf$ with $\ibf\in\ovT\beh\St_t$ are still independent,
and all of these random variables still have their original law, except the
$\sig_\ibf$ with $\ibf\in\nab\St_t$, whose laws are conditioned on the events
$\sig_\ibf>t-\tau^\ast_\ibf$. From this observation, using the memoryless
property of the exponential distribution, the claim of the lemma
follows.
\epro

For each $s\geq 0$ and $\ibf\in\nab\St_s$, within the marked branching tree
$\big(\St^\ibf,(\om^\ibf_\jbf,\sig^{\ibf,s}_\jbf)_{\jbf\in\St^\ibf}\big)$ rooted
at $\ibf$, we define the birth and death times $\tau^{\ibf,\ast}_\jbf$ and
$\tau^{\ibf,\dgg}_\jbf$ as in (\ref{taudef}), with $\sig_\jbf$ replaced by 
$\sig^{\ibf,s}_\jbf$, and we use this to define $\St^{\ibf,s}_t$ and
$\nab\St^{\ibf,s}_t$ $(t\geq 0)$ as in (\ref{Stdef}). Finally, we define
$G^{\ibf,s}_t=G_{\St^{\ibf,s}_t}$ as in (\ref{GUdef}) and (\ref{Gtdef}).\med

\bpro[of Lemma~\ref{L:consis}]
We fix a marked branching tree as in (\ref{MBT}) and times $0\leq s\leq t$.
Conditional on $\Fi_t$, we assign i.i.d.\ $(X_\ibf)_{\ibf\in\nab\St_t}$ with
common law $\mu_0$ to the leaves of $\St_t$ and define $(X_\ibf)_{\ibf\in\St_t}$
inductively as in (\ref{finrec}).

We observe that $\nab\St_t$ is given by the disjoint union
\be
\nab\St_t=\bigcup_{\ibf\in\nab\St_s}\{\ibf\jbf:\jbf\in\nab\St^{\ibf,s}_{t-s}\}.
\ee
Conditioning on $\Fi_t$ is the same as first conditioning on
\be\label{stree}
\big(\nab\St_s,(\omb_\jbf,\sig_\jbf)_{\jbf\in\St_s}\big),
\ee
and then conditioning on
\be
\big(\nab\St^{\ibf,s}_{t-s},
(\om^\ibf_\jbf,\sig^{\ibf,s}_\jbf)_{\jbf\in\St^{\ibf,s}_{t-s}}\big)_{\ibf\in\nab\St_s},
\ee
which by Lemma~\ref{memless} are conditionally independent given the random
variable in (\ref{stree}). Set
\be
X^\ibf_\jbf:=X_{\ibf\jbf}
\qquad(\jbf\in\St^{\ibf,s}_{t-s}\cup\nab\St^{\ibf,s}_{t-s},\ \ibf\in\nab\St_s).
\ee
Then
\be
X_\ibf=X^\ibf_\wurz=G^{\ibf,s}_{t-s}\big((X^\ibf_\jbf)_{\jbf\in\nab\St^{\ibf,s}_{t-s}}\big).
\ee
In view of this, by Theorem~\ref{T:meanrep}, conditional on the the random
variable in (\ref{stree}), i.e., conditional on $\Fi_s$, the random variables
$(X_\ibf)_{\ibf\in\nab\St_s}$ are i.i.d.\ with common law $\mu_{t-s}$, where 
$(\mu_s)_{s\geq 0}$ denotes the solution of the mean-field equation
(\ref{mean}) with initial state~$\mu_0$.
\epro

\bpro[of Proposition~\ref{P:contRTP}]
Since $(\sig_\ibf)_{\ibf\in\ovT}$ and $(\omb_\ibf,X_\ibf)_{\ibf\in\ovT}$ are
independent, the conditional law of $(\omb_\ibf,X_\ibf)_{\ibf\in\ovT}$ given
$(\sig_\ibf)_{\ibf\in\ovT}$ is the same as the unconditional law.
We claim that under the conditional law given $(\sig_\ibf)_{\ibf\in\ovT}$,
the random finite subtree $\St_t$ is a stopping tree in the sense of
(\ref{stoptree}). Indeed, $\St_t=\V$ if and only if for each $\ibf\in\V$ and
$j\in\N_+$ (resp.\ $j\in[d]$, depending on how $\ovT$ is chosen), one has
$\ibf j\in\V$ if and only if
\be
{\rm(i)}\ j\leq\kappa(\ibf)
\quand
{\rm(ii)}\ \tau^\dgg_{\ibf j}\leq t.
\ee
Here the event in (i) is clearly measurable
w.r.t.\ $\sig((\omb_\ibf)_{\ibf\in\V})$ while under the conditional law given
$(\sig_\ibf)_{\ibf\in\ovT}$, (ii) is just a deterministic condition.
We can therefore apply Lemma~\ref{L:stoptree} to conclude that 
conditional on $(\sig_\ibf)_{\ibf\in\ovT}$, $\St_t$, and
$(\omb_\ibf)_{\ibf\in\St_t}$, the random variables $(X_\ibf)_{\ibf\in\poa\St_t}$
are i.i.d.\ with common law $\nu$.

We observe that
\be
\nab\St_t=\big\{\ibf j:\ibf\in\St_t,\ \ibf j\not\in\St_t,
\ j\leq\kappa(\omb_\ibf)\big\}
\ee
is a function of $\St_t$, and $(\omb_\ibf)_{\ibf\in\St_t}$. Therefore, if we
condition on $\Fi_t=\sig(\nab\St_t,(\omb_\ibf,\sig_\ibf)_{\ibf\in\St_t})$,
the random variables $(X_\ibf)_{\ibf\in\nab\St_t}$ are i.i.d.\ with common law
$\nu$. This proves (\ref{contRTP})~(i). Condition (\ref{contRTP})~(ii) is also
clearly fulfilled by the definition of an RTP.
\epro

\subsection{Endogeny, bivariate uniqueness, and the higher-level equation}
\label{S:endlev}

In this subsection, we prove Theorem~\ref{T:bivar2} and
Propositions~\ref{P:levmom} and \ref{P:minconv}.

Recall that an RTP $(\omb_\ibf,X_\ibf)_{\ibf\in\ovT}$ is endogenous if
$X_\wurz$ is measurable with respect to the \si-field generated by the random
variables $(\omb_\ibf)_{\ibf\in\ovT}$. In general, if $X$ is a random variable
taking values in a Polish space and $\Fi$ is a sub-\si-field, then it is not
hard to see that $X$ is a.s.\ equal to a $\Fi$-measurable function if and only
if the conditional law $\P[X\in\,\cdot\,|\Fi]$ is a.s.\ a delta-measure. In
view of this, the following lemma implies that an RTP is endogenous if and only
if $X_\wurz$ is a.s.\ measurable w.r.t.\ the \si-field generated by the random
variables $\St$ and $(\omb_\ibf)_{\ibf\in\St}$.

\bl[Relevant randomness]
Let\label{L:endT} $(\omb_\ibf,X_\ibf)_{\ibf\in\ovT}$ be an RTP corresponding to
a solution $\nu$ of the RDE (\ref{RDE}). Let $\ov\Fi$ be
the \si-field generated by the random variables $(\omb_\ibf)_{\ibf\in\ovT}$
and let $\Fi$ be the \si-field generated by the random variables $\St$ and
$(\omb_\ibf)_{\ibf\in\St}$. Then
\be\label{relevant}
\P[X_\wurz\in\,\cdot\,|\ov\Fi]=\P[X_\wurz\in\,\cdot\,|\Fi]\quad{\rm a.s.}
\ee
\el
\bpro
Since $\ov\Fi$ is generated by $\Fi$ and the random variables
$(\omb_\ibf)_{\ibf\in\ovT\beh\St}$, formula (\ref{relevant}) says that
conditional on on $\Fi$, the random variables
$(\omb_\ibf)_{\ibf\in\ovT\beh\St}$ are independent of $X_\wurz$.
Let $\U^{(n)}$ be deterministic finite rooted subtrees of $\ovT$ that
 increase to $\ovT$. Let $\ov\Fi^{(n)}$ be the \si-field generated by
$(\omb_\ibf)_{\ibf\in\U^{(n)}}$ and let $\Fi^{(n)}$ be the \si-field
generated by $\St\cap\U^{(n)}$ and $(\omb_\ibf)_{\ibf\in\St\cap\U^{(n)}}$.
Conditional on $\Fi^{(n)}$, the state at the root $X_\wurz$ is a deterministic
function of $(X_\ibf)_{\ibf\in\nab(\St\cap\U^{(n)})}$. Therefore, by point~(ii)
in the definition of an RTP in Lemma~\ref{L:RTP}, $X_\wurz$ is conditionally
independent of $(\omb_\ibf)_{\ibf\in(\ovT\beh\St)\cap\U^{(n)}}$ given $\Fi^{(n)}$, 
or equivalently,
\be
\P[X_\wurz\in A\,|\ov\Fi^{(n)}]
=\P[X_\wurz\in A\,|\Fi^{(n)}]\quad{\rm a.s.}
\ee
for each measurable $A\sub S$. Letting $n\to\infty$, using martingale
convergence, we arrive at (\ref{relevant}).
\epro

The following lemma prepares for the proof of Theorem~\ref{T:bivar2}.

\bl[Successful coupling]
Let\label{L:auxl} $(\omb_\ibf,X_\ibf)_{\ibf\in \T}$ be an endogenous RTP
corresponding to a solution $\nu$ of the RDE \eqref{RDE} and let 
$(\sig_\ibf)_{\ibf\in\ovT}$ be an independent i.i.d.\ collection of
exponential random variables with mean $|\rb|^{-1}$.
Furthermore, let $(Y_\ibf)_{\ibf\in \T}$ be an i.i.d.\ collection of
$S$-valued random variables with common law $\nu$, independent of
$(\omb_\ibf,X_\ibf,\sig_\ibf)_{\ibf\in \T}$. For each $t>0$, define random
variables $(X^t_\ibf)_{\ibf\in\St_t \cup \nab \St_t}$ by
\be\ba{rl}
{\rm(i)}&X^t_\ibf:=Y_\ibf\quad (\ibf\in\nab \St_t)\\[5pt]
{\rm(ii)}&\dis X^t_\ibf=\ga[\omb_\ibf]\big(X^t_{\ibf1},\ldots,X^t_{\ibf\kappa(\omb_\ibf)}\big) \qquad(\ibf\in\St_t).
\ec
Then
\be
\label{convprob}
X_\wurz^t \asto{t} X_\wurz \textrm{ in probability}.
\ee
\el
\bpro
The following argument is a continuous-time version of the proofs of
\cite[Thm~11~(c)]{AB05} and \cite[Lemma~6]{MSS18}. Let $\Fi_t$ be the
filtration defined in \eqref{Fidef}. We add a final element
$\Fi_\infty:=\sig(\bigcup_{t\geq 0}\Fi_t)$ to the filtration, which is the
\si-algebra generated by the random tree $\St$ and the random variables
$(\omb_\ibf,\sig_\ibf)_{\ibf\in\St}$. Let $f,g:S \to \R$ be bounded and
measurable functions. Since $X_\wurz$ and $X^t_\wurz$ are conditionally
independent and identically distributed given $\Fi_t$, we have
\be\ba{l}\label{copro}
\E[f(X_\wurz)g(X^t_\wurz)]
=\E\big[\E[f(X_\wurz)|\Fi_t]\E[g(X^t_\wurz)|\Fi_t]\big]
=\E\big[\E[f(X_\wurz)|\Fi_t]\E[g(X_\wurz)|\Fi_t]\big]\\[5pt]
\dis\qquad\asto{t}
\E\big[\E[f(X_\wurz)|\Fi_\infty]\E[g(X_\wurz)|\Fi_\infty]\big]
=\E[f(X_\wurz)g(X_\wurz)],
\ec
where we used the martingale convergence and in the last equality also
endogeny and Lemma \ref{L:endT}. Since \eqref{copro} holds in particular for
any bounded continuous $f$ and $g$, we conclude that the law of
$(X_\wurz,X_\wurz^t)$ converges weakly to the law of $(X_\wurz,X_\wurz)$, which
implies \eqref{convprob}.
\epro

\bpro[of Theorem~\ref{T:bivar2}]
If (ii) holds, then $\ov\nu^{(2)}$ is the only fixed point in $\Pc(S^2)_\nu$ of the bivariate
mean-field equation. Since a measure is a fixed point of the bivariate mean-field equation if and
only if it is a fixed point of the map $T^{(2)}$, by Theorem~\ref{T:bivar}, it
follows that the RTP corresponding to $\nu$ is endogenous.

Assume, conversely, that the RTP corresponding to $\nu$ is endogenous. Let
$(Y_\ibf^1,\ldots,Y_\ibf^n)_{\ibf\in \T}$ be a collection of
i.i.d. $S^n$-valued random variables with common law $\mu^{(n)}_0$,
independent of the RTP $(\omb_\ibf,X_\ibf)_{\ibf\in\ovT}$ and the exponential
lifetimes $(\sig_\ibf)_{\ibf\in\ovT}$. For each $1\leq m\leq n$ and $t>0$,
define random variables $(X^{m,t}_\ibf)_{\ibf\in\St_t \cup \nab \St_t}$ by
\be\ba{rl}
{\rm(i)}&X^{m,t}_\ibf:=Y^m_\ibf\quad (\ibf\in\nab \St_t)\\[5pt]
{\rm(ii)}&\dis X^{m,t}_\ibf=\ga[\omb_\ibf]\big(X^{m,t}_{\ibf1},\ldots,X^{m,t}_{\ibf\kappa(\omb_\ibf)}\big) \qquad(\ibf\in\St_t).
\ec
Then, by Theorem~\ref{T:meanrep} applied to the $n$-variate map $\ga^{(n)}$,
we see that $(X_\wurz^{1,t},\ldots,X_\wurz^{n,t})$ has law $\mu^{(n)}_t$.
By endogeny we get from Lemma \ref{L:auxl} that 
\be\label{Xnconv}
(X_\wurz^{1,t},\ldots,X_\wurz^{n,t}) \asto{t} (X_\wurz,\ldots,X_\wurz) \textrm{ in probability}.
\ee
This completes the proof since the right-hand side of \eqref{Xnconv} has law $\ov\nu^{(n)}$ as defined in \eqref{nmu}.
\epro

\bpro[of Proposition~\ref{P:levmom}]
The fact that $(\rho_t)_{t\geq 0}$ solves the higher-level mean-field equation
(\ref{levelmean}) means that
\be\label{levelmeanweak}
\dif{t}\li\rho_t,\phi\re=\int_\Om\rb(\di\om)
\big\{\li T_{\ch\ga[\om]}(\rho_t),\phi\re-\li\rho_t,\phi\re\big\}\qquad(t\geq 0),
\ee
for any bounded measurable $\phi:\Pc(S)\to\R$. In particular, we can apply
this to functions of the form
\be
\phi(\mu)=\int_S\mu(\di x^1)\cdots\int_S\mu(\di x^n)\,f(x^1,\ldots,x^n),
\ee
where $f:S^n\to\R$ is bounded and measurable. Then
\be
\li\rho_t,\phi\re=\int_{S^n}\!\rho^{(n)}_t(\di x)\,f(x)
\quand
\li T_{\ch\ga[\om]}(\rho_t),\phi\re
=\int_{S^n}\!\big(T_{\ch\ga[\om]}(\rho_t)\big)^{(n)}(\di x)\,f(x),
\ee
where $(T_{\ch\ga[\om]}(\rho_t))^{(n)}$ denotes the $n$-th moment measure
of $T_{\ch\ga[\om]}(\rho_t)$. By \cite[Lemma~2]{MSS18},
\be
\big(T_{\ch\ga[\om]}(\rho_t)\big)^{(n)}=T_{\ga^{(n)}[\om]}(\rho^{(n)}_t).
\ee
Inserting this into (\ref{levelmeanweak}), we see that $(\rho^{(n)}_t)_{t\geq
  0}$ solves the $n$-variate mean-field equation.
\epro

The following lemma prepares for the proof of Proposition~\ref{P:minconv}.

\bl[Conditional law of the root]
Let\label{L:condroot} $(\omb_\ibf,X_\ibf)_{\ibf\in\ovT}$ be an RTP
corresponding to a solution $\nu$ of the RDE (\ref{RDE}), let
$(\sig_\ibf)_{\ibf\in\ovT}$ be an independent i.i.d.\ collection of
exponentially distributed random variables with mean $|\rb|^{-1}$, and let
$(\Fi_t)_{t\geq 0}$ be the filtration defined in (\ref{Fidef}). Then the measures
\be\label{condroot}
\rho_t:=\P\big[\P[X_\wurz\in\,\cdot\,|\Fi_t]\in\,\cdot\,\big]
\qquad(t\geq 0)
\ee
solve the higher-level mean-field equation (\ref{levelmean}) with initial state
$\rho_0=\de_\nu$.
\el
\bpro
Conditional on $\Fi_t$, the map $G_t:S^{\nab\St_t}\to S$ is a deterministic map,
and $(X_\ibf)_{\ibf\in\nab\St_t}$ are i.i.d.\ with common law $\nu$. Therefore,
applying \cite[Lemma~8]{MSS18} to the case that the \si-fields $\Hi_k$ there
are all trivial and the probability measure $\P$ there is replaced by the
conditional law given $\Fi_t$, we see that
\bc
\dis T_{\ch G_t}(\de_\nu)
&=&\dis\P\big[\P[G_t\big((X_\ibf)_{\ibf\in\nab\St_t}\big)\in\,\cdot\,\big|\,\Fi_t]
\in\,\cdot\,\big|\,\Fi_t\big]\\[5pt]
&=&\dis\P\big[\P[X_\wurz\in\,\cdot\,\big|\,\Fi_t]\in\,\cdot\,\big|\,\Fi_t\big]
\ec
Now by Theorem~\ref{T:meanrep},
\be
\rho_t:=\E[T_{\ch G_t}(\de_\nu)]
=\P\big[\P[X_\wurz\in\,\cdot\,\big|\,\Fi_t]\in\,\cdot\,\big]
\qquad(t\geq 0)
\ee
solves the higher-level mean-field equation (\ref{levelmean}) with initial state
$\rho_0=\de_\nu$.
\epro

\bpro[of Proposition~\ref{P:minconv}]
Let $(\rho^i_t)_{t\geq 0}$ $(i=1,2)$ be solutions to the
higher-level mean-field equation (\ref{levelmean}) such that
$\rho^1_0\leq_{\rm cv}\rho^2_0$. Define $\rho^{i}_{t,(n)}$ as in
(\ref{mun}), with $T$ replaced by the higher-level map $\ch T$
from (\ref{RDElevel}). It has been shown in \cite[Prop~3]{MSS18} that $\ch T$
is monotone w.r.t.\ the convex order, so by induction we obtain from
(\ref{mun}) that $\rho^{1}_{t,(n)}\leq_{\rm cv}\rho^{2}_{t,(n)}$ for all $n\geq 1$
and $t\geq 0$. Letting $n\to\infty$, using (\ref{munconv}), we see that
$\rho^1_t\leq_{\rm cv}\rho^2_t$ for all $t\geq 0$.

Let $\nu$ be a solution of the RDE (\ref{RDE}). It has been shown in
\cite[Prop.~3]{MSS18} that $\ov\nu$ solves the higher-level RDE
(\ref{RDElevel}) and there exists a (necessarily unique) solution $\un\nu$ of
(\ref{RDElevel}) such that (\ref{sandwich}) holds. It has moreover been shown
in \cite[Prop.~4]{MSS18} that $\un\nu$ is given by (\ref{unnu}). In view of
this, to complete the proof, it suffices to show that the solution
$(\rho_t)_{t\geq 0}$ to the higher-level mean-field equation (\ref{levelmean})
with initial state $\rho_0=\de_\nu$ converges to the measure in (\ref{unnu}).

We apply Lemma~\ref{L:condroot}. As in the proof of Lemma~\ref{L:auxl}, we add
a final element $\Fi_\infty:=\sig(\bigcup_{t\geq 0}\Fi_t)$ to the filtration,
which is the \si-algebra generated by the random tree $\St$ and the random
variables $(\omb_\ibf,\sig_\ibf)_{\ibf\in\St}$. Then, by martingale convergence,
\be
\P[X_\wurz\in\,\cdot\,|\Fi_t]\Asto{t}\P[X_\wurz\in\,\cdot\,|\Fi_\infty]
\quad{\rm a.s.},
\ee
and hence the measures $\rho_t$ in (\ref{condroot}) satisfy
\be\label{rhoun}
\rho_t\Asto{t}\P\big[\P[X_\wurz\in\,\cdot\,|\Fi_\infty]\in\,\cdot\,\big],
\ee
where $\Rightarrow$ denotes weak convergence of probability measures on
$\Pc(S)$, which is in turn equipped with the topology of weak convergence of
probability measures on $S$. Since the exponentially distributed random
variables $(\sig_\ibf)_{\ibf\in\ovT}$ are independent of the RTP
$(\omb_\ibf,X_\ibf)_{\ibf\in\ovT}$, we have
\be
\P[X_\wurz\in\,\cdot\,|\Fi_\infty]=\P[X_\wurz\in\,\cdot\,|\Fi]
=\P[X_\wurz\in\,\cdot\,|(\omb_\ibf)_{\ibf\in\ovT}],
\ee
where as in Lemma~\ref{L:endT} $\Fi$ denotes the \si-field generated by the
random variables $\St$ and $(\omb_\ibf)_{\ibf\in\St}$ and the last equality
follows from that lemma. Inserting this into (\ref{rhoun}) we see that
$\rho_t$ converges weakly to $\un\nu$ as defined in (\ref{unnu}).
\epro

\section{Further results}\label{S:fur}

In this section, we prove some additional results about RTPs.  In
Subsection~\ref{S:mono}, we prove Proposition~\ref{P:lowup} about the upper
and lower solutions of a monotonous RDE. In Subsection~\ref{S:uni} we prove
Lemmas~\ref{L:rootdet}, \ref{L:detree}, and \ref{L:monuni} which
give conditions for uniqueness of solutions to an RDE. Subsection~\ref{S:dual}
is devoted to the proof of Lemma~\ref{L:good}.

\subsection{Monotonicity}\label{S:mono}

In this subsection, we prove Proposition~\ref{P:lowup}.
We start with a number of simple lemmas.

\bl[A continuous monotone function]
Let\label{L:fxy} $S$ be a compact metrizable space that is equipped with a
closed partial order in the sense of (\ref{compat}), and let $d$ be a
metric that generates the topology. Then
\be\label{fxy}
f(x,y):=\inf\big\{d(x',y'):x'\leq x,\ y'\geq y\big\}
\qquad(x,y\in S)
\ee
defines a continuous function $f:S^2\to\half$ such that $f(x,y)=0$ if and only
if $x\geq y$ and moreover $f(x,y)$ is decreasing in $x$ and increasing in $y$.
\el

\bpro
Since for any $(x,y),(x',y')\in S^2$,
\be
\big|d(x,y)-d(x',y')\big|\leq d(x,x')+d(y,y'),
\ee
the function $d:S^2\to\half$ is continuous. Assume that $(x_n,y_n)\in S^2$
converge to a limit $(x,y)$. Since the infimum of a family of
continuous functions is upper semi-continuous, we have
\be
\limsup_{n\to\infty}f(x_n,y_n)\leq f(x,y).
\ee
To prove that $f$ is actually continuous, assume the converse. Then there
exists a sequence such that
\be\label{notcon}
\liminf_{n\to\infty}f(x_n,y_n)\leq f(x,y)-\eps.
\ee
for some $\eps>0$. By the definition of $f$, there exist $x'_n\leq x_n$ and
$y'_n\geq y_n$ such that $d(x'_n,y'_n)\leq f(x_n,y_n)+\eps/2$. Since $S$ is
compact, we can select a subsequence such that (\ref{notcon}) still holds and
the $(x'_n,y'_n)$ converge to a limit $(x',y')$. Since the partial order is
closed in the sense of (\ref{compat}), we have $x'\leq
x$ and $y'\geq y$, so
\be
f(x,y)\leq d(x',y')
=\lim_{n\to\infty}d(x'_n,y'_n)\leq\eps/2+\liminf_{n\to\infty}f(x_n,y_n),
\ee
which contradicts (\ref{notcon}). We conclude that $f:S^2\to\half$ is
continuous.

If $x\geq y$, then setting $(x',y')=(x,x)$ shows that $f(x,y)=0$. Conversely,
if $f(x,y)=0$ then there exist $x_n\leq x$ and $y_n\geq y$ such that
$d(x_n,y_n)\to 0$. Using the compactness of $S$, by going to a subsequence, we
can assume that the $(x_n,y_n)$ converge to a limit $(z,z)$. Since the partial
order is closed in the sense of (\ref{compat}), $y\leq
z\leq x$ and hence $x\geq y$.

If $x\leq x_\ast$ and $y\geq y_\ast$, then
\be
f(x_\ast,y_\ast):=\inf\big\{d(x',y'):x'\leq x_\ast,\ y'\geq y_\ast\big\}
\leq\inf\big\{d(x',y'):x'\leq x,\ y'\geq y\big\}=f(x,y)
\ee
since the second infimum is taken over a smaller set, showing that
$f(x,y)$ is decreasing in $x$ and increasing in $y$.
\epro

\bl[Comparison principle]
Let\label{L:compar} $S$ be a compact metrizable space that is equipped with a
partial order that is closed in the sense of
(\ref{compat}). Let $X,Y$ be $S$-valued random variables such that $X\leq Y$
a.s.\ and $\P[X\in\,\cdot\,]\geq\P[Y\in\,\cdot\,]$. Then $X=Y$ a.s.
\el

\bpro[of Lemma~\ref{L:compar}]
Set $f_z(x):=f(z,x)$ with $f$ as in Lemma~\ref{L:fxy}. Then, for each $z\in
S$, $f_z:S\to\half$ is continuous and monotone increasing, and $f_z(x)=0$ if and only if
$x\leq z$. Let
\be
S^2_<:=\big\{(x,y)\in S^2:x\leq y,\ x\neq y\big\}.
\ee
We will prove the lemma by showing that if $X,Y$ are $S$-valued random
variables such that $\P[(X,Y)\in S^2_<]>0$, then $\E[f_z(X)]<\E[f_z(Y)]$ for
some $z\in S$ contradicting $\P[X\in\,\cdot\,]\geq\P[Y\in\,\cdot\,]$.
For each $z\in S$ and $\de>0$, we define an open set
$O_{z,\de}\sub S^2$ by
\be\label{Oz}
O_{z,\de}:=\{(x,y)\in S^2:f_z(y)-f_z(x)>\de\}.
\ee
Since for each $(x,y)\in S^2_<$, one has $f_x(x)=0$ but $f_x(y)>0$, we see
that
\be
\bigcup\{O_{z,\de}:z\in S,\ \de>0\}\supset S^2_<.
\ee
We now use the inner regularity of measures on Polish spaces w.r.t.\ compacta, 
which follows from the regularity and tightness of any probability measure on
a Polish space \cite[Thm.~1.2 and 3.2]{Par05}. Thus, we can find a compact set
$K\sub S^2_<$ such that $\P[(X,Y)\in K]>0$. Since $K$ is compact, it is
covered by finitely many sets of the form (\ref{Oz}), so there must exists a
$z\in S$ and $\de>0$ such that $\P[(X,Y)\in O_{z,\de}]>0$. Since $f_z$ is
monotone increasing and $X\leq Y$ it follows that $\E[f_z(X)]<\E[f_z(Y)]$.
\epro

\bl[Compatibility of the stochastic order]
Assume\label{L:compat} that $S$ is equipped with a partial order that is
closed in the sense of (\ref{compat}). Then the stochastic order on
$\Pc(S)$ is closed with respect to the topology of weak convergence.
\el
\bpro
We need to show that if $\mu^1_n\leq\mu^2_n$ for all $n\in\N$ and the
$\mu^i_n\in\Pc(S)$ converge weakly as $n\to\infty$ to a limit $\mu^i_\infty$
$(i=1,2)$, then $\mu^1_\infty\leq\mu^2_\infty$. Since $\mu^1_n\leq\mu^2_n$,
for each $n$, we can couple $X^i_n$ with laws $\mu^i_n$ $(i=1,2)$ such that
$X^1_n\leq X^2_n$. Since $\mu^1_n$ and $\mu^2_n$ converge as $n\to\infty$,
the joint laws of $(X^1_n,X^2_n)$ are tight, so by going to a subsequence we
may assume that they converge. Then, by Skorohod's representation
theorem, we can couple the random variables $(X^1_n,X^2_n)$ for different $n$
in such a way that they converge a.s.\ to a limit $(X^1_\infty,X^2_\infty)$.
Since the partial order on $S$ is closed, we have 
$X^1_\infty\leq X^2_\infty$ a.s., proving that $\mu^1_\infty\leq\mu^2_\infty$.
\epro

\bl[Monotonicity of $T$]
Assume\label{L:gamon} that $S$ is equipped with a partial order that is
closed and that $\ga[\om]$ is monotone for all
$\om\in\Om$. Then the operator $T$ in (\ref{Tdef}) is monotone w.r.t.\ the
stochastic order.
\el
\bpro
If $\mu_1\leq\mu_2$, then we can couple random variables $X^1$ and $X^2$ with
laws $\mu_1,\mu_2$ such that $X^1\leq X^2$. Let $(X^1_i,X^2_i)_{i\geq 1}$ be
i.i.d.\ copies of $(X^1,X^2)$. Then
\be
\ga[\om](X^1_1,\ldots,X^1_{\kappa(\om)})
\leq\ga[\om](X^2_1,\ldots,X^2_{\kappa(\om)})
\ee
for all $\om\in\Om$ and hence $T(\mu_1)\leq T(\mu_2)$ by (\ref{Tdef}).
\epro

In practice, Lemma~\ref{L:gamon} is the usual way to prove monotonocity of a
map of the form (\ref{Tdef}). Nevertheless, it is known that there are  
maps of the form (\ref{Tdef}), in particular, probability kernels, that are
monotone yet cannot be represented in terms of monotone maps
\cite[Example~1.1]{FM01}.

\bl[Monotonicity in the initial state]
Assume\label{L:mumon} that $S$ is equipped with a partial order that is
closed and that the operator $T$ in (\ref{Tdef}) is
monotone w.r.t.\ the stochastic order. Then solutions $(\mu^i_t)_{t\geq 0}$
$(i=1,2)$ of the mean-field equation (\ref{mean}) started in initial
states $\mu^1_0\leq\mu^2_0$ satisfy $\mu^1_t\leq\mu^2_t$ $(t\geq 0)$.
\el
\bpro
Inductively define $\mu^i_{t,(n)}$ as in (\ref{mun}) with $\mu_0$ replaced by
$\mu^i_0$ $(i=1,2)$. Then $\mu^1_{t,(n)}\leq\mu^2_{t,(n)}$ for all $n\geq 1$
and $t\geq 0$. Letting $n\to\infty$, we see as in the proof of
Proposition~\ref{P:meanrep} that $\mu^i_{t,(n)}\Rightarrow\mu^i_t$ as
$n\to\infty$. By Lemma~\ref{L:compat}, we conclude that $\mu^1_t\leq\mu^2_t$
$(t\geq 0)$.
\epro

In the next two lemmas we need to assume compactness of $S$.

\bl[Increasing limits]
Assume\label{L:inclim} that $S$ is a compact metrizable space equipped with a
partial order that is closed. Then every increasing
sequence in $S$ converges to a limit.
\el
\bpro
Let $(x_n)_{n\in\N}$ be a sequence in $S$ such that $x_n\leq x_{n+1}$ for all
$n\in\N$. By compactness, it suffices to prove that all subsequential limits
are the same. Let $(x_m)_{m\in M}$ and $(x_k)_{k\in K}$ be subsequences that
converge to limits $x$ and $x'$, respectively. For all $k\in K$, let
$k_-:=\sup\{m\in M:m\leq k\}$. Then $x_{k_-}\leq x_k$ for all $k\in K$ and
letting $k\to\infty$, using the compatibility condition
(\ref{compat}), we see that $x\leq x'$. The same argument gives $x'\leq x$ and
hence $x=x'$.
\epro

\bl[Increasing limits in the stochastic order]
Assume\label{L:stochinc} that $S$ is a compact metrizable space equipped with a
partial order that is closed. Then every sequence in
$\Pc(S)$ that is increasing in the stochastic order converges weakly to a
limit.
\el
\bpro
Let $(\mu_n)_{n\geq 0}$ be increasing in the stochastic order. Then, for each
$n\geq 1$, we can couple random variables $X_{n-1}$ and $X_n$ with laws
$\mu_{n-1}$ and $\mu_n$ such that $\P[X_{n-1}\leq X_n]=1$. Let $K_n(x,\di
y):=\P[X_n\in\,\di y\,|\,X_{n-1}=x]$ and let $(Y_n)_{n\geq 0}$ be a
time-inhomogeneous Markov chain with initial law $\mu_0$ and transition
kernels $\P[Y_n\in\,\di y\,|\,Y_{n-1}=x]=K_n(x,\di y)$. Then $Y_{n-1}\leq Y_n$
a.s.\ for all $n\geq 1$ and hence the $Y_n$ a.s.\ increase to a limit $Y_\infty$
by Lemma~\ref{L:inclim}. It follows that the $\mu_n$ converge weakly to the
law of $Y_\infty$.
\epro

\noindent
We now turn to the proof of Proposition~\ref{P:lowup}. 

\bl[Lower and upper solutions]
All\label{L:lowup} conclusions of Proposition~\ref{P:lowup} except for the
statement about endogeny hold when the assumption that $\ga[\om]$ is monotone
for all $\om\in\Om$ is replaced by the weaker condition that $T$ is monotone.
\el

\bpro
The proof is similar to the proof of \cite[Lemma~15]{AB05}, which in turn
is based on well-known principles \cite[Thm~III.2.3]{Lig85}. By symmetry, it
suffices to prove the statement for $\nu_{\rm low}$.

Since for each $s\geq 0$, $(\mu^{\rm low}_{s+t})_{t\geq 0}$ solves
(\ref{mean}) with initial state $\mu^{\rm low}_s\geq\de_0$, we conclude
from Lemmas~\ref{L:gamon} and \ref{L:mumon} that
$\mu^{\rm low}_{s+t}\geq\mu^{\rm low}_t$ for each $s\geq 0$ and hence
$t\mapsto\mu^{\rm low}_t$ is increasing w.r.t.\ the stochastic order. By
Lemma~\ref{L:stochinc}, it follows that $\mu^{\rm low}_t\Rightarrow\nu_{\rm
  low}$ for some probability measure $\nu_{\rm low}$ on $S$. Since
$\mu^{\rm low}_{s+t}\Rightarrow\nu_{\rm low}$ for all $s\geq 0$, using
Lemma~\ref{L:mucon} and the continuity of $T$, we see that $\nu_{\rm low}$ is
a fixed point of the mean-field equation (\ref{mean}) and hence solves the RDE
(\ref{RDE}).

If $\nu$ is any solution of the RDE (\ref{RDE}), then $\mu^{\rm low}_t\leq\nu$
for all $t\geq 0$ by Lemma~\ref{L:mumon} and the fact that $\nu$ is a fixed
point of (\ref{mean}). Letting $t\to\infty$, using Lemma~\ref{L:compat},
we see that $\nu_{\rm low}\leq\nu$.
\epro

\bl[Random maps applied to extremal elements]
Under\label{L:extrmap} the assumptions of Proposition~\ref{P:lowup}, if
$(\om_\ibf)_{\ibf\in\T}$ are i.i.d.\ with common law $|\rb|^{-1}\rb$, then
there exist random variables $X^{\rm upp}_\wurz$ and $X^{\rm low}_\wurz$ with
laws $\nu_{\rm upp}$ and $\nu_{\rm low}$ that are given by the decreasing,
resp.\ increasing limits
\be\label{Xulo}
X^{\rm upp}_\wurz=\lim_{\U^{(n)}\up\St}G_{\U^{(n)}}(1,\ldots,1)
\quand
X^{\rm low}_\wurz:=\lim_{\U^{(n)}\up\St}G_{\U^{(n)}}(0,\ldots,0),
\ee
where the limit does not depend on the choice of the sequence $\U^{(n)}\in\Ti$
such that $\U^{(n)}\up\St$. Here $\Ti$ denotes the set of all finite
subtrees $\U\sub\ovT$ such that either $\wurz\in\U$ or $\U=\emptyset$,
and for each $\U\in\Ti$, the random map $G_\U:S^{\nab\U}\to S$ is defined in
(\ref{GUdef}).
\el

\bpro
By symmetry, it suffices to prove the statement for $X^{\rm low}_\wurz$.
Since $\ga[\om]$ is monotone for each $\om\in\Om$, the map $G_\U$ is monotone
for each $\U\in\Ti$. Define
\be\label{XUp}
X^\U_\wurz:=G_\U(0,\ldots,0)\qquad(\U\in\Ti).
\ee
Then $X^\U_\wurz\leq X^\V_\wurz$ for all $\U\sub\V$ and hence
if $\U^{(n)}\in\Ti$ increase to $\St$, then the $X^{\U^{(n)}}_\wurz$ increase
to a limit $X^{\rm low}_\wurz$ that does not depend on the choice of the
sequence $\U^{(n)}$.

Let $(\sig_\ibf)_{\ibf\in\ovT}$ be an independent i.i.d.\ collection of
exponential random variables with mean $|\rb|^{-1}$. Define $\St_t\in\Ti$ as
in (\ref{Stdef}). Then by Theorem~\ref{T:meanrep}, $X^{\St_t}_\wurz$ has law
$\mu^{\rm low}_t$ while by what we have already proved $X^{\St_t}_\wurz$
increases to $X^{\rm low}_\wurz$. Since $\mu^{\rm
  low}_t\Rightarrow\nu_{\rm low}$, it follows that $X^{\rm low}_\wurz$
has law $\nu_{\rm low}$.
\epro

\bpro[of Proposition~\ref{P:lowup}]
In view of Lemma~\ref{L:lowup}, it only remains to prove the
statement about endogeny. Let $(\omb_\ibf,X_\ibf)_{\ibf\in\ovT}$ be an RTP
corresponding to $\ga$ and some solution $\nu$ to the RDE (\ref{RDE}). Then
\be
X_\wurz:=G_\U\big((X_\ibf)_{\ibf\in\nab\U}\big)\geq X^\U_\wurz
\qquad(\U\in\Ti),
\ee
with $X^\U_\wurz$ as in (\ref{XUp}).
So letting $\U\up\ovT$, using the fact that the partial order is closed,
we obtain that $X_\wurz\geq X^{\rm low}_\wurz$. In
particular, if $\nu=\nu_{\rm low}$, then since $X^{\rm low}_\wurz$ also has law
$\nu_{\rm low}$, Lemma~\ref{L:compar} tells us that 
$X_\wurz=X^{\rm low}_\wurz$ a.s. Since the latter is measurable w.r.t.\ the
\si-field generated by the $(\omb_\ibf)_{\ibf\in\ovT}$, this proves the
endogeny of the RTP corresponding to $\ga$ and $\nu_{\rm low}$.
\epro

\subsection{Conditions for uniqueness}\label{S:uni}

In this subsection, we prove Lemmas~\ref{L:rootdet}, \ref{L:detree},
and \ref{L:monuni}.\med

\bpro[of Lemma~\ref{L:rootdet}]
If $G_t$ is constant  then $\St_t$ is a root determining subtree,
proving the implication (i)$\volgt$(ii). Conversely, if there a.s.\ exists a
root determining subtree $\U$, then, since $\St_t\up\St$, there a.s.\ exists a
(random) $t<\infty$ such that $\St_t \supset\U$ and hence $G_s$ is
constant for all $s\geq t$.
The implication (iii)$\volgt$(ii) is trivial. Conversely, if $\St$ contains a
root determining subtree $\U$, then by the finiteness of the latter we can
keep removing elements from $\U$ as long as this is still possible while
retaining the property that $\U$ is root determining.
\epro

\bpro[of Lemma~\ref{L:detree}]
(i)$\volgt$(ii): This is clear, since a finite uniquely determined subtree is
root determining.

(ii)$\volgt$(iii): For each $\ibf\in\St$, let $\St^\ibf$, defined in
(\ref{Tibf}), denote the subtree of $\St$ that is rooted at $\ibf$. Since
$\St^\ibf$ is equally distributed with $\St$, by (ii), for each $\ibf\in\St$,
there a.s.\ exists a root determining subtree $\U^\ibf\sub\St^\ibf$. Since
$x\in\Xi_\St$ implies
$x_\ibf=G_{\U^\ibf}\big((x_{\ibf\jbf})_{\jbf\in\nab\U^\ibf}\big)$ and
$G_{\U^\ibf}$ is constant, it follows that $\St$ is a.s.\ uniquely determined.

(ii)$\volgt$(v): Since $G_{\U^\ibf}$ is constant, we can define
\be\label{Xif}
X_\ibf:=G_{\U^\ibf}\big((x_{\ibf\jbf})_{\jbf\in\nab\U^\ibf}\big),
\ee
where the right-hand side does not depend on the choice of
$(x_{\ibf\jbf})_{\jbf\in\nab\U^\ibf}$. It is straightforward to check that
$(\omb_\ibf,X_\ibf)_{\ibf\in\ovT}$ satisfies conditions (i)--(iii) of
Lemma~\ref{L:RTP} and hence is an RTP corresponding to $\ga$. It follows that
$\nu:=\P[X_\wurz\in\,\cdot\,]$ solves the RDE (\ref{RDE}).

Let $(Y_\ibf)_{\ibf\in\ovT}$ be an independent i.i.d.\ collection of
$S$-valued random variables with common law $\mu_0$, let
$(\sig_\ibf)_{\ibf\in\ovT}$ be an independent i.i.d.\ collection of
exponential random variables with mean $|\rb|^{-1}$, and define
$X^t_\wurz:=G_t\big((Y_\ibf)_{\ibf\in\nab\St_t}\big)$. Then $X^t_\wurz$ has law
$\mu_t$ by Theorem~\ref{T:meanrep}. Since $G_t=G_{\St_t}$ with $\St_t\up\St$, 
we see from (\ref{Xif}) that $\P[X^t_\wurz\neq X_\wurz]\to 0$ as $t\to\infty$,
proving that $\|\mu_t-\nu\|\to 0$.

(iii)$\volgt$(iv): We note the following general principle: if $S_1,S_2,S_3$
are Polish spaces and $(X_1,X_2)$ and $(X'_1,X_3)$ are random variables taking
values in $S_1\times S_2$ resp.\ $S_1\times S_3$ such that $X_1$ and $X'_1$
are equal in law, then we can couple $(X_1,X_2)$ and $(X'_1,X_3)$ such that
$X_1=X'_1$. To see this, let $\mu$ denote the law of $X_1$, let $K_i(x_1,\di
x_i)$ denote a regular version of the conditional law of $X_i$ given $X_1$
resp.\ $X'_1$ $(i=1,2)$, and define the joint law of $X_1,X_2,X_3$ as
\be
\P[X_1\in A_1,\ X_2\in A_2,\ X_3\in A_3]
:=\int_{A_1}\mu(\di x_1)\int_{A_2}K_2(x_1,\di x_2)\int_{A_3}K_{3}(x_1,\di x_3),
\ee
i.e., make $X_2$ and $X_3$ conditionally independent given $X_1$. Applying
this general principle, we see that if $\nu_1,\nu_2$ are solutions to the RDE
(\ref{RDE}), then we can couple the associated RTPs
$(\omb_\ibf,X^1_\ibf)_{\ibf\in\ovT}$ and $(\omb'_\ibf,X^2_\ibf)_{\ibf\in\ovT}$
in such a way that $\omb_\ibf=\omb'_\ibf$ for all $\ibf\in\ovT$.
Since $\St$ is a.s.\ uniquely determined, it follows that
$X^1_\wurz=X^2_\wurz$ a.s.\ and hence $\nu_1=\nu_2$.
The same argument also shows that any solution to the bivariate
RDE is concentrated on the diagonal, which by Theorem~\ref{T:bivar} implies
endogeny.

(iii) and $S$ finite imply (ii): Since $\V\sub\U$ and $\V$ root determining
imply that $\U$ is root determining, we see that $\P[G_t\mbox{ not constant}]$
decreases to $\P[G_t\mbox{ not constant }\forall t\geq 0]$. Assume that this
event has positive probability and condition on it. Choose
$t(n)\to\infty$. Then there exist $x^n,y^n\in\Xi_{\St_{t(n)}}$ such that
$x^n_\wurz\neq y^n_\wurz$. Since $S$ is finite, the sequences $x^n$ and $y^n$
have subsequences that converge pointwise for each $\ibf\in\St$ to limits
$x^\infty,y^\infty$. It is easy to see that
$x^\infty,y^\infty\in\Xi_\St$. Moreover, $x^\infty_\wurz\neq y^\infty_\wurz$.
This shows that on the event $\{G_t\mbox{ not constant }\forall t\geq 0\}$,
the tree $\St$ is not uniquely determined.

(ii) and $S=\{0,1\}$ imply (i): It suffices to show that each root determining
subtree $\U$ of $\St$ contains a uniquely determined subtree. For any
$\ibf\in\ovT$, let $\ovT^{(\ibf)}:=\{\ibf\jbf:\jbf\in\ovT\}$ denote $\ibf$
and its descendants, and let $\Xi_{\U,\ibf}$ denote the set of all
$(x_\jbf)_{\jbf\in(\U\cup\nab\U)\cap\ovT^{(\ibf)}}$ that satisfy
\be
x_\jbf=\ga[\omb_\jbf](x_{\jbf 1},\ldots,x_{\jbf\kappa(\omb_\jbf)})
\qquad(\jbf\in\U\cap\ovT^{(\ibf)}).
\ee
Define $(\chi_\ibf)_{\ibf\in\U\cup\nab\U}$ by
\be\label{chibf}
\chi_\ibf:=\{x_\ibf:x\in\Xi_{\U,\ibf}\}\qquad(\ibf\in\U\cup\nab\U).
\ee
We claim that
\be\ba{l}\label{chiextend}
\mbox{If $\V\sub\U$ is a subtree, then any $x\in\Xi_\V$ that satisfies}\\
\mbox{$x_\ibf\in\chi_\ibf$ for all $\ibf\in\nab\V$ can be extended to an
 $x\in\Xi_\U$.}
\ec
Indeed, this follows from the fact that the sets
$(\U\cup\nab\U)\cap\ovT^{(\ibf)}$ for different $\ibf\in\nab\V$ are mutually
disjoint, which allows us to choose $x\in\Xi_{\U,\ibf}$ independently for each
$\ibf\in\nab\V$.

Note that $\chi_\ibf=S$ for $\ibf\in\nab\U$ and
$|\chi_\wurz|=1$ if $\U$ is root determining. Let $\V$ be the
connected component of $\{\ibf\in\U:|\chi_\ibf|=1\big\}$
that contains $\wurz$. Since $S$ has only two elements, $\chi_\ibf=S$ for
$\ibf\in\nab\V$. Using (\ref{chiextend}), it follows that $\V$ is uniquely
determined.
\epro

For completeness, we give three examples to show that the implications
(ii)$\volgt$(i), (iii)$\volgt$(ii), and (v)$\volgt$(ii) do not hold in
general. In all of these examples, $\kappa(\om)=1$ for all $\om\in\Om$, which
means $\St=\{\wurz,1,11,111,\ldots\}=\{1_{(n)}:n\geq 0\}$, where $1_{(n)}$
denotes the word of length $n$ made from the alphabet $\{1\}$. It follows that
the operator $T$ from (\ref{Tdef}) is just the linear operator associated with
the transition kernel of a Markov chain. In all our examples, we take
$\ga[\om]=g$ for all $\om\in\Om$, where $g:S\to S$ is a fixed map.

\begin{example}[(ii)$\not\volgt$(i)]
Let $S=\{0,1,2\}$ and $g(x):=(x-1)\vee 0$. Then $\St$ a.s.\ contains a root
determining subtree but $\St$ a.s.\ does not contain a uniquely determined
subtree.
\end{example}
\bpro
The subtree $\U:=\{\wurz,1\}$ is root determining, since $G_\U(x)=(x-2)\vee
0=0$ for all $x\in S^{\nab\U}=S$. On the
other hand, if $\V=\{\wurz,1,11,\ldots,1_{(n)}\}$ is a finite subtree of $\St$
that contains the root, then there exist $x,y\in\Xi_\V$ with $x_{1_{(n)}}=0$
and $y_{1_{(n)}}=1$, which shows $\V$ is not uniquely determined.
\epro

\begin{example}[(iii)$\not\volgt$(ii)] Let $S=\N$ and
\be
g(x):=\left\{\ba{ll}
0\quad&\mbox{if }x=0,\\
x+1\quad&\mbox{if }x\geq 1.
\ea\right.
\ee
Then $\St$ is a.s.\ uniquely determined but $\St$ a.s.\ contains no root
determining subtree.
\end{example}
\bpro
If $x\in\Xi_\St$ satisfies $x_{1_{(n)}}=m\neq 0$ for some $n$, then
$x_{1_{(n+k)}}\neq 0$ and $x_{1_{(n+k)}}=m-k$ for all $k\geq 0$, which leads
to a contradiction. It follows that $\Xi_\St$ contains a single element, which
is given by $x_{1_{(n)}}=0$ for all $n\geq 0$. In particular, $\St$ is uniquely
determined. On the other hand, for each finite subtree
$\U=\{\wurz,1,11,\ldots,1_{(n-1)}\}$ that contains the root, the function
$G_\U$ is of the form $G_\U(0)=0$ and $G_\U(x)=x+n$ $(x\geq 1)$, which is
clearly not constant.
\epro

\begin{example}[(v)$\not\volgt$(ii)]
Let $S=\{0,1\}$ and $g(x):=1-x$ with $\pi(\{g\})>0$. Then the RDE
(\ref{RDE}) has a solution $\nu$ that is globally attractive in the sense
that any solution $(\mu_t)_{t\geq 0}$ to (\ref{mean}) satisfies
$\|\mu_t-\nu\|\asto{t}0$, where $\|\,\cdot\,\|$ denotes the total variation
norm. Nevertheless, $\St$ contains no root determining subtree.
\end{example}
\bpro
Since the continuous-time Markov chain that jumps from $x$ to $1-x$ 
with rate $\pi(\{g\})$ is ergodic, the RDE (\ref{RDE}) has a solution
$\nu$ that is globally attractive. On the other hand, if
$\U=\{\wurz,1,11,\ldots,1_{(n-1)}\}$ is a finite subtree that contains the
root, then $G_\U(x)=x$ if $n$ is even and $G_\U(x)=1-x$ if $n$ is odd, so
$G_\U$ is not constant.
\epro

\bpro[of Lemma~\ref{L:monuni}]
By Lemma~\ref{L:detree}, it suffices to prove that if the RDE (\ref{RDE}) has
a unique solution, then $G_t$ is constant for $t$ large enough. By
Proposition~\ref{P:lowup}, the RDE (\ref{RDE}) has a unique solution if and
only if $\nu_{\rm low}=\nu_{\rm upp}$. Let $0$ and $1$ denote the minimal and
maximal elements of $S$. By Lemma~\ref{L:extrmap}, $G_t(0,\ldots,0)$
and $G_t(1,\ldots,1)$ converge as $t\to\infty$ to a.s.\ limits with laws
$\nu_{\rm low}$ and $\nu_{\rm upp}$, respectively. Since $\ga[\om]$ is
monotone for each $\om$, the maps $G_t$ are monotone, and hence
\be\label{squeeze}
G_t(0,\ldots,0)\leq G_t(x)\leq G_t(1,\ldots,1)
\ee
for all $x\in S^{\nab \St_t}$. Since $S$ is finite, if the laws of the left- and
right-hand sides of (\ref{squeeze}) converge to the same limit, then
$\lim_{t\to\infty}\P[G_t(0,\ldots,0)=G_t(1,\ldots,1)]=1$, proving that $G_t$
is constant for $t$ large enough.
\epro

\subsection{Duality}\label{S:dual}

In this subsection, we prove Lemma~\ref{L:good}. For a start, we
will generalize quite a bit and assume that $S$ is a finite partially ordered
set and that $\ga[\om]:S^{\kappa(\om)}\to S$ is monotone for all $\om\in\Om$,
where $S^{\kappa(\om)}$ is equipped with the product partial order. As in
Subsection~\ref{S:RTPconstr}, we let $\Ti$ denote the set of all finite
subtrees $\U\sub\ovT$ such that either $\wurz\in\U$ or $\U=\emptyset$. For
each $\U\in\Ti$, we define $G_\U:S^{\nab\U}\to S$ as in (\ref{GUdef}), where
$\nab\U:=\{\wurz\}$ if $\U=\emptyset$.

For any $\U\in\Ti$, we let $\Sig_\U$ denote the set of all
$(y_\ibf)_{\ibf\in\U\cup\nab\U}$ that satisfy
\be\label{Sigdef}
(y_{\ibf 1},\ldots,y_{\ibf\kappa(\omb_\ibf)})
\mbox{ is a minimal element of }
\big\{x\in S^{\kappa(\omb_\ibf)}:\ga[\omb_\ibf](x)\geq y_\ibf\big\}
\qquad(\ibf\in\U).
\ee

\bl[Monotone duality]
For\label{L:mondual} any $\U\in\Ti$, $x\in S^{\nab\U}$, and $z\in S$, one has
$G_\U(x)\geq z$ if and only if there exists a $y\in\Sig_\U$ such that
$y_\wurz=z$ and $x\geq y$ on $\nab\U$.
\el
\bpro
Fix $z\in S$. For each $\U\in\Ti$, let us write
\be\label{YU}
Y_\U:=\big\{(y_\ibf)_{\ibf\in\nab\U}:y\in\Sig_\U,\ y_\wurz=z\big\}.
\ee
Then we need to show that
\be
G_\U(x)\geq z\quad\desd\quad
\exists y\in Y_\U\mbox{ s.t.\ }x\geq y.
\ee
The proof is by induction on the number of elements of $\U$. If
$\U=\emptyset$, then $G_\U$ is the identity map, $Y_\U=\{z\}$, and the
statement is trivial.

We will show that if the statement is true for $\U$ and if $\jbf \in\nab\U$,
then the statement is also true for $\V:=\U\cup\{\jbf\}$. Let $x\in S^{\nab\V}$
and inductively define $x_\ibf$ for $\ibf\in\V$ as in (\ref{xinduc}).
By the induction hypothesis, $x_\wurz\geq z$ if and only if 
\be\label{Udual}
(x_\ibf)_{\ibf\in\nab\U}\geq y\mbox{ for some }y\in Y_\U.
\ee
Here $\nab\V=(\nab\U\beh\{\jbf\})\cup\{\jbf 1,\ldots,\jbf\kappa(\omb_{\jbf})\}$
and
\be\ba{r@{\,}l}\label{dualevol}
\dis Y_\V=\Big\{&\dis y\in S^{\nab\V}:\exists y'\in \Sig^{\nab\U}\mbox{ s.t.\ }
y_{\ibf}=y'_{\ibf}\ \forall \ibf\in\nab\U\beh\{\jbf\}\\
\dis&\mbox{and}\quad
(y_{\jbf 1},\ldots,y_{\jbf\kappa(\omb_{\jbf})})
\mbox{ is a minimal element of }
\big\{x\in S^{\kappa(\omb_{\jbf})}:\ga[\omb_{\jbf}](x)\geq y_{\jbf}\big\}
\Big\}.
\ec
It follows that (\ref{Udual}) is equivalent to
\be
(x_\ibf)_{\ibf\in\nab\V}\geq y\mbox{ for some }y\in Y_\V,
\ee
which completes the induction step of the proof.
\epro

\bl[Minimal elements]
Assume\label{L:minprop} that for all $\om\in\Om$, there do not exist
$z,z'\in S$ and minimal elements $y,y'$ of $\{y:\ga[\om](y)\geq z\}$
resp.\ $\{y':\ga[\om](y')\geq z'\}$ such that $z\not\leq z'$
but $y\leq y'$. Fix $z\in S$. For any $\U\in\Ti$, define $Y_\U$ as in
(\ref{YU}) dependent on $z$. Then
\be
Y_\U=\big\{y \in S^{\nab\U}:y\mbox{ is a minimal element of }
G_\U^{-1}(\{z\})\big\}.
\ee
\el
\bpro
By Lemma~\ref{L:mondual},
\be
G_\U^{-1}\big(\{z\}\big)
=\big\{x\in S^{\nab\U}:x\geq y\mbox{ for some }y\in Y_\U\big\}.
\ee
In view of this, it suffices to prove that
\be\label{hypot}
\mbox{$Y_\U$ does not contain two elements $y,y'$ with $y\neq y'$ and $y\leq
  y'$.}
\ee
The proof is by induction on the number of elements of $\U$. If
$\U=\emptyset$, then $\nab\U=\{\wurz\}$ and $Y_\U$ consists of a single element
that has $y_\wurz=z$, so (\ref{hypot}) is satisfied. Assume that (\ref{hypot})
holds for $\U$ and let $\V:=\U\cup\{\ibf\}$ for some $\ibf\in\nab\U$.
Then (\ref{dualevol}) and the assumption of the lemma imply that
(\ref{hypot}) holds for $\V$.
\epro

\bl[Sets with two elements]
Assume\label{L:01} that $S=\{0,1\}$ and that $\ga[\om](0,\ldots,0)=0$ for all
$\om\in\Om$. Then the assumption of Lemma~\ref{L:minprop} is satisfied.
\el
\bpro
If $z\not\leq z'$ then we must have $z=1$ and $z'=0$, so we must show that
there do not exist minimal elements $y,y'$ of $\{y:\ga[\om](y)\geq 1\}$
resp.\ $\{y':\ga[\om](y')\geq 0\}$ such that $y\leq y'$. Clearly,
$\{y':\ga[\om](y')\geq 0\}=\{0,1\}^{\kappa(\om)}$ has only one minimal element,
which is the configuration $(0,\ldots,0)\in\{0,1\}^{\kappa(\om)}$,
so we must show that there does not exist a minimal element $y$ of
$\{y:\ga[\om](y)\geq 1\}$ such that $y\leq(0,\ldots,0)$. Equivalently, this
says that $\ga[\om](0,\ldots,0)\not\geq 1$ which is satisfied since
$\ga[\om](0,\ldots,0)=0$.
\epro

\bl[Lower and upper solutions]
Assume\label{L:lowdu} that $S$ is a finite partially ordered set that contains
minimal and maximal elements, denoted by 0 and 1. Assume that $\ga[\om]$ is
monotone for all $\om\in\Om$. Then, for all $z\in S$,
\bc
\dis\nu_{\rm upp}\big(\{x:x\geq z\}\big)
&=&\dis\P\big[\exists y\in\Sig_\St\mbox{ s.t.\ }y_\wurz=z\big]\\[5pt]
\dis\nu_{\rm low}\big(\{x:x\geq z\}\big)&=&\dis
\P\big[\exists y\in\Sig_\St\mbox{ s.t.\ }y_\wurz=z\mbox{ and }
\{\ibf\in\St:y_\ibf\neq 0\}\mbox{ is finite}\big].
\ec
\el
\bpro
By Lemma~\ref{L:mondual}, $G_\U(1,\ldots,1)\geq z$ if and only if
$\Sig^z_\U:=\{y\in\Sig_\U:y_\wurz=z\}$ is not empty. If $\U\sub\V$, then
$\Sig^z_\V\neq\emptyset$ implies $\Sig^z_\U\neq\emptyset$, so
the events $\{\Sig^z_{\U^{(n)}}\neq\emptyset\}$ decrease to a limit.
We claim that this is the event $\{\Sig^z_\St\neq\emptyset\}$.
Since the restriction of an element $y\in\Sig^z_\St$ to $\U$ yields an element
of $\Sig^z_\U$, it is clear that
\be
\big\{\Sig^z_{\U^{(n)}}\neq\emptyset\ \forall n\big\}
\supset\{\Sig^z_\St\neq\emptyset\}.
\ee
Conversely, if for each $n$ there exists some $y(n)\in\Sig^z_{\U^{(n)}}$, then
by the finiteness of $S$ we can select a subsequence of the $y(n)$ that
converges pointwise to a limit $y$. Since $y\in\Sig^z_\St$, this proves the
other inclusion. By Lemma~\ref{L:extrmap}, it follows that
\be
\P[X^{\rm upp}_\wurz\geq z]=\P\big[\Sig^z_\St\neq\emptyset\big].
\ee

By Lemma~\ref{L:mondual}, $G_\U(0,\ldots,0)\geq z$ if and only if
$\Sig^z_\U$ contains an element $y$ such that $y_\ibf=0$ for all
$\ibf\in\nab\U$. Since for each $\om\in\Om$, the zero configuration
$(0,\ldots,0)$ is the unique minimal element of $\{x\in
S^{\kappa(\om)}:\ga[\om](x)\geq 0\}$, we observe that if
$y\in\Sig_\U$ satisfies $y_\ibf=0$ for all $\ibf\in\nab\U$, then 
$y$ can uniquely by extended to an element of $\Sig_\V$ for any $\V\supset\U$
by putting $y_\ibf:=0$ for $\ibf\in(\V\cup\nab\V)\beh(\U\cup\nab\U)$.
In view of this, by Lemma~\ref{L:extrmap},
\be
\P[X^{\rm low}_\wurz\geq z]=\P\big[\exists y\in\Sig^{z}_\St\mbox{ s.t.\ }
\{\ibf\in\St:y_\ibf\neq 0\}\mbox{ is finite}\big].
\ee
\epro

\bpro[of Lemma~\ref{L:good}]
Recall the definition of $\Sig_\U$ in (\ref{Sigdef}).
We observe that $\Op$ is an open subtree of $\U$ if and
only if its indicator function $1_\Op$ satisfies $1_\Op\in\Sig_\U$ and
$1_\Op(\wurz)=1$. In view of this, (\ref{dualow}) is just a special case of
Lemma~\ref{L:lowdu}. Formula (\ref{YGS}) follows from Lemmas~\ref{L:minprop}
and \ref{L:01} applied to $\U=\St_t$.
\epro

\section{Cooperative branching}\label{S:coop}

In this section we prove all results that deal specifically with our running
example of a system with cooperative branching and deaths. In
Subsection~\ref{S:cobiv}, we prove Proposition~\ref{P:copbiv} about the
bivariate mean-field equation. In Subsection~\ref{S:cohigh}, we prove
Theorem~\ref{T:coblev} and Lemma~\ref{L:nontriv} about the higher-level
mean-field equation. In Subsection~\ref{S:cotree}, finally, we prove Lemmas~\ref{L:cobU}
and \ref{L:birth} which illustrate the concepts of minimal root determining
subtrees and open subtrees in the concrete set-up of our example.

\subsection{The bivariate mean-field equation}\label{S:cobiv}

In this subsection we prove Proposition~\ref{P:copbiv}. 
We identify a measure $\mu^{(2)}$ on $\{0,1\}^2$ with the function
$\mu^{(2)}:\{0,1\}^2\to\R$ defined as $\mu^{(2)}(0,0):=\mu(\{(0,0)\})$,
$\mu^{(2)}(0,1):=\mu(\{(0,1)\})$, etc. We parametrize a measure
$\mu^{(2)}\in\Pc_{\rm sym}(\{0,1\}^2)$ by the parameters
\be\label{prdef}
p:=\mu^{(2)}(1,0)+\mu^{(2)}(1,1)
\quand
r:=\mu^{(2)}(0,1)+\mu^{(2)}(1,0)+\mu^{(2)}(1,1).
\ee
We observe that $\mu^{(2)}(0,0)=1-r$, $\mu^{(2)}(1,0)=\mu^{(2)}(0,1)=r-p$,
and hence $\mu^{(2)}(1,1)=1-(1-r)-2(r-p)=2p-r$. It follows that $p$ and $r$
determine $\mu^{(2)}$ uniquely and indeed, the map
\be\label{prspace}
\Pc_{\rm sym}(\{0,1\}^2)\ni\mu^{(2)}\mapsto(p,r)\in
\Di:=\big\{(p,r)\in\R^2:0\leq p\leq 1,\ p\leq r\leq 1\wedge 2p\big\}
\ee
is a bijection. Moreover, $\mu^{(2)}$ is concentrated on the diagonal
$\{(0,0),(1,1)\}$ if and only if $p=r$. A function $(\mu^{(2)}_t)_{t\geq 0}$
with values in $\Pc_{\rm sym}(\{0,1\}^2)$ gives through (\ref{prdef}) rise to
a function $(p_t,r_t)_{t\geq 0}$ taking values in $\Di$.

\bl[Change of parameters]
A\label{L:mupr} function $(\mu^{(2)}_t)_{t\geq 0}$ with values in $\Pc_{\rm
  sym}(\{0,1\}^2)$ solves (\ref{pibivar}) if and only if the associated
function $(p_t,r_t)_{t\geq 0}$ solves
\be\left.\ba{r@{\quad}r@{\,}c@{\,}l}\label{prODE}
{\rm(i)}&\dis\dif{t}p_t&=&\dis\al p_t^2(1-p_t)-p_t,\\[5pt]
{\rm(ii)}&\dis\dif{t}r_t&=&\dis\al\big[r_t^2-2(r_t-p_t)^2\big](1-r_t)-r_t,
\ea\quad\right\}\ (t\geq 0).
\ee
\el
\bpro
As noted in Section~\ref{S:end}, if $\mu^{(2)}_t$ solves the bivariate
mean-field equation, then its one-dimensional marginals solve the mean-field equation
(\ref{mean}). Since $\mu^{(2)}_t$ is symmetric, both
marginals are the same. We denote these by $\mu_t$. Then
$p_t:=\mu_t(\{1\})$ and the equation we find for $p_t$ is the same as in
(\ref{coopmean}).

We will now obtain the equation for the parameter $r_t$. 
By definitions \eqref{Tgdef} and \eqref{cobdth} we have for any $\mu^{(2)} \in
\Pc(\{0,1\}^2)$ that $T_{\cob^{(2)}}(\mu^{(2)})$ is the law of the random
variable
\be
\big(X^1_1\vee(X^1_2\wedge X^1_3),X^2_1\vee(X^2_2\wedge X^2_3)\big),
\ee
where $(X^1_i,X^2_i)$ $(i=1,2,3)$ are i.i.d. with law $\mu^{(2)}$.
It follows that
\be\ba{l}
\dis T_{\cob^{(2)}}(\mu^{(2)})(0,0)
=\P\big[(X^1_1,X^2_1)=(0,0)\big]
\Big(1-\P\big[(X^1_2,X^2_2)\neq(0,0)\big]
\P\big[(X^1_3,X^2_3)\neq(0,0)\big]\\[5pt]
\dis\quad
+\P\big[(X^1_2,X^2_2)=(0,1)]\P\big[(X^1_{3},X^2_{3})=(1,0)]
+\P\big[(X^1_2,X^2_2)=(1,0)]\P\big[(X^1_{3},X^2_{3})=(0,1)]\Big)\\[5pt]
\dis\quad=(1-r)\big(1-r^2+2(r-p)^2\big).
\ec
Similar, but simpler considerations give
\be
T_{\dth^{(2)}}(\mu^{(2)})(0,0)=1.
\ee
Equation (\ref{pibivar}) in the point $(0,0)$ now gives
\be
\dif{t}(1-r_t)=\dif{t}\mu^{(2)}_t(0,0)
=\al\big\{(1-r_t)\big(1-r_t^2+2(r_t-p_t)^2\big)-(1-r_t)\big\}
+\big\{1-(1-r_t)\big\},
\ee
which simplifies to the second equation in (\ref{prODE}).
\epro

In view of Lemma~\ref{L:mupr} and the remarks that precede it,
Proposition~\ref{P:copbiv} follows from the following proposition.

\bp[Bivariate differential equation]
For $\al>4$, the equation (\ref{prODE}) has four fixed points in
the space $\Di$ defined in \eqref{prspace}, which are of the form
\be\label{bivfix}
(z_{\rm low},z_{\rm low}),\quad(z_{\rm mid},z_{\rm mid}),
\quad(z_{\rm mid},r_{\rm mid}),\quand(z_{\rm upp},z_{\rm upp}),
\ee
with $z_{\rm low},z_{\rm mid},z_{\rm upp}$ as in (\ref{fixpts}) and
$z_{\rm mid}<r_{\rm mid}$. Solutions to (\ref{prODE}) started in $\Di$
converge to one of these fixed points, the domains of attraction being
\be\label{prdomains}
\{(p,r):p<z_{\rm mid}\},\quad\{(z_{\rm mid},z_{\rm mid})\},\quad
\{(z_{\rm mid},r):r>z_{\rm mid}\},\quand\{(p,r):p>z_{\rm mid}\},
\ee
respectively. For $\al=4$, the equation (\ref{prODE}) has two fixed points in
the space $\Di$, which are
\be
(z_{\rm low},z_{\rm low})\quand(z_{\rm mid},z_{\rm mid})=(z_{\rm upp},z_{\rm upp}),
\ee
with domains of attraction
\be
\{(p,r):p<z_{\rm mid}\}\quand\{(p,r):p\geq z_{\rm mid}\}.
\ee
For $\al<4$, $(z_{\rm low},z_{\rm low})$ is the only fixed point in $\Di$ and its
domain of attraction is the whole space $\Di$.
\ep
\bpro
In Section~\ref{S:flow} we have found all fixed points of (\ref{prODE})~(i)
and determined their domains of attraction. It is clear from (\ref{prODE}) that if
$z$ is a fixed point of (\ref{prODE})~(i), then $(z,z)$ is a fixed point of
(\ref{prODE}), so $(z_{\rm low},z_{\rm low})$ and for $\al\geq 4$ also
$(z_{\rm mid},z_{\rm mid})$ and $(z_{\rm upp},z_{\rm upp})$ are fixed points
of (\ref{prODE}).

If $\al\geq 4$ and $p_0<z_{\rm mid}$ or if $\al<4$ and $p_0$ is arbitrary,
then we have seen in Section~\ref{S:flow} that solutions to (\ref{prODE})~(i) satisfy
$p_t\to 0=z_{\rm low}$ as $t\to\infty$. Since $0\leq r_t\leq 2p_t$, it follows
that also $r_t\to 0$. This proves the statements of the proposition about the
domain of attraction of $(z_{\rm low},z_{\rm low})$ for all values of $\al$.

Let
\be\label{PRdef}
P_\al(p):=\al p^2(1-p)-p
\quand
R_{\al,p}(r):=\al\big[r^2-2(r-p)^2\big](1-r)-r
\ee
denote the drift functions of $p_t$ and $r_t$, respectively. We observe that
$R_{\al,p}(r)\leq P_\al(r)$ $(p,r\in\R)$ and $P_\al(r)<0$ for all $z_{\rm upp}<r\leq 1$,
which implies that
\be\label{Runi}
\sup_{p\in\R}R_{\al,p}(r)<0\qquad(z_{\rm upp}<r\leq 1).
\ee
It follows that solutions of (\ref{prODE}) satisfy
\be\label{rsup}
\limsup_{t\to\infty}r_t\leq z_{\rm upp}.
\ee
If $\al>4$ and $p_0>z_{\rm mid}$ or if $\al=4$ and $p_0\geq z_{\rm mid}$, we
have seen in Section~\ref{S:flow} that solutions to (\ref{prODE})~(i) satisfy
$p_t\to z_{\rm upp}$ as $t\to\infty$. Combining this with (\ref{rsup}) and
the fact that $p_t\leq r_t$, we see that $(p_t,r_t)\to(z_{\rm upp},z_{\rm upp})$.

To complete the proof, we must investigate the long-time behavior of solutions
of (\ref{prODE}) when $\al>4$ and $p_0=z_{\rm mid}$. In this case $p_t=z_{\rm
  mid}$ for all $t\geq 0$ and $r_t$ takes values in $[z_{\rm mid},2z_{\rm mid}]$
and solves the differential equation
\be\label{rR}
\dif{t}r_t=R_{\al,z_{\rm mid}}(r_t)\qquad(t\geq 0).
\ee
It is clear $r_t=z_{\rm mid}$ for all $t\geq 0$ is a solution. Since
  $z_{\rm mid}<1/2$, in view of (\ref{prspace}), we must prove
that all solutions with $z_{\rm mid}<r_0\leq 2z_{\rm mid}$ converge to a
nontrivial fixed point. We write
\be
R_{\al,z_{\rm mid}}(r)=P_\al(r)-2\al(r-z_{\rm mid})^2(1-r).
\ee
Since the first term has a positive slope at $r=z_{\rm mid}$ while the second
term has zero slope, we conclude that $R_{\al,z_{\rm mid}}$ has a positive
slope at $r=z_{\rm mid}$. Since solutions to (\ref{prODE}) do not leave the
domain $\Di$, we must have $R_{\al,z_{\rm mid}}(2z_{\rm mid})\leq 0$.
Since $R_{\al,z_{\rm mid}}(r)=\al r^3+O(r^2)$ as $r\to\infty$, we must have
$R_{\al,z_{\rm mid}}(r)>0$ for $r$ sufficiently large. These observations
imply that the cubic function $R_{\al,z_{\rm mid}}$ has three zeros
$r_{\rm low}<r_{\rm mid}<r_{\rm up}$ with
\be
r_{\rm low}=z_{\rm mid}<r_{\rm mid}\leq 2z_{\rm mid}<r_{\rm upp},
\ee
and $R_{\al,z_{\rm mid}}>0$ on $(z_{\rm mid},r_{\rm mid})$ and $R_{\al,z_{\rm mid}}<0$ on
$(r_{\rm mid},r_{\rm upp})$. It follows that solutions to (\ref{rR}) started
with $z_{\rm mid}<r_0\leq 2z_{\rm mid}$ satisfy $r_t\to r_{\rm mid}$ as
$t\to\infty$.
\epro

\subsection{The higher-level mean-field equation}\label{S:cohigh}

In this subsection we prove Theorem~\ref{T:coblev} and
Lemma~\ref{L:nontriv}. We start with two preparatory lemmas.

\bl[Convex order and second moments]
Let\label{L:cosec} $S$ be a Polish space $S$ and let
$\rho_1,\rho_2\in\Pc(\Pc(S))$ satisfy $\rho_1\leq_{\rm cv}\rho_2$ and
$\rho^{(2)}_1=\rho^{(2)}_2$. Then $\rho_1=\rho_2$.
\el
\bpro
This follows from \cite[Lemma~14]{MSS18}.
\epro

In the next lemma, we use the notation $\ov\mu:=\P[\de_X\in\,\cdot\,]$
  defined in Subsection~\ref{S:high}.

\bl[Maximal measure in convex order]
Let\label{L:comax} $S$ be a Polish space and let $\mu\in\Pc(S)$. Then a
measure $\rho\in\Pc(\Pc(S))$ satisfies $\rho=\ov\mu$ if and only if
$\rho^{(2)}=\ov\mu^{(2)}$.
\el
\bpro
The condition $\rho^{(2)}=\ov\mu^{(2)}$ implies that the first moment measure
of $\rho$ is $\mu$. By (\ref{cvextr}), it follows that $\rho\leq_{\rm
  cv}\ov\mu$, so the statement follows from Lemma~\ref{L:cosec}.
\epro

\bpro[of Theorem~\ref{T:coblev}]
It follows from their definition that the measures
$\ov\nu_{\rm low},\ov\nu_{\rm mid},\ov\nu_{\rm upp}$
and $\un\nu_{\rm low},\un\nu_{\rm mid},\un\nu_{\rm upp}$
solve the higher-level RDE and their first moment measures are
$\nu_{\rm low},\nu_{\rm mid},\nu_{\rm upp}$, respectively.

By \cite[Thm~5]{MSS18}, one has $\un\nu=\ov\nu$ if and only if the RTP
corresponding to $\nu$ is endogenous. By Theorem~\ref{T:bivar}, endogeny is
equivalent to bivariate uniqueness, so we obtain from
Proposition~\ref{P:copbiv} that $\un\nu_{\rm low}=\ov\nu_{\rm low}$,
$\un\nu_{\rm mid}\neq\ov\nu_{\rm mid}$, and $\un\nu_{\rm upp}=\ov\nu_{\rm
  upp}$.

Since the second moment measures of
$\ov\nu_{\rm low},\ov\nu_{\rm mid},\ov\nu_{\rm upp}$ are of the form
(\ref{nmu}), we see that the measures
$\ov\nu^{(2)}_{\rm low},\ov\nu^{(2)}_{\rm mid},\ov\nu^{(2)}_{\rm upp}$ from
Proposition~\ref{P:copbiv} are indeed the second moment measures of
$\ov\nu_{\rm low},\ov\nu_{\rm mid},\ov\nu_{\rm upp}$.

By Proposition~\ref{P:levmom}, the second moment measure of $\un\nu_{\rm mid}$
solves the bivariate RDE. Since $\un\nu_{\rm mid}\neq\ov\nu_{\rm mid}$,
Lemma~\ref{L:comax} tells us that the second moment measure of $\un\nu_{\rm
  mid}$ is different from $\ov\nu^{(2)}_{\rm mid}$. It follows that the
measure $\un\nu^{(2)}_{\rm mid}$ from Proposition~\ref{P:copbiv} is indeed the
second moment measure of $\un\nu_{\rm mid}$.

Let $(\rho_t)_{t\geq 0}$ be a solution to the higher-level mean-field
equation. Assume that $\al>4$. Then Propositions~\ref{P:copbiv} and
\ref{P:levmom} tell us that $\rho^{(2)}_t$ converges to one of the fixed
points $\ov\nu^{(2)}_{\rm low},\un\nu^{(2)}_{\rm mid},\ov\nu^{(2)}_{\rm mid},
\ov\nu^{(2)}_{\rm upp}$, depending on whether
\be\ba{ll}
\dis{\rm(i)}\ \rho^{(1)}(\{1\})<z_{\rm mid},\qquad
&\dis{\rm(ii)}\ \rho^{(1)}(\{1\})=z_{\rm mid}
\mbox{ and }\rho^{(2)}\neq\ov\nu^{(2)}_{\rm mid},\\[5pt]
\dis{\rm(iii)}\ \rho^{(2)}=\ov\nu^{(2)}_{\rm mid},\qquad
&\dis{\rm(iv)}\ \rho^{(1)}(\{1\})>z_{\rm mid}.
\ec
By Lemma~\ref{L:comax}, these four cases correspond exactly to the four domains
of attraction in (\ref{domains}). To prove that in fact $\rho_t$ converges to
$\ov\nu_{\rm low},\un\nu_{\rm mid},\ov\nu_{\rm mid}$, or $\ov\nu_{\rm upp}$,
respectively, in each of these cases, by the compactness of
$\Pc(\Pc(\{0,1\}))$, it suffices to prove that if
$\rho_{t_n}\Rightarrow\rho_\ast$ along a sequence of times $t_n\to\infty$,
then $\rho_\ast$ is the right limit point.
In the cases~(i), (iii) and (iv) this is clear from Lemma~\ref{L:comax}.

To prove the statement also in case~(ii), let $(\rho'_t)_{t\geq 0}$ be the
solution to the higher-level mean-field equation started in
$\rho'_0=\de_{\nu_{\rm mid}}$. Then (\ref{cvextr}) and
Proposition~\ref{P:minconv} tells us that $\rho'_t\leq_{\rm cv}\rho_t$ for all
$t\geq 0$ and $\rho'_t\Rightarrow\un\nu_{\rm mid}$. Taking the limit
$t_n\to\infty$, using condition~(i) of Theorem~\ref{T:Stras}, we conclude that
$\un\nu_{\rm mid}\leq_{\rm cv}\rho_\ast$. Since moreover $\un\nu^{(2)}_{\rm
  mid}=\rho^{(2)}_\ast$, we can apply Lemma~\ref{L:cosec} to conclude that
$\un\nu_{\rm mid}=\rho_\ast$.

This completes the proof for $\al>4$. The cases $\al=4$ and $\al<4$ are
similar, but simpler.
\epro

\bpro[of Lemma~\ref{L:nontriv}]
We note that if $\eta_1,\eta_2,\eta_3\in[0,1]$, then
\bc
\dis 1_{\{\widehat{\cob}(\eta_1,\eta_2,\eta_3)=1\}}
=1_{\{\eta_1+(1-\eta_1)\eta_2\eta_3=1\}}
&=&\dis1_{\{\eta_1=1\}}\vee\big(1_{\{\eta_2=1\}}\wedge1_{\{\eta_3=1\}}\big),\\[5pt]
\dis 1_{\{\widehat{\cob}(\eta_1,\eta_2,\eta_3)>0\}}
=1_{\{\eta_1+(1-\eta_1)\eta_2\eta_3>0\}}
&=&\dis1_{\{\eta_1>0\}}\vee\big(1_{\{\eta_2>0\}}\wedge1_{\{\eta_3>0\}}\big).
\ec
Combining this with \eqref{cobdth} and Proposition \ref{P:levmom}, we see that
if $\rho$ solves the higher-level RDE \eqref{RDElevel}, then
\be
\int\rho(\di\eta)\eta,\quad\rho(\{1\}),\quand\rho((0,1])
\ee
must all solve the RDE \eqref{RDE}. Applying this to $\un\nu_{\rm mid}$ which has
$\int\un\nu_{\rm mid}(\di\eta)\eta=z_{\rm mid}$, we see that
\be
\un\nu_{\rm mid}(\{1\})\in\{0,z_{\rm mid},z_{\rm upp}\}\quand\un\nu_{\rm mid}((0,1])\in\{0,z_{\rm mid},z_{\rm upp}\}.
\ee
We first observe that since $\int\un\nu_{\rm mid}(\di\eta)\eta=z_{\rm mid}$,
we can have $\un\nu_{\rm mid}(\{1\})\geq z_{\rm mid}$ only if
$\un\nu_{\rm mid}=\ov\nu_{\rm mid}$, which we know is not the case, so we
conclude that $\un\nu_{\rm mid}(\{1\})=0$. If $\un\nu_{\rm mid}((0,1])\leq
z_{\rm mid}$ then $\int\un\nu_{\rm mid}(\di\eta)\eta=z_{\rm mid}$ forces
$\un\nu_{\rm mid}(\{1\})=z_{\rm mid}$, which we know is not the case, so we
conclude that $\un\nu_{\rm mid}((0,1])=z_{\rm upp}$ and hence
$\un\nu_{\rm mid}(\{0\})=1-z_{\rm upp}=z_{\rm mid}$, where the last equality
follows from (\ref{fixpts}).

To calculate $\int\!\un\nu_{\rm mid}(\di\eta)\,\eta^2$, we use that
$1-\un\nu^{(2)}_{\rm mid}(0,0)=r_{\rm mid}$, where $r_{\rm mid}$ is the
second largest solution $r$ of the equation $R_{\al,z_{\rm mid}}(r)=0$,
with $R_{\al,z_{\rm mid}}$ defined as in (\ref{PRdef}). The smallest
solution of the cubic equation $R_{\al,z_{\rm mid}}(r)=0$ is
$r=z_{\rm mid}$. Dividing by $(r-z_{\rm mid})$ yields a quadratic equation of
which $r_{\rm mid}$ is the smallest solution. Since these are straightforward,
but tedious calculations, we omit them.
\epro

\subsection{Root-determining and open subtrees}\label{S:cotree}

In this subsection we prove Lemmas~\ref{L:cobU} and \ref{L:birth}.\med

\bpro[of Lemma~\ref{L:cobU}]
Since $\kappa(\om)=3$ if $\ga[\om]=\cob$ and $\kappa(\om)=0$ if
$\ga[\om]=\dth$, we see that
\be
\int\rb(\di\om)(\kappa(\om)-1)=2\al-1,
\ee
which is $\leq 1$ if and only if $\al\leq\ha$. At the end of
Subsection~\ref{S:flow}, we have seen that in our example the RDE (\ref{RDE})
has a unique solution if and only if $\al<4$. By Lemma~\ref{L:monuni} this is
equivalent to condition~(ii) of Lemma~\ref{L:detree}. Since $S=\{0,1\}$,
Lemma~\ref{L:detree} tells is that in our example, conditions (i)--(iii)
are equivalent.

We claim that a finite subtree $\U\sub\St$ satisfying (\ref{minroot}) is
uniquely determined and in fact $x\in\Xi_\U$ implies $x_\ibf=0$ for all
$\ibf\in\U$. To prove this, let
$A=\{\ibf\in\U:x_\ibf=0\ \forall x \in\Xi_\U\}$. Since $\U$ is finite, if
$\U\beh A$ is not empty then we can find some $\ibf\in\U\beh A$ such that
$\ibf j\not\in A$ for $j=1,2,3$. (Here we take $\ovT$ to be the set of
all words made from the alphabet $\{1,2,3\}$.) If $\ga[\omb_\ibf]=\dth$, then
$x_\ibf=0$ for all $x\in\Xi_\U$ which contradicts the fact that
$\ibf\in\U\beh A$. But if $\ga[\omb_\ibf]=\cob$, then (\ref{minroot}) and the
fact that $\ibf j\not\in A$ for $j=1,2,3$ again imply
$x_\ibf=0$ for all $x\in\Xi_\U$, so we see that $\U\beh A$ must be empty.
In particular, this shows that $\U$ is root determining.

To see that $\U$ is a minimal root determining subtree, assume that $\V\sub\U$
is a smaller one. Then there must be be some $\ibf\in\V$ such that
$\ga[\omb_\ibf]=\cob$ and either $\ibf 1\not\in\V$ or $\V\cap\{\ibf 2,\ibf
3\}=\emptyset$. (Here we use that by definition, minimal root determining
subtrees contain the root, so $\V$ is not empty.) But then either $\ibf
1\in\nab\V$ or $\{\ibf 2,\ibf 3\}\sub\nab\V$. Define $x\in\Xi_\V$ inductively
by (\ref{xinduc}) with $x_{\jbf}=1$ for all $\jbf \in\nab\V$. Then $x_\ibf=1$.
Either $\ibf$ is the root or its predecessor $\lib$ satisfies
$x_{\,\lib}=1$ by (\ref{minroot}), so by induction we see that
$x_\ibf=1$. Since the all-zero configuration is also an element of $\Xi_\V$,
this proves that $\V$ is not root determining and hence $\U$ is minimal.
\epro

\bpro[of Lemma~\ref{L:birth}]
We observe that
\bc
\dis Y_\cob
&=&\dis\big\{y\in\{0,1\}^3:
y\mbox{ is a minimal element of }\cob^{-1}(\{1\})\big\}\\[5pt]
&=&\dis\big\{y\in\{0,1\}^3:
\cob(y)=1\mbox{ and }\cob(y')=0\ \forall y'\leq y,\ y'\neq y\big\}\\[5pt]
&=&\dis\big\{(1,0,0),(0,1,1)\big\}.
\ec
Since the set $\dth^{-1}(\{1\})=\{y\in\{0,1\}^0:\dth(y)=1\}$ is empty, it has
no minimal elements, and hence $Y_\dth=\emptyset$. On the other hand, 
$\bth(y)=1$ for all $y\in\{0,1\}^{0}$. In fact, $\{0,1\}^{0}=\{\wurz\}$ is a set with
only one element, the empty word, so $\bth^{-1}(\{1\})=\{\wurz\}$ and hence
the set of its minimal elements is $Y_{\bth}=\{\wurz\}$. Now (\ref{goodef}) with
the convention that $1_{A_\ibf}:=\wurz$ if $\kappa(\omb_\ibf)=0$ says
that $\Op$ is an open subtree of $\U$ if and only if:
\begin{enumerate}
\item $\{j\in\{1,2,3\}:\ibf j\in\Op\}=\{1\}$ or $=\{2,3\}$ for each
  $\ibf\in\Op\cap\U$ such that $\ga[\omb_\ibf]=\cob$,
\item $\wurz\in\emptyset$ for each $\ibf\in\Op\cap\U$ such that
  $\ga[\omb_\ibf]=\dth$,
\item $\wurz\in\{\wurz\}$ for each $\ibf\in\Op\cap\U$ such that
  $\ga[\omb_\ibf]=\bth$,
\end{enumerate}
which corresponds to the condition in (\ref{cobgood}).
\epro

\end{document}